\documentclass[leqno,12pt]{amsart}
\usepackage{amssymb, amsmath}
\usepackage{color}

\definecolor{darkgreen}{rgb}{0.0, 0.2, 0.13}
\definecolor{darkmagenta}{rgb}{0.55, 0.0, 0.55}

\newcommand{\eeq}{\end{equation}}
\newcommand{\ben}{\begin{eqnarray}}
\newcommand{\een}{\end{eqnarray}}
\newcommand{\beno}{\begin{eqnarray*}}
\newcommand{\eeno}{\end{eqnarray*}}
\newcommand{\norm}[1]{\left\Vert#1\right\Vert}

\def\div{\, \mbox{div}\,  }
\newcommand{\eps}{\varepsilon}

\newcommand{\cS}{\mathcal{S}}

\usepackage{amsmath,amssymb,amsfonts}
\usepackage{graphicx}

\usepackage[utf8]{inputenc}
\usepackage[colorlinks=true,linkcolor=black,citecolor=black]{hyperref}

\newcommand{\RR}{\mathbb{R}}

\newcommand{\ZZ}{\mathbb{Z}}
\newcommand{\NN}{\mathbb{N}}

\newcommand{\DD}{\Delta}

\newcommand{\Bcal}{\mathcal{B}}

\newcommand{\Fcal}{\mathcal{F}}

\newcommand{\Kcal}{\mathcal{K}}

\newcommand{\Rcal}{\mathcal{R}}
\newcommand{\Scal}{\mathcal{S}}
\newcommand{\Tcal}{\mathcal{T}}

\newcommand{\dd}{\partial}
\newcommand{\divv}{\mbox{div\,}}

\newcommand{\supp}{\mbox{supp }}

\newcommand{\abs}[1]{\left\vert#1\right\vert}
\newcommand{\set}[1]{\left\{#1\right\}}
\newcommand{\psca}[1]{\left\langle#1\right\rangle}

\newcommand{\pare}[1]{\left(#1\right)}
\newcommand{\f}{\frac}

\newcommand{\pa}{\partial}
\newcommand{\p}{\partial}
\newcommand\na{\nabla}

\newtheorem{theorem}{Theorem}[section]
\newtheorem{definition}[theorem]{Definition}
\newtheorem{lemma}[theorem]{Lemma}
\newtheorem{proposition}[theorem]{Proposition}

\newtheorem{remark}[theorem]{Remark}

\title{\textsc{Hydrostatic approximation of the 2D primitive equations in a thin strip}}

\author{Nacer Aarach}
\address{IMB, Universit\'e de Bordeaux, 351, cours de la Lib\'eration, 33405 Talence , France}
\email{nacer.aarach@math.u-bordeaux.fr}

\author{Van-Sang Ngo}
\address{Laboratoire de Math\'ematiques Rapha\"el Salem, UMR 6085 CNRS, Universit\'e de Rouen, 76801 Saint-Etienne du Rouvray Cedex, France}
\email{van-sang.ngo@univ-rouen.fr}

\date{}

\begin{document}

\maketitle

\begin{abstract}

We prove the global well-posedness of the 2D non-rotating primitive equations with no-slip boundary conditions in a thin strip of width $\eps$ for small data which are analytic in the tangential direction. We also prove that the hydrostatic limit (when $\eps \to 0$) is a couple of a Prandtl-like system for the velocity with a transport-diffusion equation for the temperature. 

\end{abstract}

\section{Introduction}

{\color{black}

\subsection{Primitive equations in a thin strip}

The primitive equations describe the motion of large scale fluids on the earth (hundreds to thousands of kilometers), typically an ocean or the atmosphere. The study of geophysical fluids of that scale involves two important phenomena: the effect of the vertical stratification due to the gravity and the effect of the rotation of the earth at large scale. The first effect naturally appears when we consider a fluid of variable density (hot and and cold air for instance). The fluid has then a vertical distribution where heavier layers stay under lighter ones and the gravity constantly maintains this structure against the internal movements of the fluid. 

The second effect becomes important at large scale when the fluid can ``sense'' the rotation of the earth. More precisely, two additional force terms will appear in the equations: the centrifugal force, which is included in the gravity gradient term, and the Coriolis force, which becomes large when the rotation is rapid (or when the scale is large) and which induces a vertical rigidity in the fluid movements. The latter is well known as the phenomenon of Taylor-Proudman columns. The estimate of the importance of this rigidity leads to the comparison between the typical time scale of the system and the Brunt-V\"ais\"al\"a frequency by means of the introduction of the Froude number $Fr$. For more details about physical considerations, we refer to \cite{CDGG4}, \cite{Cushman}, \cite{Pedlosky}, and \cite{Sadourny}.

In \cite{EmMa}, Embid and Majda considered specific scale simplifications by choosing the same scale for the rotation and the stratification. More precisely, the Froude number is supposed to be $Fr = \eps F$, where $F > 0$ is a constant (which will also be called ``Froude number'') and $\eps$ is the Rossby number. Using moment, energy, mass conservation laws, they obtained the following primitive equations in $\RR^3$
\begin{equation*}
	\left\{
	\begin{aligned}
		&\dd_t U_\eps + v_\eps\cdot\nabla U_\eps - LU_\eps + \frac{1}{\eps} \mathcal{A} U_\eps = \frac{1}{\eps} (-\nabla \Phi_\eps,0),\\
		&\divv v_\eps = 0,\\
		&U_\eps\vert_{t=0} = U_0,
	\end{aligned}
	\right.
\end{equation*} 
where the unknowns are $U_\eps = (v_\eps,\rho_\eps)$, which represents the velocity and the temperature of the fluid, and $\Phi_\eps$ the pressure. The operator $L$ is defined by $LU_\eps = (\nu \Delta v_\eps, \nu' \Delta \rho_\eps)$ and the skew-symmetric matrix $\mathcal{A}$ is given by
\begin{displaymath}
	\mathcal{A} = \begin{pmatrix}
		0 & -1 & 0 & 0 \\
		1 & 0 & 0 & 0 \\
		0 & 0 & 0 & F^{-1} \\
		0 & 0 & - F^{-1} & 0
	\end{pmatrix}.
\end{displaymath}

Another important remark concerning geophysical fluids is the difference between the horizontal scale (generally of order of hundreds to thousands of kilometers) and the vertical scale (generally a few kilometers for oceans and 10-20 kilometers for the atmosphere). In order to take into account this anisotropy, one can consider anisotropic viscosities as in the works of Chemin, Desjardins, Gallagher and Grenier (see \cite{CDGG4} and the references therein for instance) or in the work of Charve and Ngo (see \cite{CN2011} for instance). Another direction consists in studying the fluids in thin domains where a dimension is supposed to be very small and goes to zero, which is the main consideration of our work. 

In this paper, we will neglect the effect of the rotation and only focus on the effect of the vertical stratification. The combined effect of the rotation and the stratification in the full primitive equations will be studied in a forthcoming paper. In the two-dimensional case, these considerations lead to the following non-rotating primitive equations in the thin strip $\Scal^\eps = \set{(x,y) \in \RR^2 \;:\; 0 < y < \eps}$, where $\eps$ is supposed to be very small.
\begin{equation}
	\label{NS2}
	\left\{\;
	\begin{aligned}
		&\partial_t U^\eps + U^\eps .\nabla U^\eps-\eps^2 \Delta U^\eps + \nabla P^\eps = \frac{1}{\eps} (0,T^\eps) && \qquad \mbox{in }\, ]0,\infty[ \times \mathcal{S}^{\eps} \\
		&\partial_t T^\eps + U^\eps.\nabla T^\eps -  \Delta_{\eps} T^\eps = 0 && \qquad \mbox{in }\, ]0,\infty[ \times \mathcal{S}^{\eps}\\
		&\divv U^\eps = 0 && \qquad \mbox{in }\, ]0,\infty[ \times \mathcal{S}^{\eps},		
	\end{aligned}
	\right.
\end{equation}
Here, $U^\eps(t,x,y) = \pare{U_1^\eps(t,x,y), U_2^\eps(t,x,y)}$ denotes the velocity of the fluid and $P^\eps(t,x,y)$ the pressure which guarantees the divergence-free property of the velocity field $U^\eps$; $T^\eps(t,x,y)$ is the temperature of the fluid and the Froude number is supposed to be $F=1$. The system $\eqref{NS2}$ is complemented by the no-slip boundary condition
\begin{equation*}
	U^\eps_{\vert y=0} = U^\eps_{\vert y=\eps} = 0 \quad\mbox{and}\quad T^\eps_{\vert y=0} = T^\eps_{\vert y=\eps} = 0.
\end{equation*}
In the equation of the velocity, the Laplacian is $\Delta = \partial_x^2 + \partial_y^2$ and in the equation of the temperature, the anisotropic Laplacian $\Delta_{\eps} = \partial_{x}^2 + \eps^{2} \partial_{y}^2$ reflects the difference between the horizontal and the vertical scales. We will explain our choice of the diffusion term $\Delta_\eps$ in the next section.

\subsection{Hydrostatic limit of non-rotating primitive equations}

In this framework, it is believed that the fluid behaviors tend towards a geostrophic balance (see \cite{G1982}, \cite{H2004} or \cite{PZ2005}). In a formal way, as in \cite{PZZ2019}, taking into account this anisotropy, we consider the initial data in the following form, 
$$ U^\eps_{\vert t=0} = U_{0}^{\eps} = \left( u_0\pare{x,\frac{y}{\eps}},\eps v_{0}\pare{x,\frac{y}{\eps}}\right) \text{ \ \ in \ \ } \mathcal{S}^{\eps}$$ and $$T^\eps_{\vert t=0} = T^\eps_0 \pare{x,\frac{y}{\eps}}$$ 
and we look for solutions in the form
\begin{equation} 
	\label{S1eq3}
	\left\{\;
	\begin{aligned}
		&U^\eps(t,x,y)=\pare{u^\eps \pare{t,x,\frac{y}{\eps}}, \eps v^\eps \pare{t,x,\frac{y}{\eps}}}\\
		&T^\eps(t,x,y)=T^\eps \pare{t,x,\frac{y}{\eps}}\\
		&P^\eps(t,x,y)=p^\eps \pare{t,x,\frac{y}{\eps}}.
	\end{aligned}
	\right.
\end{equation}

In order to study the limit when $\eps \to 0$, we perform the rescaling $z = \frac{y}{\eps}$ and bring our problem to the fixed domain $\mathcal{S} = \left\{(x,z)\in\mathbb{R}^2:\ 0<z<1\right\}$. We rewrite the system $\eqref{NS2}$ as follows

\begin{equation}
	\label{eq:hydroPE}
	\left\{\;
	\begin{aligned}
		&\partial_t u^\eps +u^\eps\partial_x u^\eps + v^\eps\partial_zu^\eps-\eps^2\partial_x^2u^\eps-\partial_z^2u^\eps+\partial_x p^\eps=0 && \mbox{in } \; ]0,\infty[ \times \mathcal{S}\\
		&\eps^2\left(\partial_tv^\eps+u^\eps\partial_x v^\eps+v^\eps\partial_zv^\eps-\eps^2\partial_x^2v^\eps-\partial_z^2v^\eps\right)+\partial_zp^\eps=T^\eps && \mbox{in } \; ]0,\infty[ \times \mathcal{S}\\
		&\partial_t T^\eps+u^\eps\partial_x T^\eps+ v^\eps\partial_z T^\eps- \pa_x^2 T^\eps - \pa_z^2 T^\eps=0 && \mbox{in } \; ]0,\infty[ \times \mathcal{S}\\
		&\partial_x u^\eps+\partial_zv^\eps=0 && \mbox{in } \; ]0,\infty[ \times \mathcal{S}\\
		&\left(u^\eps, v^\eps, T^\eps \right)|_{t=0}=\left(u^\eps_0, v^\eps_0,T^\eps_0 \right) && \mbox{in } \; \mathcal{S}\\
		&\left(u^\eps, v^\eps,T^\eps \right)|_{z=0}=\left(u^\eps, v^\eps,T^\eps \right)|_{z=1}=0.
	\end{aligned}
	\right.
\end{equation}

Formally taking $\eps\to 0$ in the system $\eqref{eq:hydroPE}$, we obtain the following hydrostatic limit for primitive equations, which are combination of a Prandtl-like system with a transport-diffusion equation of the temperature
\begin{equation}
	\label{eq:hydrolimit}
	\left\{\;
	\begin{aligned}
		&\p_tu+u\p_x u+v\p_zu-\p_z^2u+\p_xp=0 && \mbox{in } \; ]0,\infty[ \times \mathcal{S}\\
		&\p_zp= T && \mbox{in } \; ]0,\infty[ \times \mathcal{S}\\
		&\p_t T+u\partial_{x} T+v\partial_{z} T-\pa_x^2 T-\pa_z^2 T=0 && \mbox{in } \; ]0,\infty[ \times \mathcal{S}\\
		&\p_xu+\p_zv=0  && \mbox{in } \; ]0,\infty[ \times \mathcal{S}\\
		&(u,T)|_{t=0}=(u_0,T_0)  && \mbox{in } \; \cS\\
		&\left(u, v,T\right)|_{z=0}=\left(u, v,T\right)|_{z=1}=0.
	\end{aligned}
	\right.
\end{equation}

\medskip

We want to recall some results on the well-posedness of the system $\eqref{NS2}$. This system was studied by Lions, Temam and Wang in \cite{Temam,Wang,Lions}, where the authors considered full viscosity and diffusivity, and established the global existence of weak solutions. In the two-dimensional case, the local existence of strong solutions was proved by Guill\'en-Gonz\'alez, Masmoudi and Rodriguez-Bellido \cite{30}, while the global existence was achieved by Bresch, Kazhikhov and Lemoine in \cite{6} and by Temam and Ziane in \cite{Ziane}. In our work, we will study the global well posedness of the system \eqref{NS2} in the 2D thin strip $\mathcal{S}^\eps$ when $\eps$ is close to zero. The equivalent result is also available for the rescaled system \eqref{eq:hydroPE}.

Concerning the hydrostatic limit system \eqref{eq:hydrolimit}, we remark the same difficulty as for Prandtl equations due to its degenerate form and the non-local nonlinear term $v\dd_zu$ which leads to the loss of one derivative in the tangential direction when one wants to perform energy estimates.
For a more complete survey on this challenging problem, we refer the reader to the works \cite{awxy, e-1,e-2,GV-D,M-W} and the references therein. The main ideas to overcome this difficulty consist in imposing a monotonicity hypothesis on the normal derivative of the velocity or an analytic regularity on the velocity.
In the pioneering works \cite{oleinik}, Oleinik and Samokhin used the Crocco transformation under the monotonicity assumption to transform Prandtl equations into a new quasilinear system and established the local existence of solutions in the Sobolev functional framework. However, the nonlinearity of the Crocco variables induces certain difficulties for understanding the nature of Prandtl equations.

Later, in \cite{samm}, Sammartino and Caflisch solved the problem for analytic solutions (analytic in both tangential and normal directions) on a half space without using the monotonicity assumption and the Crocco transformation. The analyticity in normal variable was then removed by Lombardo, Cannone and Sammartino in \cite{LCS2003}. The main argument used in \cite{samm, LCS2003} is to apply the abstract Cauchy-Kowalewskaya (CK) theorem. We also mention the well-posedness results of the Prandtl equations in Gevrey classes \cite{GvM2015, GvMV2019}. Under the monotonicity assumption, recently, Alexandre, Wang, Xu and Yang \cite{awxy} and Masmoudi and Wong \cite{M-W} obtained the existence of local smooth solutions for Prandtl equations by performing direct energy estimates in weighted Sobolev spaces and by exploiting the cancelation properties of the ``bad terms'', without using Crocco transformation.

We remark that unlike the case of Prandtl equations, in the system \eqref{eq:hydrolimit}, the pressure term is not defined by the outer flow using Bernoulli's law but by temperature via the relation $\dd_z p = T$. One of the novelties of the paper is to find a way to treat the pressure term using the temperature equation to obtain the global well posedness of our system. In the case where the temperature is constant, the well-posedness of the hydrotatic Navier-Stokes equations was studied in Gevrey classes by G\'erard-Varet, Masmoudi and Vicol in \cite{GvM2015} and recently in \cite{PZZ2019}, Paicu, Zhang and Zhang proved the global well-posedness of the hydrostatic Navier-Stokes system for small analytic data. We remark that the well-posedness of hydrostatic limit systems is still open in Sobolev settings. In this work, we will apply the method of \cite{PZZ2019} in the case of the hydrostatic primitive equations, where the temperature is not constant.

\begin{remark}
	\begin{enumerate}
		\item We first remark that the tangential pressure term is of order $1$ in the first equation of \eqref{eq:hydrolimit}. Let us suppose for now that we have the necessary regularity to perform the next calculations. Taking the $L^2$-scalar product of $\dd_x p$ with $u$ and performing integrations by parts, we obtain
		\begin{align*}
			\psca{\pa_x p,u}_{L^2} &= -\psca{ p,\pa_x u}_{L^2} = \psca{ p,\pa_z v}_{L^2} \\
			& = -\psca{ \pa_z p, v}_{L^2} = -\psca{ T,\int_0^z -\pa_x u \ dz'}_{L^2}\\
			& = -\psca{ \pa_x T,\int_0^z u \ dz'}_{L^2}.
		\end{align*}
		So, in order to control $\p_x p$, we need a control of $\dd_x T$ (or more precisely, we need at least a control of $\abs{D_x}^{\frac{1}2}T$ as explained in the next calculation). That is the reason why we consider an anisotropic Laplacian term $\Delta_\eps$ in the temperature equation.		
		%We also remark that in our energy estimates, we do not need to control the normal derivative $\dd_z T$, that means that we can simply suppose that $\eps = 0$ in the Laplacian term $\Delta_\eps$.
		
		\medskip
		
		\item We remark a particular case where we can still get a control of $\dd_x p$ without the consideration of an anisotropic diffusion on the temperature. Indeed, if we suppose that the pressure satisfies a hydrostatic law of the type $\pa_z p = f(t) T$, where the function $f(t)$ is integrable, then we can perform the following estimate
		\begin{align*}
			\psca{\pa_x p,u}_{L^2} &= -\psca{ p,\pa_x u}_{L^2} = \psca{ p,\pa_z v}_{L^2} \\
			& = -\psca{ \pa_z p, v}_{L^2} = \psca{ f(t) T,\int_0^z \pa_x u \ dz'}_{L^2}\\
			& \cong -f(t) \psca{ \abs{D_x}^\f12 T,\int_0^z \abs{D_x}^\f12 u \ dz'}_{L^2}.
		\end{align*}
		The method of \cite{PZZ2019} (that we will resume in the next section) and the functional settings used in our paper allow to gain an addition ``diffusion-type'' term, which implies the necessary controls on  $\abs{D_x}^\f12 T$ and $\abs{D_x}^\f12 u$. 
	\end{enumerate}
\end{remark}

\subsection{Functional framework}

%{\color{black}{}
	
In order to introduce our results, we will briefly recall some elements of the Littlewood-Paley theory and introduce the function spaces and techniques that we are going to use throughout this paper. Let $\psi$ be an even smooth function in $C^{\infty}_0(\RR)$ such that the support is contained in the ball $B_{\RR}(0,\frac{4}3)$ and $\psi$ is equal to 1 on a neighborhood of the ball $B_{\RR}(0,\frac{3}4)$. Let $\chi(z) = \psi\pare{\frac{z}2} - \psi(z)$. Thus, the support of $\chi$ is contained in the ring $\set{z \in \RR: \frac{3}4 \leq \abs{z} \leq \frac{8}3}$, and $\chi$ is identically equal to 1 on the ring $\set{z \in \RR: \frac{4}3 \leq \abs{z} \leq \frac{3}2}$. The functions $\psi$ and $\chi$ enjoy the following important properties
\begin{align} 
	\label{nonhomodecomp}
	&\forall z\in \RR, \quad \psi(z) + \sum_{j\in\NN} \chi (2^{-j} z) = 1,\\
	\label{decomposition}
	&\forall z\neq 0, \quad \sum_{j\in\ZZ} \chi (2^{-j} z) = 1.
\end{align}
Here, \eqref{nonhomodecomp} is used to define non-homogeneous Sobolev (and Besov) spaces and \eqref{decomposition} is important for the definition of homogeneous Sobolev (and Besov) spaces. We also remark that $\chi$ satisfies the following quasi-orthogonal relation
\begin{equation*}
	\forall\, j,j' \in \NN, \; \abs{j-j'}\geq 2, \quad \supp \chi(2^{-j} \cdot) \cap \supp \chi(2^{-j'} \cdot) = \emptyset.
\end{equation*}
Let $\Fcal_h$ and $\Fcal_h^{-1}$ be the Fourier transform and the inverse Fourier transform respectively in the horizontal direction. We will also use the notation $\widehat{u} = \mathcal{F}_hu$. We introduce the following definitions of the homogeneous dyadic cut-off operators.
\begin{definition}
	\label{def:dyadic} For all tempered distribution u in the horizontal direction (of $x$ variable) and for all $q \in \ZZ$, we set
	\begin{align*}
		\Delta^h_q u(x,z) &= \mathcal{F}_h^{-1}\pare{\chi(2^{-q}\abs{\xi})\widehat{u}(\xi,z)},\\
		S^h_q u(x,z) &= \sum_{l \leq q-1} \Delta^h_l u(x,z).
	\end{align*}
\end{definition}
We refer to \cite{BCD2011book} and \cite{B1981} for a more detailed construction of the dyadic decomposition. This definition, combined with the equality \eqref{decomposition}, implies that all tempered distributions can be decomposed with respect to the horizontal frequencies as $$u = \sum_{q \in \ZZ} \Delta_q^h u.$$ 
The following Bernstein lemma gives important properties of a distribution $u$ when its Fourier transform is well localized. We refer the reader to \cite[Lemma 2.1.1]{C1995book} for the proof and for other comments.
\begin{lemma}
	\label{le:Bernstein}
	Let $k\in\NN$, $d \in \NN^*$ and $r_1, r_2 \in \RR$ satisfy $0 < r_1 < r_2$. There exists a constant $C > 0$ such that, for any $a, b \in \RR$, $1 \leq a \leq b \leq +\infty$, for any $\lambda > 0$ and for any $u \in L^a(\RR^d)$, we have
	\begin{equation*} %\label{eq:Bernstein1}
		\supp\pare{\widehat{u}} \subset \set{\xi \in \RR^d \;\vert\; \abs{\xi} \leq r_1\lambda} \quad \Longrightarrow \quad \sup_{\abs{\alpha} = k} \norm{\dd^\alpha u}_{L^b} \leq C^k\lambda^{k+ d \pare{\frac{1}a-\frac{1}b}} \norm{u}_{L^a},
	\end{equation*}
	and
	\begin{equation*} %\label{eq:Bernstein2}
		\supp\pare{\widehat{u}} \subset \set{\xi \in \RR^d \;\vert\; r_1\lambda \leq \abs{\xi} \leq r_2\lambda} \quad \Longrightarrow \quad C^{-k} \lambda^k\norm{u}_{L^a} \leq \sup_{\abs{\alpha} = k} \norm{\dd^\alpha u}_{L^a} \leq C^k \lambda^k\norm{u}_{L^a}.
	\end{equation*}
\end{lemma}

We now introduce the function spaces used throughout the paper. As in \cite{PZZ2019}, we define the Besov-type spaces $\Bcal^s$, $s\in\RR$ as follows.
\begin{definition}
	Let $s \in \RR$ and $\Scal = \RR \times ]0,1[$. For any $u \in S'_h(\Scal)$, \emph{i.e.}, $u$ belongs to $S'(\Scal)$ and $\lim_{q\to -\infty} \norm{S^h_q u}_{L^\infty} = 0$, we set
	\begin{equation*}
		\norm{u}_{\Bcal^s} \stackrel{\tiny def}{=} \; \sum_{q \in \ZZ} 2^{qs} \norm{\DD_q^h u}_{L^2}.
	\end{equation*}
	\begin{enumerate}
		\item[(i)] For $s \leq \frac{1}2$, we define 
		\begin{equation*}
			\Bcal^s(\Scal) \stackrel{\tiny def}{=} \; \set{u \in S'_h(\Scal) \ : \ \norm{u}_{\Bcal^s} < +\infty}.
		\end{equation*}
		\item[(ii)] For $s \in \; ]k-\frac12, k+\frac12]$, with $k \in \NN^*$, we define $\Bcal^s(\Scal)$ as the subset of distributions $u$ in $S'_h(\Scal)$ such that $\dd_x^ku \in \Bcal^{s-k}(\Scal)$.
	\end{enumerate}
\end{definition}  
For a better use of the smoothing effect given by the diffusion terms, we will need the following time-weighted Chemin-Lerner-type spaces.
\begin{definition}
	\label{def:CLweight}
	Let $p \in [1,\infty]$ and let $f \in L^1_{\tiny loc}(\RR_+)$ be a non-negative function. Then, the space $\tilde{L}^p_{t,f(t)}(\Bcal^s(\Scal))$ is the closure of $C([0,t];C^\infty_0(\Scal))$ under the norm
	\begin{equation*}
		\norm{u}_{\tilde{L}^p_{t,f(t)}(\Bcal^s(\Scal))} \stackrel{\tiny def}{=} \; \sum_{q\in\ZZ} 2^{qs} \pare{\int_0^t f(t') \norm{\DD^h_q u(t')}_{L^2}^p dt'}^{\frac{1}{p}},
	\end{equation*}
	with the usual change if $p = \infty$. In the case where $f \equiv 1$, we will simply use the notation $\tilde{L}^p_t(\Bcal^s(\Scal))$.
\end{definition}
The following estimates are a direct consequence of the above definition.
\begin{proposition}
	\label{prop:dqweight}
	Let $p \in [1,\infty]$ and let $f \in L^1_{\tiny loc}(\RR_+)$ be a non-negative function. If $u \in \tilde{L}^p_{t,f(t)}(\Bcal^s(\Scal))$ then, there exists a sequence $(d_q(u))_{q\in\ZZ}$ such that $\sum d_q(u) = 1$ and
	\begin{equation*}
		\pare{\int_0^t f(t') \norm{\DD^h_q u(t')}_{L^2}^p dt'}^{\frac{1}{p}} \leq d_q(u) 2^{-qs} \norm{u}_{\tilde{L}^p_{t,f(t)}(\Bcal^s(\Scal))}, \quad \forall \, q\in \ZZ.
	\end{equation*}
\end{proposition}

}

\subsection{Main results}

{\color{black}
	
The difficulty when one wants to estimate the nonlinear terms relies in finding a way to compensate the loss of one derivative. The main idea here is to exploit the smoothing effect given by the above function spaces, which allow to gain ``one half derivative''. 
Using the method introduced by Chemin in \cite{C2004} (see also \cite{CGP2011}, \cite{PZ2011} or \cite{PZZ2019}), for any $f \in L^2(\Scal)$, we define the following auxiliary function, which allows to control the analyticity of $f$ in the horizontal variable $x$, 
\begin{align} \label{eq:AnPhi}
	f_{\phi}(t,x,z) = e^{\phi(t,D_x)} f(t,x,z) \stackrel{\tiny def}{=} \mathcal{F}^{-1}_h(e^{\phi(t,\xi)} \widehat{f}(t,\xi,z)) \text{ \ \ with \ \ } \phi(t,\xi) = (a- \lambda \rho(t))|\xi|,
\end{align}
where the quantity $\rho(t)$, which describes the evolution of the analytic band of $f$, satisfies
\begin{equation}
	\label{eq:AnTheta} \forall\, t > 0, \ \dot{\rho}(t) \geq 0 \text{ \ \ and \ \ } \rho(0) = 0. 
\end{equation}
Here, the constants $a > 0$ and $\lambda > 0$ are given and in particular, $a$ represents the width of the initial analytic band. We will chose $\lambda$ according to the necessity of the ``bootstrap'' process in the next sections.

We remark that if we differentiate, with respect to the time variable, a function of the type $e^{\phi(t,D_x)} f(t,x,z)$, we obtain an additional ``good term'' which allows to gain ``one half derivative''. More precisely, we have
\begin{equation}
	\label{eq:derivativePhi}
	\frac{d}{dt} \pare{e^{\phi(t,D_x)} f(t,x,z)} = -\dot{\rho}(t) \abs{D_x} e^{\phi(t,D_x)} f(t,x,z) + e^{\phi(t,D_x)} \dd_t f(t,x,z),
\end{equation}
where $-\dot{\rho}(t) \abs{D_x} e^{\phi(t,D_x)} f(t,x,z)$ plays a smoothing role if $\dot{\rho}(t) \geq 0$. This smoothing effect allows to obtain our global existence and stability results in the analytic framework.

%We recall that in the Prandtl case, only have the local %existence and the convergence is still an open question! %Besides, Prandtl system is known to be very unstable.

\bigskip

From now on, let $\mathcal{K} > 0$ be the Poincar\'e constant on the strip $\mathcal{S}$, in the sens that, for any $f \in L^2(\mathcal{S})$, $f_{\vert_{\dd \mathcal{S}} = 0}$ and $\dd_zf \in L^2(\mathcal{S})$, we have $$\norm{f}_{L^2(\mathcal{S})} \leq \mathcal{K}\norm{\dd_z f}_{L^2(\mathcal{S})}.$$ Our main results are the following theorems.
\begin{theorem}[Global well-posedness of the hydrostatic limit system]
	\label{th:hydrolimit}
	Let $a > 0$, $s > 0$ and assume that $e^{a|D_x|} (u_0,T_0) \in \mathcal{B}^\f12 \cap \mathcal{B}^{\frac{3}2} \cap \mathcal{B}^s$. There exist positive constants $c_0$, $C$ and a decreasing function $\tilde{\phi}: \RR_+ \to [\frac{2a}3,a]$
	such that, if we suppose that the initial data $(u_0,T_0)$ satisfy the compatibility condition $\int_0^1 u_0 dz =0$ and the smallness assumption 
	\begin{align*}
		\norm{e^{a|D_x|} u_0}_{\mathcal{B}^\f12} +  \norm{e^{a|D_x|}  T_{0}}_{\mathcal{B}^\f12} \leq c_0a \quad \mbox{and}\quad \norm{e^{a|D_x|} u_0}_{\mathcal{B}^{\frac{3}2}} +  \norm{e^{a|D_x|}  T_{0}}_{\mathcal{B}^{\frac{3}2}} \leq c_0,
	\end{align*}
	then the system \eqref{eq:hydrolimit} has a unique global solution 
	\begin{displaymath}
		(u,T) \in \tilde{L}^{\infty}(\mathbb{R}_+;\mathcal{B}^s) \cap C(\mathbb{R}_+;\mathcal{B}^s) \quad \mbox{with}\quad \dd_zu \in \tilde{L}^2(\mathbb{R}_+;\mathcal{B}^s),
	\end{displaymath} 
	satisfying
	\begin{align}
		\label{eq:hydrolimitU}
		\Vert e^{\mathcal{R}t} (u_{\phi},T_{\phi}) \Vert_{\tilde{L}^{\infty}(\mathbb{R}_+;\mathcal{B}^s)} + \frac{1}{4} \norm{e^{\Rcal t} \pa_z u_{\phi}}_{\tilde{L}^2_t(\Bcal^s)} + \frac{1}{2}\norm{e^{\Rcal t} \nabla T_{\phi}}_{\tilde{L}^2_t(\Bcal^s)} 
		\leq 2C \Vert e^{a|D_x|} (u_0,T_{0}) \Vert_{\mathcal{B}^s},
	\end{align}
	for any $0 \leq \mathcal{R} \leq \frac{1}{2\mathcal{K}}$, where $\phi(t,\xi) = \tilde{\phi}(t) \abs{\xi}$ and where for any $f\in L^2(\Scal)$,
	\begin{align*}
		f_{\phi}(t,x,z) = e^{\phi(t,D_x)} f(t,x,z) = \mathcal{F}^{-1}_h(e^{\phi(t,\xi)} \widehat{f}(t,\xi,z)).
	\end{align*}
	
	Furthermore, we have, 
	\begin{multline}
		\label{eq:hydrolimitdtU}
		\Vert e^{\mathcal{R}t} (\pa_tu)_\phi \Vert_{\tilde{L}^2_t(\mathcal{B}^s)} + \Vert e^{\mathcal{R}t} \pa_z u_\phi \Vert_{\tilde{L}^\infty_t(\mathcal{B}^s)}
		\\
		\leq C\pare{ \norm{e^{a|D_x|} (u_0,T_0)}_{\mathcal{B}^s} + \norm{e^{a|D_x|} (u_0,T_0)}_{\mathcal{B}^\f12} \norm{e^{a|D_x|} (u_0,T_0)}_{\mathcal{B}^{s+1}} }.
	\end{multline}
\end{theorem}

\begin{remark}
	We remark that the normal component $v$ is uniquely determined from the incompressibility and the boundary condition
	\begin{align} \label{vv}
		v(t,x,z) = \int_0^z \pa_z v(t,x,z')dz' = -\int_0^z \pa_x u(t,x,z') dz' .
	\end{align} 
\end{remark}

\medskip

\begin{theorem}[Global well-posedness of the primitive system]
	\label{th:primitive}
	Let $a > 0$, $s > 0$, $\eps > 0$ and assume that $e^{a|D_x|} (u_0^\eps,v_0^\eps,T_0^\eps) \in \mathcal{B}^\f12 \cap \mathcal{B}^{\frac{3}2} \cap \mathcal{B}^s$. There exist positive constants $c_1$, $C$ (independent of $\eps$) and a decreasing function $\tilde{\Theta}: \RR_+ \to [\frac{2a}3,a]$
	such that, if we suppose that the initial data $(u_0^\eps,v_0^\eps,T_0^\eps)$ satisfy 
	\begin{align*}
		\norm{e^{a|D_x|} (u_0^\eps,v_0^\eps,T_0^\eps)}_{\mathcal{B}^\f12} \leq c_1a \quad \mbox{and}\quad \norm{e^{a|D_x|} (u_0^\eps,v_0^\eps,T_0^\eps)}_{\mathcal{B}^{\frac{3}2}} \leq c_1,
	\end{align*}
	then, for any $0 < \eps < \frac{1}{2C}$ the system \eqref{eq:hydroPE} has a unique global solution $(u^\eps,v^\eps,T^\eps)$ satisfying,
	\begin{align*}
		\Vert e^{\mathcal{R}t}(u^\eps_{\Theta},\eps v^\eps_{\Theta},T^\eps_{\Theta}) \Vert_{\tilde{L}^{\infty}_t(\mathcal{B}^s)} + \Vert e^{\mathcal{R}t} \pa_z (u^\eps_{\Theta},\eps v^\eps_{\Theta},T^\eps_{\Theta})\Vert_{\tilde{L}^{2}_t(\mathcal{B}^s)} + \Vert e^{\mathcal{R}t} \dd_x T^\eps_{\Theta} \Vert_{\tilde{L}_t(\mathcal{B}^s)} \leq  C\Vert e^{a|D_x|} (u^\eps_{0},\eps v^\eps_{0},T^\eps_{0}) \Vert_{\mathcal{B}^s},
	\end{align*}
	for any $0 \leq \mathcal{R} \leq \frac{1}{2\mathcal{K}}$. Here, $\Theta(t,\xi) = \tilde{\Theta}(t) \abs{\xi}$ and for any $f \in L^2(\Scal)$, 
	\begin{align*} 
		f_{\Theta}(t,x,z) = e^{\Theta(t,D_x)} f(t,x,z) = \mathcal{F}^{-1}_h(e^{\Theta(t,\xi)} \widehat{f}(t,\xi,z)).
	\end{align*}
\end{theorem}

\begin{theorem}[Convergence to the hydrostatic limit system]
	\label{th:Convergence}
	Let $a > 0$ and $0 < \eps \leq 1$. We suppose that the initial data $(u_0,v_0,T_0)$ and $(u^\eps_0,v^\eps_0,T^\eps_0)$ satisfy the assumptions of Theorems \ref{th:hydrolimit} and \ref{th:primitive}. Let $(u,v,T)$ and $(u^\eps,v^\eps,T^\eps)$ be the respective solutions of the systems \eqref{eq:hydrolimit} and \eqref{eq:hydroPE}. Then, there exist a constant $M > 0$ independent of $\eps$ and a decreasing function $\tilde{\varphi}: \RR_+ \to [\frac{a}3,a]$ such that
	\begin{multline*} 
		\Vert (u^\eps_\varphi - u_\varphi,\eps v^\eps_\varphi -\eps v_\varphi) \Vert_{\tilde{L}^\infty_t(\mathcal{B}^s)} + \Vert \pa_z(u^\eps_\varphi - u_\varphi,\eps v^\eps_\varphi -\eps v_\varphi) \Vert_{\tilde{L}^2_t(\mathcal{B}^s)} +\eps \Vert (u^\eps_\varphi - u_\varphi,\eps v^\eps_\varphi -\eps v_\varphi) \Vert_{\tilde{L}^2_t(\mathcal{B}^{s+1})} \\
		\leq C \left( \Vert e^{a|D_x|}(u_0^\eps - u_0, \eps(v_0^\eps - v_0)) \Vert_{\mathcal{B}^s}+ C\Vert e^{a|D_x|} (T_0^\eps - T_0) \Vert_{\mathcal{B}^s} + M\eps  \right). 
	\end{multline*}
	where $\varphi(t,\xi) = \tilde{\varphi}(t) \abs{\xi}$ and where, for any $f \in L^2(\Scal)$,
	\begin{align*} 
		f_{\varphi}(t,x,z) = e^{\varphi(t,D_x)} f(t,x,z) = \mathcal{F}^{-1}_h(e^{\varphi(t,\xi)} \widehat{f}(t,\xi,z)).
	\end{align*}
\end{theorem}

\subsection{Organisation of the paper} 

Our paper will be divided into several sections as follows. In the next section, we establish the needed nonlinear estimates which will be used throughout this paper. 	In Section \ref{se:hydrolim}, we prove the global wellposedness of the system \eqref{eq:hydrolimit} for small data in analytic framework. Section \ref{se:PEhydro} is devoted to the study of the system \eqref{eq:hydroPE} and the proof of Theorem \ref{th:primitive}. In Section \ref{se:convergence}, we prove the convergence of the system \eqref{eq:hydroPE} towards the system \eqref{eq:hydrolimit} when $\eps$ goes to $0$. Finally, in the appendix, we give the proofs of Lemmas \ref{lem:uww} and \ref{lem:vww}.

We end this introduction by the notations that will be used in all that follows. For $f\lesssim g$, we mean that there is a positive constant $C$, which may be different from line to line, such that $f\leq Cg$. We denote by $\psca{f,g}_{L^2}$ the inner product of $f$ and $g$ in $L^2(\mathcal{S})$. Finally, we denote by $(d_q)_{q\in\ZZ}$(resp. $(d_q(t))_{q\in\ZZ}$) to be a generic element of $ \ell^1(\ZZ)$ so that $\sum_{q\in\ZZ} d_q = 1$ (resp. $\sum_{q\in\ZZ} d_q(t) = 1$).

}

\section{Nonlinear estimates}

%{\color{darkmagenta}}

{\color{black}
	
The proofs of our main theorems rely on the following lemmas \ref{lem:uww} and \ref{lem:vww} which give controls of the nonlinear terms in our analytic norms.

\begin{lemma} \label{lem:uww}
	Let $s > 0$ and $\phi: \RR_+\times \RR \to \RR_+$. There exists a constant $C \geq 1$ such that, for any functions $u$, $w$ and $\overline{w}$ which are defined on $\RR_+\times \mathcal{S}$, $u\vert_{\dd\Scal} = w\vert_{\dd\Scal} = \overline{w}_{\dd\Scal} = 0$ and satisfy 
	\begin{displaymath}
		\dd_z u_\phi(t), \dd_z w_\phi(t) \in \mathcal{B}^{\frac{1}{2}}, \ \forall t \geq 0, \quad \mbox{and}\quad w, \overline{w} \in \tilde{L}^2_{t,f(t)} (\mathcal{B}^{s+\frac{1}{2}}),
	\end{displaymath}
	we have, for any $0 \leq \mathcal{R} \leq \frac{1}{2\mathcal{K}}$ and for any $q\in \ZZ$,
	\begin{multline*}
		\int_0^t \abs{\psca{e^{\Rcal t'} \Delta_q^h (u\p_x w)_{\phi}, e^{\Rcal t'} \Delta_q^h w_{\phi}}_{L^2}} dt' 
		\\
		\leq C d_q^2 2^{-2qs} \pare{\norm{e^{\Rcal t} u_{\phi}}_{\tilde{L}^2_{t,f(t)}(\mathcal{B}^{s+\frac{1}{2}})} + \norm{e^{\Rcal t} w_{\phi}}_{\tilde{L}^2_{t,f(t)}(\mathcal{B}^{s+\frac{1}{2}})}} \norm{e^{\Rcal t} \overline{w}_{\phi}}_{\tilde{L}^2_{t,f(t)}(\mathcal{B}^{s+\frac{1}{2}})},
	\end{multline*}
	where $u_\phi$ is determined by \eqref{eq:AnPhi}, where 
	\begin{displaymath}
		f \in L^1_{\tiny loc}(\RR_+) \quad\mbox{and}\quad f(t) \geq \max\set{\norm{\dd_z u_\phi(t)}_{\mathcal{B}^{\frac{1}{2}}}, \norm{\dd_z w_\phi(t)}_{\mathcal{B}^{\frac{1}{2}}}}, \ \forall t \geq 0,
	\end{displaymath}
	and $(d_q)_{q\in\ZZ}$ is a positive sequence with $\sum_{q\in\ZZ} d_q = 1$. 
\end{lemma}

\begin{remark}
	\label{rem:s<1}
	In the case where $0 < s < 1$, we can relax the condition $f(t) \geq \norm{\dd_z w_\phi(t)}_{\mathcal{B}^{\frac{1}{2}}}$ and obtain, for any $w, \overline{w} \in \tilde{L}^2_{t,f(t)} (\mathcal{B}^{s+\frac{1}{2}})$, the following estimate (as in \cite{PZZ2019})
	\begin{equation*}
		\int_0^t \abs{\psca{e^{\Rcal t'} \Delta_q^h (u\p_x w)_{\phi}, e^{\Rcal t'} \Delta_q^h \overline{w}_{\phi}}_{L^2}} dt' 
		\\
		\leq C d_q^2 2^{-2qs} \norm{e^{\Rcal t} w_{\phi}}_{\tilde{L}^2_{t,f(t)}(\mathcal{B}^{s+\frac{1}{2}})} \norm{e^{\Rcal t} \overline{w}_{\phi}}_{\tilde{L}^2_{t,f(t)}(\mathcal{B}^{s+\frac{1}{2}})}.
	\end{equation*}
\end{remark}

Before starting the calculations, we remark that one should be careful when dealing with the product of the type $(fg)_\phi$, since in general, we have $(fg)_\phi \neq f_\phi g_\phi$. However, we still can ``compare'' this two products as in the following lemma.
\begin{lemma}
	\label{lem:f+}
	For any $f,g\in L^2_x$, we set
	\begin{displaymath}
		f^+ = \mathcal{F}_\xi^{-1}(|\mathcal{F}_x(f)|) \quad \mbox{and} \quad g^+ = \mathcal{F}_\xi^{-1}(|\mathcal{F}_x(g)|).
	\end{displaymath}
	Then, we have 
	\begin{displaymath}
		|(\widehat{fg})_\phi(\xi)| \leq \widehat{f^+_\phi g^+_\phi }(\xi).
	\end{displaymath}
\end{lemma} 
\begin{proof}
	We have 
	\begin{align*}
		|(\widehat{fg})_\phi(\xi)| = e^{\phi(\xi)} |\widehat{f}(.)*\widehat{g}(.)(\xi)| \leq e^{\phi(\xi)} \int |\widehat{f}(\xi - \eta)||\widehat{g}(\eta)|d\eta.
	\end{align*}
	From the definition of the function $\phi$, we have $e^{\phi(\xi)} \leq e^{\phi(\xi-\eta)}e^{\phi(\eta)}$, then 
	\begin{align*}
		|(\widehat{fg})_\phi(\xi)| &\leq \int e^{\phi(\xi-\eta)}|\widehat{f}(\xi - \eta)| e^{\phi(\eta)}|\widehat{g}(\eta)|d\eta \\
		&\leq \int |\widehat{f}_\phi(\xi - \eta)||\widehat{g}_\phi(\eta)|d\eta \\
		&\leq |\widehat{f}_\phi| *|\widehat{g}_\phi|(\xi) = \widehat{f}^+_\phi * \widehat{g}^+_\phi = \widehat{f^+_\phi g^+_\phi }(\xi).
	\end{align*}
\end{proof}
\begin{remark}
	\label{rem:f+}
	For any $f\in L^2_x(\mathbb{R})$, we have $\|f^+\|_{L^2_x} = \|f\|_{L^2_x}$. Since the norms of the function spaces used in our work are based on the $L^2$-norm (with respect to the space variables), without loss of generality, from now on, we can assume that $\widehat{f} \geq 0$.
\end{remark}

\noindent\textit{Proof of Lemma \ref{lem:uww}.}

We first recall the Bony homogeneous decomposition into paraproducts and remainders (see for instance \cite[Chapter 2]{BCD2011book}) in the tangential direction. For two tempered distributions $X$ and $Y$, we have
$$ XY = \Tcal^h_{X}Y + \Tcal^h_{Y}X + \Rcal^h(X,Y),$$
where
\begin{align*}
	\Tcal^h_XY = \sum_{q\in\ZZ} S_{q-1}^h X \DD_q^h Y \quad\mbox{ and }\quad \Rcal^h(X,Y) = \sum_{\abs{q'-q} \leq 1} \DD_q^h X \DD_{q'}^h Y.
\end{align*}

Using this decomposition, we write
\begin{align*}
	\int_0^t \abs{\psca{e^{\Rcal t'} \Delta_q^h (u\p_x w)_{\phi}, e^{\Rcal t'} \Delta_q^h \overline{w}_{\phi}}_{L^2}} dt' \leq A_{1,q} + A_{2,q} + A_{3,q}, 
\end{align*}
where
\begin{align*}
	A_{1,q} &= \int_0^t \abs{\psca{ e^{\Rcal t'} \Delta_q^h ( \Tcal^h_{u}\pa_x w)_{\phi}, e^{\Rcal t'} \Delta_q^h \overline{w}_{\phi}}_{L^2}} dt'\\
	A_{2,q} &= \int_0^t \abs{\psca{ e^{\Rcal t'} \Delta_q^h ( \Tcal^h_{\pa_xw} u )_{\phi}, e^{\Rcal t'} \Delta_q^h \overline{w}_{\phi}}_{L^2}} dt'\\
	A_{3,q} &= \int_0^t \abs{\psca{ e^{\Rcal t'} \Delta_q^h (\Rcal^h(u,\pa_xw))_{\phi}, e^{\Rcal t'} \Delta_q^h \overline{w}_{\phi}}_{L^2}} dt'.
\end{align*}

Using the support properties of the Fourier transform of the dyadic operators $\Delta^h_q$ (see for example \cite[Proposition 2.10]{B1981}), Lemma \ref{lem:f+} and Remark \ref{rem:f+}, we get
\begin{align*}
	A_{1,q} \leq \sum_{|q-q'| \leq 4} \int_0^t e^{2\Rcal t'} \Vert S^h_{q'-1} u_{\phi} (t') \Vert_{L^{\infty}}\Vert \Delta_{q'}^h \partial_xw_{\phi} (t') \Vert_{L^2} \Vert \Delta_q^h \overline{w}_{\phi} (t') \Vert_{L^2} dt'.
\end{align*} 
For any $l \in \ZZ$, Sobolev inclusion $\dot{H}^1_z([0,1]) \hookrightarrow L^\infty_z([0,1])$ and the definition of $\Bcal^{\frac{1}2}$ yield
\begin{equation} \label{eq:normB12}
	\Vert \Delta^h_l u_{\phi}(t') \Vert_{L^{\infty}} \lesssim 2^{\frac{l}{2}} \Vert \Delta^h_l u_{\phi}(t') \Vert_{L^2_x(L^{\infty}_z)} \lesssim 2^{\frac{l}{2}} \Vert \Delta^h_l \pa_z u_{\phi}(t') \Vert_{L^2} \lesssim d_q(u_\phi) \Vert \pa_z u_{\phi}(t') \Vert_{\mathcal{B}^{\f12}},
\end{equation}
where $\set{d_q(u_\phi)}$ is a positive sequence with $\sum d_q(u_\phi) = 1$ and then,
\begin{align*}
	\Vert S^h_{q'-1} u_{\phi} (t') \Vert_{L^{\infty}} = \Big\Vert \sum_{l\leq q'-2} \Delta^h_l u_\phi(t') \Big\Vert_{L^{\infty}} \lesssim  \Vert \pa_z u_{\phi}(t') \Vert_{\mathcal{B}^{\f12}}.
\end{align*}
Using Cauchy-Schwarz inequality and Proposition \ref{prop:dqweight}, we obtain
\begin{align}
	\label{eq:A1q} A_{1,q} &\lesssim \sum_{|q-q'| \leq 4} 2^{q'} \int_0^t \Vert \pa_z u_{\phi}(t') \Vert_{\mathcal{B}^{\f12}} e^{\Rcal t'} \Vert \Delta_{q'}^h w_{\phi}(t') \Vert_{L^2} e^{\Rcal t'} \Vert \Delta_q^h \overline{w}_{\phi}(t') \Vert_{L^2} dt'\\
	&\lesssim \sum_{|q-q'| \leq 4} 2^{q'} \left(\int_0^t \Vert \pa_z u_{\phi} \Vert_{\mathcal{B}^{\f12}} e^{2\Rcal t'} \Vert \Delta_{q'}^h w_{\phi} \Vert_{L^2}^2 dt'\right)^{\f12} \left(\int_0^t \Vert  \pa_z u_{\phi} \Vert_{\mathcal{B}^{\f12}} e^{2\Rcal t'} \Vert \Delta_q^h \overline{w}_{\phi} \Vert_{L^2}^2 dt'\right)^{\f12} \notag\\
	&\lesssim \sum_{|q-q'| \leq 4} 2^{q'} \left(\int_0^t f(t') e^{2\Rcal t'} \Vert \Delta_{q'}^h w_{\phi} \Vert_{L^2}^2 dt'\right)^{\f12} \left(\int_0^t f(t') e^{2\Rcal t'} \Vert \Delta_q^h \overline{w}_{\phi} \Vert_{L^2}^2 dt'\right)^{\f12} \notag\\
	&\lesssim 2^{-2qs} d_q^2 \Vert e^{\Rcal t} w_{\phi} \Vert_{\tilde{L}^{2}_{t,f(t)}(\mathcal{B}^{s+\frac{1}{2}})} \Vert e^{\Rcal t} \overline{w}_{\phi} \Vert_{\tilde{L}^{2}_{t,f(t)}(\mathcal{B}^{s+\frac{1}{2}})}, \notag
\end{align}
where 
\begin{align*}
	d_q^2 = d_q(\overline{w}_\phi) \sum_{|q-q'| \leq 4} d_{q'}(w_\phi) 2^{(q-q')(s-\f12)}.
\end{align*}

\smallskip

Using the symmetry of $\Tcal^h_{X}Y$ and $\Tcal^h_{Y}X$, we can estimate $A_{2,q}$ in the similar way as $A_{1,q}$. First, Bernstein lemma \ref{le:Bernstein} and similar calculations as in \eqref{eq:normB12} give
\begin{displaymath}
	\Vert S^h_{q'-1} \dd_x w_{\phi} (t') \Vert_{L^{\infty}} \lesssim 2^{q'} \Vert S^h_{q'-1} w_{\phi} (t') \Vert_{L^{\infty}} \lesssim 2^{q'} \Vert \pa_z w_{\phi}(t') \Vert_{\mathcal{B}^{\f12}}.
\end{displaymath}
Then, we obtain
\begin{align}
	\label{eq:A2q}
	A_{2,q} &\lesssim \sum_{|q-q'| \leq 4} \int_0^t e^{2\Rcal t'} \Vert S^h_{q'-1} \partial_xw_{\phi} (t') \Vert_{L^{\infty}}\Vert \Delta_{q'}^h  u_{\phi} (t') \Vert_{L^2} \Vert \Delta_q^h \overline{w}_{\phi} (t') \Vert_{L^2} dt'\\
	&\lesssim \sum_{|q-q'| \leq 4} 2^{q'} \left(\int_0^t \Vert \pa_z w_{\phi} \Vert_{\mathcal{B}^{\f12}} e^{2\Rcal t'} \Vert \Delta_{q'}^h u_{\phi} \Vert_{L^2}^2 dt'\right)^{\f12} \left(\int_0^t \Vert  \pa_z w_{\phi} \Vert_{\mathcal{B}^{\f12}} e^{2\Rcal t'} \Vert \Delta_q^h \overline{w}_{\phi} \Vert_{L^2}^2 dt'\right)^{\f12} \notag\\
	&\lesssim \sum_{|q-q'| \leq 4} 2^{q'} \left(\int_0^t f(t') e^{2\Rcal t'} \Vert \Delta_{q'}^h u_{\phi} \Vert_{L^2}^2 dt'\right)^{\f12} \left(\int_0^t f(t') e^{2\Rcal t'} \Vert \Delta_q^h \overline{w}_{\phi} \Vert_{L^2}^2 dt'\right)^{\f12} \notag\\
	&\lesssim 2^{-2qs} d_q^2 \Vert e^{\Rcal t} u_{\phi} \Vert_{\tilde{L}^{2}_{t,f(t)}(\mathcal{B}^{s+\frac{1}{2}})} \Vert e^{\Rcal t} \overline{w}_{\phi} \Vert_{\tilde{L}^{2}_{t,f(t)}(\mathcal{B}^{s+\frac{1}{2}})}, \notag
\end{align}
where 
\begin{align*}
	d_q^2 = d_q(\overline{w}_\phi) \sum_{|q-q'| \leq 4} d_{q'}(u_\phi) 2^{(q-q')(s-\f12)}.
\end{align*}

To end the proof, it remains to estimate $A_{3,q}$. Using the support localization properties given in \cite[Proposition 2.10]{B1981}, the definition of $\Rcal^h(u,\pa_x w)$, Sobolev inclusion $\dot{H}^1_z([0,1]) \hookrightarrow L^\infty_z([0,1])$, Young inequality and Bernstein lemma \ref{le:Bernstein}, we can write
\begin{align*}
	\abs{\psca{\Delta_q^h (\Rcal^h(u,\pa_x w))_{\phi}, \Delta_q^h \overline{w}_{\phi}}_{L^2}} &\lesssim \sum_{q'\geq q-3} \Vert \DD^h_{q'} u_{\phi}\Vert_{L^{2}_xL^{\infty}_z} \Vert \Delta_{q'}^h \pa_x w_{\phi}  \Vert_{L^2} \Vert \Delta_{q}^h \overline{w}_{\phi} \Vert_{L^\infty_xL^{2}_z} \\
	&\lesssim \sum_{q'\geq q-3}  2^{\frac{q+q'}2} \pare{2^{\frac{q'}{2}} \norm{\DD^h_{q'} \pa_z u_{\phi}}_{L^{2}}} \Vert \Delta_{q'}^h w_{\phi}  \Vert_{L^{2}} \Vert \Delta_{q}^h \overline{w}_{\phi} \Vert_{L^{2}}.
\end{align*}
The definition of the $\Bcal^{\frac{1}2}$-norm implies that
\begin{displaymath}
	2^{\frac{q'}{2}} \norm{\DD^h_{q'} \pa_z u_{\phi}}_{L^{2}} \leq \norm{\pa_z u_{\phi}}_{\Bcal^{\frac{1}2}}, \quad \forall \ q' \in \ZZ.
\end{displaymath} 
Thus, as for the term $A_{1,q}$, we can perform the following estimates
\begin{align}
	\label{eq:A3q}
	A_{3,q} &\lesssim \sum_{q'\geq q-3} 2^{\frac{q+q'}2} \int_0^t \Vert \pa_z u_{\phi}(t') \Vert_{\mathcal{B}^{\f12}} e^{\Rcal t'} \Vert \Delta_{q'}^h w_{\phi}(t') \Vert_{L^2} e^{\Rcal t'} \Vert \Delta_q^h \overline{w}_{\phi}(t') \Vert_{L^2} dt'\\
	&\lesssim \sum_{q'\geq q-3} 2^{\frac{q+q'}2} \left(\int_0^t \Vert \pa_z u_{\phi} \Vert_{\mathcal{B}^{\f12}} e^{2\Rcal t'} \Vert \Delta_{q'}^h w_{\phi} \Vert_{L^2}^2 dt'\right)^{\f12} \left(\int_0^t \Vert  \pa_z u_{\phi} \Vert_{\mathcal{B}^{\f12}} e^{2\Rcal t'} \Vert \Delta_q^h \overline{w}_{\phi} \Vert_{L^2}^2 dt'\right)^{\f12} \notag\\
	&\lesssim \sum_{q'\geq q-3} 2^{\frac{q+q'}2} \left(\int_0^t f(t') e^{2\Rcal t'} \Vert \Delta_{q'}^h w_{\phi} \Vert_{L^2}^2 dt'\right)^{\f12} \left(\int_0^t f(t') e^{2\Rcal t'} \Vert \Delta_q^h \overline{w}_{\phi} \Vert_{L^2}^2 dt'\right)^{\f12} \notag\\
	&\lesssim 2^{-2qs} d_q^2 \Vert e^{\Rcal t} w_{\phi} \Vert_{\tilde{L}^{2}_{t,f(t)}(\mathcal{B}^{s+\frac{1}{2}})} \Vert e^{\Rcal t} \overline{w}_{\phi} \Vert_{\tilde{L}^{2}_{t,f(t)}(\mathcal{B}^{s+\frac{1}{2}})}, \notag
\end{align}
where
\begin{align*}
	d_q^2 = d_q(\overline{w}_\phi) \sum_{q'\geq q-3} d_{q'}(w_\phi) 2^{(q-q')s}.
\end{align*}
Here, we remark that we can write the above sum as a convolution product and using Young inequality, if we set 
$$\overline{d}_q = \sum_{q'\geq q-3} d_{q'}(w_\phi) 2^{(q-q')s},$$ then $(\overline{d}_q)_{q\in \ZZ}$ is an element of $\ell^1$.

\smallskip

Lemma \ref{lem:uww} is then a consequence of Estimates \eqref{eq:A1q}, \eqref{eq:A2q} and \eqref{eq:A3q}. \hfill $\square$

\medskip

\begin{lemma} \label{lem:vww}
	Let $s > 0$ and $\phi: \RR_+\times \RR \to \RR_+$. There exists a constant $C \geq 1$ such that, for any $(u,v,w,\overline{w})$, which are defined on $\RR_+\times \mathcal{S}$, $(u,w,\overline{w})\vert_{\dd\Scal} = 0$ and satisfy, for any $t\geq 0$, 
	\begin{equation*}
		u_\phi(t) \in \mathcal{B}^{\frac{3}{2}}, \quad \dd_z u_\phi(t), \, \dd_z w_\phi(t) \in \mathcal{B}^{\frac{1}{2}}, \quad \int_0^1 \pa_x u \, dz = 0, \quad \mbox{and} \quad u, \, w, \, \overline{w} \in \tilde{L}^2_{t,f(t)} (\mathcal{B}^{s+\frac{1}{2}}),
	\end{equation*}
	with 
	\begin{align*}
		f \in L^1_{\tiny loc}(\RR_+), \quad f(t) \geq \max\set{\norm{\dd_z u_\phi}_{\mathcal{B}^{\frac{1}{2}}}, \norm{\dd_z w_\phi}_{\mathcal{B}^{\frac{1}{2}}}} \quad \mbox{and} \quad v(t,x,z) = -\int_0^z \pa_x u(t,x,z') dz',
	\end{align*} 
	we have, for any $\mathcal{R} \geq 0$ and for any $q\in \ZZ$,
	\begin{multline}
		\label{eq:vww}
		\int_0^t \abs{\psca{e^{\mathcal{R}t'} \Delta_q^h (v\p_z w)_{\phi}, e^{\mathcal{R}t'} \Delta_q^h \overline{w}_{\phi}}_{L^2}} dt'
		\\ 
		\leq C d_q^2 2^{-2qs} \norm{e^{\mathcal{R}t} \overline{w}_{\phi}}_{\tilde{L}^2_{t,f(t)}(\mathcal{B}^{s+\frac{1}{2}})} \pare{\norm{e^{\mathcal{R}t} u_{\phi}}_{\tilde{L}^2_{t,f(t)}(\mathcal{B}^{s+\frac{1}{2}})} + \norm{u_\phi}_{L^\infty_t(\mathcal{B}^{\frac{3}{2}})}^{\frac{1}2} \norm{e^{\mathcal{R}t} \dd_z w_{\phi}}_{\tilde{L}^2_t(\mathcal{B}^s)}},
	\end{multline}
	where $(d_q)_{q\in\ZZ}$ is a positive sequence with $\sum_{q\in\ZZ} d_q = 1$.
\end{lemma}

Before proving this lemma, we remark that the incompressibility condition $\dd_x u + \dd_z v = 0$ implies a ``transfer of regularity'' from $u$ to $v$ and so we can control $v$ by mean of $u$. Indeed, we have the following result
\begin{lemma}
	\label{lem:tranreg} Let $u$, $v$ be defined on $\Scal$ with $u\vert_{\dd\Scal} = v\vert_{\dd\Scal} = 0$ and $\dd_x u + \dd_z v = 0$ such that the following terms are well defined. Then, for any $q\in\ZZ$,
	\begin{align*}
		&\norm{\DD^h_q v_\phi}_{L^2_xL^\infty_z} \leq 2^{q} \norm{\DD^h_q u_\phi}_{L^2}
		\\
		&\norm{\DD^h_q v_\phi}_{L^\infty} \leq 2^{\frac{3q}2} \norm{\DD^h_q u_\phi}_{L^2}.
	\end{align*}
\end{lemma}

\begin{proof}
	To prove the first estimate, we write
	\begin{align*}
		v = - \int_0^z \dd_x u dz'.
	\end{align*}
	Then we have,
	\begin{align*}
		\norm{\DD^h_q v_\phi}_{L^2_xL^\infty_z} \leq \int_0^1 \norm{\dd_x \DD^h_q u_\phi}_{L^2_x} dz' \lesssim 2^q \norm{\DD^h_q u_\phi}_{L^2}.
	\end{align*}
	
	Finally, we remark that we can obtain the second estimate from the first one, using Bernstein lemma \ref{le:Bernstein}.
\end{proof}

\noindent\textit{Proof of Lemma \ref{lem:vww}.}

As in the proof of Lemma \ref{lem:uww}, we decompose the term on the left-hand side of \eqref{eq:vww} as follows
\begin{align}
	\label{eq:B123q} \int_0^t \abs{\psca{e^{\mathcal{R}t'} \Delta_q^h (v\pa_z w)_{\phi}, e^{\mathcal{R}t'} \Delta_q^h \overline{w}_{\phi}}_{L^2}} dt' \leq B_{1,q} + B_{2,q} + B_{3,q},
\end{align}
with
\begin{align*} 
	B_{1,q} &= \int_0^t \abs{\psca{ e^{\mathcal{R}t'} \Delta_q^h ( \Tcal^h_{v}\pa_z w)_{\phi}, e^{\mathcal{R}t'} \Delta_q^h \overline{w}_{\phi}}_{L^2}} dt'\\
	B_{2,q} &= \int_0^t \abs{\psca{ e^{\mathcal{R}t'} \Delta_q^h ( \Tcal^h_{\pa_zw} v )_{\phi}, e^{\mathcal{R}t'} \Delta_q^h \overline{w}_{\phi}}_{L^2}} dt'\\
	B_{3,q} &= \int_0^t \abs{\psca{ e^{\mathcal{R}t'} \Delta_q^h (\Rcal^h(v,\pa_zw))_{\phi}, e^{\mathcal{R}t'} \Delta_q^h \overline{w}_{\phi}}_{L^2}} dt'.
\end{align*}

As for the term $A_{1,q}$ in the proof of Lemma \ref{lem:uww}, we have the following estimate
\begin{align*}
	\abs{\psca{ e^{\mathcal{R}t} \Delta_q^h ( \Tcal^h_{v}\pa_z w)_{\phi}, e^{\mathcal{R}t} \Delta_q^h \overline{w}_{\phi}}_{L^2}} 
	\lesssim \sum_{|q'-q|\leq 4} e^{\mathcal{R}t} \Vert S_{q'-1}^h v_{\phi} \Vert_{L^2_xL^{\infty}_z} \Vert \Delta_{q'}^h \pa_z w_{\phi} \Vert_{L^{2}} \Vert \Delta_q^h e^{\mathcal{R}t} \overline{w}_{\phi} \Vert_{L^\infty_xL^{2}_z}.
\end{align*}
Lemma \ref{lem:tranreg}, Poincar\'e and Cauchy-Schwarz inequalities imply
\begin{multline*}
	\norm{S^h_{q'-1} v_{\phi}}_{L^2_xL^{\infty}_z} \leq \sum_{l \leq q'-2} 2^{l} \norm{\Delta^h_l u_\phi}_{L^2} 	\lesssim \sum_{l \leq q'-2} 2^{\frac{3l}4} \norm{\Delta^h_l u_\phi}_{L^2}^{\frac{1}2} 2^{\frac{l}4} \norm{\Delta^h_l \dd_z u_\phi}_{L^2}^{\frac{1}2}
	\\ 
	\lesssim \pare{\sum_{l \leq q'-2} 2^{\frac{3l}2} \norm{\Delta^h_l u_\phi}_{L^2}}^{\frac{1}2} \pare{\sum_{l \leq q'-2}  2^{\frac{l}2} \norm{\Delta^h_l \dd_z u_\phi}_{L^2}}^{\frac{1}2}
	\lesssim \norm{u_\phi}_{\Bcal^{\frac{3}2}}^{\frac{1}2} \norm{\dd_z u_\phi}_{\Bcal^{\frac{1}2}}^{\frac{1}2}.
\end{multline*}
Thus,
\begin{align}
	\label{eq:B1q}
	B_{1,q} &\lesssim \sum_{|q'-q|\leq 4} \int_0^t 2^{\frac{q}2} \norm{u_\phi}_{\Bcal^{\frac{3}2}}^{\frac{1}2} \norm{\dd_z u_\phi}_{\Bcal^{\frac{1}2}}^{\frac{1}2} \Vert \Delta_{q'}^h e^{\mathcal{R}t} \pa_z w_{\phi} \Vert_{L^{2}} \Vert \Delta_q^h e^{\mathcal{R}t} \overline{w}_{\phi} \Vert_{L^{2}}
	\\
	&\lesssim \sum_{|q'-q|\leq 4} \pare{\int_0^t \norm{u_\phi}_{\Bcal^{\frac{3}2}} \Vert \Delta_{q'}^h e^{\mathcal{R}t} \pa_z w_{\phi} \Vert_{L^{2}}^2}^{\frac{1}2} 2^{\frac{q}2} \pare{\int_0^t \norm{\dd_z u_\phi}_{\Bcal^{\frac{1}2}} \Vert \Delta_q^h e^{\mathcal{R}t} \overline{w}_{\phi} \Vert_{L^{2}}^2 }^{\frac{1}2} \notag
	\\
	&\lesssim d_q^2 2^{-2qs}\norm{u_\phi}_{L^\infty_t(\Bcal^{\frac{3}2})}^{\frac{1}2} \norm{e^{\mathcal{R}t} \pa_z w_{\phi}}_{\tilde{L}^2_t(\Bcal^s)} \norm{e^{\mathcal{R}t} \overline{w}_{\phi}}_{\tilde{L}^2_{t,f(t)}(\mathcal{B}^{s+\frac{1}{2}})} \notag
\end{align}
where 
\begin{align*}
	d_q^2 = d_q(\overline{w}_\phi) \sum_{|q'-q|\leq 4} d_{q'}(\dd_z w_\phi) 2^{(q-q')s}.
\end{align*}

For the second term on the right-hand side of \eqref{eq:B123q}, using similar calculations as in \eqref{eq:normB12}, we get
\begin{equation*}
	\Vert S_{q'-1}^h \pa_z w_{\phi} \Vert_{L^{\infty}_xL^2_z} \leq \norm{\dd_z w_\phi}_{\Bcal^{\frac{1}2}}.
\end{equation*}
Thus, Lemma \ref{lem:tranreg}, Young and Cauchy-Schwarz inequalities imply
\begin{align}
	\label{eq:B2q}
	B_{2,q} &\lesssim  \sum_{|q'-q|\leq 4} \int_0^t e^{2\mathcal{R}t'} \Vert S_{q'-1}^h \pa_z w_{\phi} \Vert_{L^{\infty}_xL^2_z} \Vert \Delta_{q'}^h v_{\phi} \Vert_{L^2_xL^{\infty}_z} \Vert \Delta_q^h \overline{w}_{\phi} \Vert_{L^2} dt' 
	\\
	&\lesssim \sum_{|q'-q|\leq 4} \int_0^t \Vert \pa_z w_{\phi} \Vert_{\mathcal{B}^{\frac{1}{2}}} 2^{q'} \Vert \Delta_{q'}^h e^{\mathcal{R}t'} u_{\phi} \Vert_{L^2} \Vert \Delta_q^h e^{\mathcal{R}t'} \overline{w}_{\phi} \Vert_{L^2} dt' \notag
	\\
	&\lesssim \sum_{|q-q'| \leq 4} 2^{q'} \left(\int_0^t f(t') e^{2\Rcal t'} \Vert \Delta_{q'}^h u_{\phi} \Vert_{L^2}^2 dt'\right)^{\f12} \left(\int_0^t f(t') e^{2\Rcal t'} \Vert \Delta_q^h \overline{w}_{\phi} \Vert_{L^2}^2 dt'\right)^{\f12} \notag
	\\
	&\lesssim d_q^2 2^{-2qs} \Vert e^{\mathcal{R} t} u_{\phi} \Vert_{\tilde{L}^2_{t,f(t)}(\mathcal{B}^{s+\frac{1}{2}})} \Vert e^{\mathcal{R} t} \overline{w}_{\phi} \Vert_{\tilde{L}^2_{t,f(t)}(\mathcal{B}^{s+\frac{1}{2}})}, \notag
\end{align}
where
\begin{align*}
	d_q^2 = d_q(\overline{w}_\phi) \sum_{|q'-q|\leq 4} d_{q'}(u_\phi) \, 2^{(q-q')(s-\frac{1}{2})}.
\end{align*}

Now, for the third term on the right-hand side of \eqref{eq:B123q}, using Lemma \ref{lem:tranreg}, we obtain
\begin{align}
	\label{eq:B3q}
	B_{3,q} &\lesssim \sum_{q'\geq q-3} \int_0^t e^{2\mathcal{R}t'} \Vert \Delta_{q'}^h v_{\phi} \Vert_{L^2_xL^{\infty}_z} \Vert \DD_{q'}^h \pa_z w_{\phi} \Vert_{L^2} \Vert \Delta_q^h \overline{w}_{\phi} \Vert_{L^\infty_xL^2_z} dt' 
	\\
	&\lesssim \sum_{q'\geq q-3} \int_0^t \pare{2^{q'} \Vert \Delta_{q'}^h e^{\mathcal{R}t'} u_{\phi} \Vert_{L^2}} \pare{2^{-\frac{q'}{2}} \Vert \pa_z w_{\phi} \Vert_{\mathcal{B}^\frac{1}{2}}} \pare{2^{\frac{q}{2}} \Vert \Delta_q^h e^{\mathcal{R}t'} \overline{w}_{\phi} \Vert_{L^2}} dt' \notag
	\\
	&\lesssim \sum_{q'\geq q-3} 2^{\frac{q+q'}2} \left(\int_0^t f(t') e^{2\Rcal t'} \Vert \Delta_{q'}^h u_{\phi} \Vert_{L^2}^2 dt'\right)^{\f12} \left(\int_0^t f(t') e^{2\Rcal t'} \Vert \Delta_q^h \overline{w}_{\phi} \Vert_{L^2}^2 dt'\right)^{\f12} \notag
	\\
	& \lesssim d_q^2 \, 2^{-2qs} \Vert e^{\mathcal{R}t} u_{\phi} \Vert_{\tilde{L}^2_{t,f(t)}(\mathcal{B}^{s+\frac{1}{2}})} \Vert e^{\mathcal{R}t} \overline{w}_{\phi} \Vert_{\tilde{L}^2_{t,f(t)}(\mathcal{B}^{s+\frac{1}{2}})}, \notag
\end{align}
where
\begin{align*}
	d_q^2 = d_q(\overline{w}_\phi) \sum_{q'\geq q-3}d_{q'}(u_\phi) \, 2^{(q-q')s}.
\end{align*}

The proof of Estimate \eqref{eq:vww} is completed.   \hfill $\square$

}

\section{Global wellposedness of the hydrostatic limit system} \label{se:hydrolim}

{\color{black}

The goal of this section is to prove Theorem \ref{th:hydrolimit}. We remark that the construction of a local smooth solution of the system \eqref{eq:hydrolimit} follows a standard parabolic regularization method, similar to the case of Prandtl system, which consists of adding an addition horizontal smoothing term of the type $\gamma \dd_x^2$ and then taking $\gamma \to 0$. The difficulty here relies on the presence of the unknown pressure term $\dd_xp$ in the first equation of \eqref{eq:hydrolimit}. However, as in \cite{CLT2019}, we can reformulate the problem by writing $v$ and $\dd_xp$ as functions of $u$ and $T$. From the Dirichlet boundary condition $(u,v)\vert_{z=0} = (u,v)\vert_{z=1} = 0$ and the incompressibility condition $\dd_x u + \dd_z v=0$, we get
\begin{align*}
	v(t,x,z) = \int_0^z \pa_z v(t,x,z')dz' = -\int_0^z \pa_x u(t,x,z') dz' .
\end{align*} 

For the pressure term, due to the Dirichlet boundary condition $(u,v,T)\vert_{z=0}=0$, we deduce from the incompressibility condition $\pa_xu+\pa_zv=0$ that
\begin{align}
	\pa_x\int_0^1u(t,x,z)\, dz' = -\int_0^1 \dd_z v(t,x,z')\, dz' = v(t,x,1) - v(t,x,0) = 0.
\end{align}
Integrating the equation $\pa_z p = T$ with respect $z$, we obtain 
\begin{align} \label{p}
	p(t,x,z) = p(t,x,0) + \int_{0}^z T(t,x,z') dz'.
\end{align}
Next, differentiating \eqref{p} with respect to $x$ and using the first equation of the system \eqref{eq:hydrolimit}, we get
\begin{align*}
	\dd_x p(t,x,0) &= - \int_{0}^z \dd_x T(t,x,z') dz' + \dd_x p(t,x,z)\\
	&= - \int_{0}^z \dd_x T(t,x,z') dz' - \pare{\p_tu+u\p_x u+v\p_zu-\p_z^2u} (t,x,z)
\end{align*}
Now, we set $c(t) = \int_0^1 u(t,x,z)dz$. Integrating the above equation with respect to $z \in [0,1]$ and performing integration by parts lead to
\begin{align*}
	\dd_x p(t,x,0) = - \int_0^1 \int_{0}^z \dd_x T(t,x,z') dz' dz + \pa_zu(t,x,1) - \pa_zu(t,x,0) - \dot{c}(t)- \pa_x \int_0^1 u^2(t,x,z)dz
\end{align*}
that yields 
\begin{multline*}
	\dd_x p(t,x,z) = \int_{0}^z \pa_x T(t,x,z')dz' - \int_0^1 \int_{0}^z \dd_x T(t,x,z') dz' dz \\ + \pa_zu(t,x,1) - \pa_zu(t,x,0) - \dot{c}(t)- \pa_x \int_0^1 u^2(t,x,z)dz.
\end{multline*}

Let  $\phi: \RR_+\times \RR \to \RR_+$ such that 
\begin{displaymath}
	\phi(t,\xi) = (a-\lambda \rho(t)) \abs{\xi}.
\end{displaymath}
Here $\lambda > 0$ and $\rho$ is a function satisfying \eqref{eq:AnTheta}, both will be determined later. Applying the operator defined in \eqref{eq:AnPhi} to the system \eqref{eq:hydrolimit} and taking into account \eqref{eq:derivativePhi}, we obtain
\begin{equation} \label{S1eq8}
	\left\{ \;
	\begin{aligned}
		&\displaystyle \p_tu_{\phi}+ \lambda \dot{\rho}(t)|D_x| u_{\phi}+(u\p_x u)_{\phi} + (v\p_zu)_{\phi}-\p_z^2u_{\phi} + \p_xp_{\phi}=0 && \mbox{in } \ \cS\times ]0,\infty[,\\
		&\displaystyle\p_zp_{\phi}= T_{\phi}\\
		&\p_t T_{\phi}+\lambda \dot{\rho}(t)|D_x| T_{\phi}+(u\pa_x T)_{\phi}+(v\pa_z T)_{\phi}-\Delta T_{\phi}=0\\
		&\p_xu_{\phi}+\p_zv_{\phi}=0,\\
		&(u_\phi,T_\phi)|_{t=0}=(e^{a\abs{D_x}}u_0,e^{a\abs{D_x}}T_0)\\
		&(u_\phi, v_\phi,T_\phi)|_{z=0}=(u_\phi, v_\phi,T_\phi)|_{z=1}=0.
	\end{aligned}
	\right.
\end{equation}
where $|D_x|$ denotes the Fourier multiplier of symbol $|\xi|$. In what follows, we recall that for the sake of the simplicity, we use ``$C$'' to denote a generic positive constant which can change from line to line and $(d_q)_{q\in\ZZ}$(resp. $(d_q(t))_{q\in\ZZ}$) a generic element of $ \ell^1(\ZZ)$ such that $\sum_{q\in\ZZ} d_q = 1$ (resp. $\sum_{q\in\ZZ} d_q(t) = 1$).

We will now perform \emph{a priori} estimates needed to prove Theorem \ref{th:hydrolimit}. We remark that the continuity with respect to time variable of the solutions in $\Bcal^s$ can be proved by a similar argument as in \cite{Paicu2005}. We will only focus on estimates in analytic norms. Applying the dyadic operator $\DD^h_q$ to the system \eqref{S1eq8}, then taking the $L^2(\Scal)$-scalar product of the first and the third equations of the obtained system with $\Delta_q^h u_{\phi}$ and $\Delta_q^h T_{\phi}$ respectively, we get
\begin{multline} \label{S4eq5}
	\frac{1}{2} \frac{d}{dt} \Vert \Delta_q^h u_{\phi}(t) \Vert_{L^2}^2 + \lambda \dot{\rho}(t) \norm{|D_x|^{\frac12} \Delta_q^h u_{\phi}}_{L^2}^2 + \Vert \Delta_q^h \pa_z u_{\phi}(t) \Vert_{L^2}^2 \\
	= -\psca{\Delta_q^h (u\p_x u)_{\phi}, \Delta_q^h u_{\phi})}_{L^2} - \psca{\Delta_q^h (v\p_z u)_{\phi}, \Delta_q^h u_{\phi}}_{L^2} - \psca{\Delta_q^h \pa_x p_{\phi}, \Delta_q^h u_{\phi}}_{L^2},
\end{multline}
and
\begin{multline} \label{S4eq6}
	\frac{1}{2} \frac{d}{dt} \Vert \Delta_q^h T_{\phi}(t) \Vert_{L^2}^2 + \lambda \dot{\rho}(t) \norm{|D_x|^{\frac12} \Delta_q^h T_{\phi}}_{L^2}^2 + \Vert \Delta_q^h \nabla T_{\phi}(t) \Vert_{L^2}^2 \\
	= -\psca{\Delta_q^h (u\p_x T)_{\phi}, \Delta_q^h T_{\phi}}_{L^2} - \psca{\Delta_q^h (v\p_z T)_{\phi}, \Delta_q^h T_{\phi}}_{L^2}. \hspace{1cm}
\end{multline}
Let $\Rcal > 0$. Multiplying \eqref{S4eq5} and \eqref{S4eq6} with $e^{2\Rcal t}$ and taking remarking that 
$$e^{2\Rcal t} \frac{d}{dt} f(t) = \frac{d}{dt} \pare{f(t) e^{2\Rcal t}} - 2\mathcal{R} f(t) e^{2\Rcal t},$$ we have
\begin{multline*}
	\frac{1}{2} \frac{d}{dt} \Vert e^{\Rcal t} \Delta_q^h u_{\phi}(t) \Vert_{L^2}^2 - \Rcal \Vert e^{\Rcal t} \Delta_q^h u_{\phi}(t) \Vert_{L^2}^2 + \lambda \dot{\rho}(t) \norm{e^{\Rcal t} |D_x|^{\frac12} \Delta_q^h u_{\phi}}_{L^2}^2 + \Vert e^{\Rcal t} \Delta_q^h \pa_z u_{\phi}(t) \Vert_{L^2}^2 \\
	= -e^{2\Rcal t} \psca{\Delta_q^h (u\p_x u)_{\phi}, \Delta_q^h u_{\phi})}_{L^2} - e^{2\Rcal t}\psca{\Delta_q^h (v\p_z u)_{\phi}, \Delta_q^h u_{\phi}}_{L^2} - e^{2\Rcal t} \psca{\Delta_q^h \pa_x p_{\phi}, \Delta_q^h u_{\phi}}_{L^2},
\end{multline*}
and
\begin{multline*}
	\frac{1}{2} \frac{d}{dt} \Vert e^{\Rcal t} \Delta_q^h T_{\phi}(t) \Vert_{L^2}^2 - \Rcal \Vert e^{\Rcal t} \Delta_q^h T_{\phi}(t) \Vert_{L^2}^2 + \lambda \dot{\rho}(t) \norm{e^{\Rcal t} |D_x|^{\frac12} \Delta_q^h T_{\phi}}_{L^2}^2 + \Vert e^{\Rcal t} \Delta_q^h \nabla T_{\phi}(t) \Vert_{L^2}^2 \\
	= - e^{2\Rcal t}\psca{\Delta_q^h (u\p_x T)_{\phi}, \Delta_q^h T_{\phi}}_{L^2} - e^{2\Rcal t}\psca{\Delta_q^h (v\p_z T)_{\phi}, \Delta_q^h T_{\phi}}_{L^2}. \hspace{1cm}
\end{multline*}
For $0 \leq \Rcal \leq \frac{1}{2\Kcal}$, we apply Poincar\'e inequality, then we integrate the above identities with respect to the time variable and get
\begin{multline} \label{S4eq5bis}
	\norm{e^{\Rcal t} \Delta_q^h u_{\phi}(t)}_{L^\infty_t (L^2)}^2 + 2\lambda \int_0^t \dot{\rho}(t') \norm{e^{\Rcal t'} |D_x|^{\frac12} \Delta_q^h u_{\phi}}_{L^2}^2 dt' + \norm{e^{\Rcal t} \Delta_q^h \pa_z u_{\phi}(t)}_{L^2_t(L^2)}^2 \\
	\leq \norm{\Delta_q^h u_{\phi}(0)}_{L^2}^2 + D_1 + D_2 + D_3, \hspace{2cm}
\end{multline}
and
\begin{multline} \label{S4eq6bis}
	\norm{e^{\Rcal t} \Delta_q^h T_{\phi}(t)}_{L^\infty_t(L^2)}^2 + 2\lambda \int_0^t \dot{\rho}(t') \norm{e^{\Rcal t} |D_x|^{\frac12} \Delta_q^h T_{\phi}}_{L^2}^2 dt' + \norm{e^{\Rcal t} \Delta_q^h \nabla T_{\phi}(t)}_{L^2_t(L^2)}^2\\
	\leq \norm{\Delta_q^h T_{\phi}(0)}_{L^2}^2 + D_4 + D_5. \hspace{2.5cm}
\end{multline}

From now on, we set 
\begin{equation}
	\label{eq:defrho}
	\dot{\rho}(t) = \norm{\dd_z u_\phi(t)}_{\mathcal{B}^{\frac{1}{2}}} + \norm{\dd_z T_\phi(t)}_{\mathcal{B}^{\frac{1}{2}}}.
\end{equation}
Lemmas \ref{lem:uww} and \ref{lem:vww} yield
\begin{align*}
	\abs{D_1} = 2\abs{\int_0^t \psca{e^{\Rcal t'}\Delta_q^h (u\p_x u)_{\phi}, e^{\Rcal t'}\Delta_q^h u_{\phi})} dt'} \leq C d_q^2 2^{-2qs} \Vert e^{\mathcal{R}t} u_{\phi} \Vert_{\tilde{L}^2_{t,\dot{\rho}(t)}(\mathcal{B}^{s+\frac{1}{2}})}^2,
\end{align*}
\begin{align*}
	\abs{D_2} &= 2\abs{\int_0^t \psca{e^{\Rcal t'}\Delta_q^h (v\p_z u)_{\phi}, e^{\Rcal t'}\Delta_q^h u_{\phi}} dt'}
	\\ 
	&\leq C d_q^2 2^{-2qs} \pare{\norm{e^{\mathcal{R}t} u_{\phi}}_{\tilde{L}^2_{t,\dot{\rho}(t)}(\mathcal{B}^{s+\frac{1}{2}})}^2 + \norm{u_\phi}_{L^\infty_t(\mathcal{B}^{\frac{3}{2}})}^{\frac{1}2} \norm{e^{\mathcal{R}t} \dd_z u_{\phi}}_{\tilde{L}^2_t(\mathcal{B}^s)} \norm{e^{\mathcal{R}t} u_{\phi}}_{\tilde{L}^2_{t,\dot{\rho}(t)}(\mathcal{B}^{s+\frac{1}{2}})}} 
	\\
	&\leq C d_q^2 2^{-2qs} \Vert e^{\mathcal{R}t} u_{\phi} \Vert_{\tilde{L}^2_{t,\dot{\rho}(t)}(\mathcal{B}^{s+\frac{1}{2}})}^2 +  d_q^2 2^{-2qs} \norm{u_\phi}_{L^\infty_t(\mathcal{B}^{\frac{3}{2}})} \norm{e^{\mathcal{R}t} \dd_z u_{\phi}}_{\tilde{L}^2_t(\mathcal{B}^s)}^2,
\end{align*}
\begin{align*}
	\abs{D_4} &= 2\abs{\int_0^t \psca{e^{\Rcal t'}\Delta_q^h (u\p_x T)_{\phi}, e^{\Rcal t'}\Delta_q^h T_{\phi}} dt'}
	\\ 
	&\leq C d_q^2 2^{-2qs} \pare{\norm{e^{\Rcal t} u_{\phi}}_{\tilde{L}^2_{t,\dot{\rho}(t)}(\mathcal{B}^{s+\frac{1}{2}})} + \norm{e^{\Rcal t} T_{\phi}}_{\tilde{L}^2_{t,\dot{\rho}(t)}(\mathcal{B}^{s+\frac{1}{2}})}} \norm{e^{\Rcal t} T_{\phi}}_{\tilde{L}^2_{t,\dot{\rho}(t)}(\mathcal{B}^{s+\frac{1}{2}})}\\
	&\leq C d_q^2 2^{-2qs} \pare{\norm{e^{\Rcal t} u_{\phi}}_{\tilde{L}^2_{t,\dot{\rho}(t)}(\mathcal{B}^{s+\frac{1}{2}})}^2 + \norm{e^{\Rcal t} T_{\phi}}_{\tilde{L}^2_{t,\dot{\rho}(t)}(\mathcal{B}^{s+\frac{1}{2}})}^2}
\end{align*}
and
\begin{align*}
	\abs{D_5} &= 2\abs{\int_0^t \psca{e^{\Rcal t'}\Delta_q^h (v\p_z T)_{\phi}, e^{\Rcal t'}\Delta_q^h T_{\phi}} dt'}
	\\ 
	&\leq C d_q^2 2^{-2qs} \norm{e^{\mathcal{R}t} T_{\phi}}_{\tilde{L}^2_{t,\dot{\rho}(t)}(\mathcal{B}^{s+\frac{1}{2}})} \pare{\norm{e^{\mathcal{R}t} u_{\phi}}_{\tilde{L}^2_{t,\dot{\rho}(t)}(\mathcal{B}^{s+\frac{1}{2}})} + \norm{u_\phi}_{L^\infty_t(\mathcal{B}^{\frac{3}{2}})}^{\frac{1}2} \norm{e^{\mathcal{R}t} \dd_z T_{\phi}}_{\tilde{L}^2_t(\mathcal{B}^s)}}
	\\
	&\leq C d_q^2 2^{-2qs} \pare{\norm{e^{\Rcal t} u_{\phi}}_{\tilde{L}^2_{t,\dot{\rho}(t)}(\mathcal{B}^{s+\frac{1}{2}})}^2 + \norm{e^{\Rcal t} T_{\phi}}_{\tilde{L}^2_{t,\dot{\rho}(t)}(\mathcal{B}^{s+\frac{1}{2}})}^2} +  d_q^2 2^{-2qs} \norm{u_\phi}_{L^\infty_t(\mathcal{B}^{\frac{3}{2}})} \norm{e^{\mathcal{R}t} \dd_z u_{\phi}}_{\tilde{L}^2_t(\mathcal{B}^s)}^2.
\end{align*}
Concerning the pressure term, using the Dirichlet boundary condition $(u,v,T)\vert_{z=0} = 0$, the incompressibility condition $\pa_x u +\pa_z v =0$, the relation $\pa_z p = T$ and Poincar\'e inequality, we can perform integrations by parts and get 
\begin{multline}
	\label{eq:pressure}
	\abs{\psca{\Delta_q^h \pa_x p_{\phi}, \Delta_q^h u_{\phi}}} = \abs{ \psca{\Delta_q^h  p_{\phi}, \Delta_q^h \pa_x u_{\phi}}} = \abs{\psca{\Delta_q^h  p_{\phi}, \Delta_q^h \pa_z v_{\phi}}} = \abs{ \psca{\Delta_q^h  \pa_z p_{\phi}, \Delta_q^h  v_{\phi}}} 
	\\
	= \abs{ \psca{\Delta_q^h   T_{\phi}, \Delta_q^h  v_{\phi}}} = \abs{\psca{\Delta_q^h  T_{\phi}, \Delta_q^h  \int_0^z \pa_x u_{\phi} dz'}} = \abs{\psca{\Delta_q^h  \pa_x T_{\phi}, \Delta_q^h  \int_0^z u_{\phi} dz'}}, 
\end{multline}
and so we can have the following bound
\begin{equation*}
	\abs{\psca{\Delta_q^h \pa_x p_{\phi}, \Delta_q^h u_{\phi}}} \leq \Vert \Delta_q^h \pa_x T_{\phi} \Vert_{L^2} \Vert \Delta_q^h u_{\phi} \Vert_{L^2}  \leq C\Vert \Delta_q^h \pa_x T_{\phi} \Vert_{L^2}^2 + \frac{1}4 \Vert \Delta_q^h \dd_z u_{\phi} \Vert_{L^2}^2.
\end{equation*}
Thus,
\begin{align*}
	\abs{D_3} = 2\abs{\int_0^t \psca{e^{\Rcal t'} \Delta_q^h \pa_x p_{\phi}, e^{\Rcal t'} \Delta_q^h u_{\phi}} dt'} \leq C d_q^2 2^{-2qs} \norm{e^{\Rcal t} \dd_x T_\phi}_{\tilde{L}^2_t(\mathcal{B}^{s})}^2 + \frac{1}{4} \norm{e^{\Rcal t} \Delta_q^h \pa_z u_{\phi}(t)}_{L^2_tL^2}^2.
\end{align*}

Now, we recall that for a positive sequence $(a_1, \ldots, a_n)$, $n\in \NN^*$, we have
\begin{equation}
	\label{squareroot}
	\frac{1}{\sqrt{n}} \sum_{j=1}^n a_n \leq \sqrt{\sum_{j=1}^n a_n^2} \leq \sum_{j=1}^n a_n.
\end{equation}
Inequality \eqref{squareroot} allows us to take the square root of each terms on the two sides of \eqref{S4eq5bis} and \eqref{S4eq6bis} (with a cost of a constant multiplier on the right-hand side, which will be included in the generic constant $C$). Summing the two obtained inequalities, we get
\begin{multline*}
	\norm{e^{\Rcal t} \Delta_q^h u_{\phi}(t)}_{L^\infty_t (L^2)} + \sqrt{2\lambda} \pare{\int_0^t \dot{\rho}(t') \norm{e^{\Rcal t'} |D_x|^{\frac12} \Delta_q^h u_{\phi}}_{L^2}^2 dt'}^{\frac{1}2} + \frac{1}2 \norm{e^{\Rcal t} \Delta_q^h \pa_z u_{\phi}(t)}_{L^2_t(L^2)} \\
	\leq \norm{e^{a\abs{D_x}} \Delta_q^h u_0}_{L^2} + C d_q 2^{-qs} \Vert e^{\mathcal{R}t} u_{\phi} \Vert_{\tilde{L}^2_{t,\dot{\rho}(t)}(\mathcal{B}^{s+\frac{1}{2}})} + C d_q 2^{-qs} \norm{e^{\Rcal t} \dd_x T_\phi}_{\tilde{L}^2_t(\mathcal{B}^{s})}
	\\
	+ d_q 2^{-qs} \norm{u_\phi}_{L^\infty_t(\mathcal{B}^{\frac{3}{2}})}^{\frac{1}2} \norm{e^{\mathcal{R}t} \dd_z u_{\phi}}_{\tilde{L}^2_t(\mathcal{B}^s)},
\end{multline*}
and
\begin{multline*}
	\norm{e^{\Rcal t} \Delta_q^h T_{\phi}(t)}_{L^\infty_t(L^2)} + \sqrt{2\lambda} \pare{\int_0^t \dot{\rho}(t') \norm{e^{\Rcal t} |D_x|^{\frac12} \Delta_q^h T_{\phi}}_{L^2}^2 dt'}^{\frac{1}2} + \norm{e^{\Rcal t} \Delta_q^h \nabla T_{\phi}(t)}_{L^2_t(L^2)}\\
	\leq \norm{e^{a\abs{D_x}} \Delta_q^h T_0}_{L^2} + C d_q 2^{-qs} \Vert e^{\mathcal{R}t} u_{\phi} \Vert_{\tilde{L}^2_{t,\dot{\rho}(t)}(\mathcal{B}^{s+\frac{1}{2}})} + C d_q 2^{-qs} \Vert e^{\Rcal t} T_{\phi} \Vert_{\tilde{L}^2_{t,\dot{\rho}(t)}(\mathcal{B}^{s+\frac{1}{2}})}
	\\
	+ d_q 2^{-qs} \norm{u_\phi}_{L^\infty_t(\mathcal{B}^{\frac{3}{2}})}^{\frac{1}2} \norm{e^{\mathcal{R}t} \dd_z T_{\phi}}_{\tilde{L}^2_t(\mathcal{B}^s)}.
\end{multline*}
Multiplying the above inequalities by $2^{qs}$ and then summing with respect to $q \in \ZZ$, we obtain 
\begin{multline} \label{S4eq5ter}
	\norm{e^{\Rcal t} u_{\phi}}_{\tilde{L}^\infty_t(\Bcal^s)} + \sqrt{2\lambda} \norm{e^{\Rcal t} u_{\phi}}_{\tilde{L}^2_{t,\dot{\rho}(t)} (\Bcal^{s + \frac12})} + \frac{1}2 \norm{e^{\Rcal t} \pa_z u_{\phi}}_{\tilde{L}^2_t(\Bcal^s)} 
	\\
	\leq \norm{e^{a\abs{D_x}} u_0}_{\Bcal^s} + C \Vert e^{\mathcal{R}t} u_{\phi} \Vert_{\tilde{L}^2_{t,\dot{\rho}(t)}(\mathcal{B}^{s+\frac{1}{2}})} + C \norm{e^{\Rcal t} \dd_x T_\phi}_{\tilde{L}^2_t(\mathcal{B}^{s})} + \norm{u_\phi}_{L^\infty_t(\mathcal{B}^{\frac{3}{2}})}^{\frac{1}2} \norm{e^{\mathcal{R}t} \dd_z u_{\phi}}_{\tilde{L}^2_t(\mathcal{B}^s)},
\end{multline}
and
\begin{multline} \label{S4eq6ter}
	\norm{e^{\Rcal t} T_{\phi}}_{\tilde{L}^\infty_t(\Bcal^s)} + \sqrt{2\lambda} \norm{e^{\Rcal t} T_{\phi}}_{\tilde{L}^2_{t,\dot{\rho}(t)} (\Bcal^{s + \frac12})} + \norm{e^{\Rcal t} \nabla T_{\phi}}_{\tilde{L}^2_t(\Bcal^s)}
	\\
	\leq \norm{e^{a\abs{D_x}} T_0}_{\Bcal^s} + C \Vert e^{\mathcal{R}t} u_{\phi} \Vert_{\tilde{L}^2_{t,\dot{\rho}(t)}(\mathcal{B}^{s+\frac{1}{2}})} + C \Vert e^{\Rcal t} T_{\phi} \Vert_{\tilde{L}^2_{t,\dot{\rho}(t)}(\mathcal{B}^{s+\frac{1}{2}})} + \norm{u_\phi}_{L^\infty_t(\mathcal{B}^{\frac{3}{2}})}^{\frac{1}2} \norm{e^{\mathcal{R}t} \dd_z T_{\phi}}_{\tilde{L}^2_t(\mathcal{B}^s)}.
\end{multline}
Without loss of generality, we can suppose that $C \geq 2$. Multiplying Inequality \eqref{S4eq6ter} by $2C$ and then adding the obtained result to Inequality \eqref{S4eq5ter}, we get
\begin{multline*}
	%\label{eq:u+T}
	\norm{e^{\Rcal t} \pare{u_{\phi},T_{\phi}}}_{\tilde{L}^\infty_t(\Bcal^s)} + \sqrt{2\lambda} \norm{e^{\Rcal t} \pare{u_{\phi},T_{\phi}}}_{\tilde{L}^2_{t,\dot{\rho}(t)} (\Bcal^{s + \frac12})} + \frac{1}{2} \norm{e^{\Rcal t} \pa_z u_{\phi}}_{\tilde{L}^2_t(\Bcal^s)} + \norm{e^{\Rcal t} \nabla T_{\phi}}_{\tilde{L}^2_t(\Bcal^s)}
	\\
	\leq 2C \norm{e^{a\abs{D_x}} (u_0,T_0)}_{\mathcal{B}^{s}} + 3C^2 \norm{e^{\mathcal{R}t} \pare{u_{\phi},T_\phi}}_{\tilde{L}^2_{t,\dot{\rho}(t)}(\mathcal{B}^{s+\frac{1}{2}})} + 2C \norm{u_\phi}_{L^\infty_t(\mathcal{B}^{\frac{3}{2}})}^{\frac{1}2} \norm{e^{\mathcal{R}t} (\dd_z u_{\phi}, \dd_z T_{\phi})}_{\tilde{L}^2_t(\mathcal{B}^s)}.
\end{multline*}
From now on, we will fix $\lambda$ and $\Rcal$ such that 
\begin{displaymath}
	0 < \Rcal < \frac{1}{2\Kcal} \quad\mbox{and}\quad \sqrt{\lambda} \geq 9C^2.
\end{displaymath}
We obtain from the above inequality that
\begin{multline}
	\label{eq:u+T}
	\norm{e^{\Rcal t} \pare{u_{\phi},T_{\phi}}}_{\tilde{L}^\infty_t(\Bcal^s)} + \sqrt{\lambda} \norm{e^{\Rcal t} \pare{u_{\phi},T_{\phi}}}_{\tilde{L}^2_{t,\dot{\rho}(t)} (\Bcal^{s + \frac12})} + \frac{1}{2} \norm{e^{\Rcal t} \pa_z u_{\phi}}_{\tilde{L}^2_t(\Bcal^s)} + \norm{e^{\Rcal t} \nabla T_{\phi}}_{\tilde{L}^2_t(\Bcal^s)}\\
	\leq  2C \norm{e^{a\abs{D_x}} (u_0,T_0)}_{\mathcal{B}^{s}} + 2C \norm{u_\phi}_{L^\infty_t(\mathcal{B}^{\frac{3}{2}})}^{\frac{1}2} \norm{e^{\mathcal{R}t} (\dd_z u_{\phi}, \dd_z T_{\phi})}_{\tilde{L}^2_t(\mathcal{B}^s)}.
\end{multline}

Let
\begin{equation*} 
	t^\star = \sup\set{t>0\ : \ \norm{u_\phi(t)}_{\mathcal{B}^{\frac{3}{2}}} \leq \frac{1}{16C^4}, \ \mbox{and} \ \rho(t) \leq  \frac{a}{3\lambda}}.
\end{equation*}
For small initial data such that 
\begin{equation*}
	\left\{
	\begin{aligned}
		&\Vert e^{a|D_x|} u_0 \Vert_{\mathcal{B}^{\frac{1}2}} + \Vert e^{a|D_x|} T_0 \Vert_{\mathcal{B}^{\frac{1}2}} < \frac{a\sqrt{2\Rcal}}{16C\lambda}\\
		&\norm{e^{a|D_x|} u_0}_{\mathcal{B}^{\frac{3}2}} + \norm{e^{a|D_x|} T_0}_{\mathcal{B}^{\frac{3}2}} \leq \frac{1}{64C^5},
	\end{aligned}
	\right.
\end{equation*}
the continuity with respect to the time variable in $\Bcal^{\frac{3}2}$ and the fact that $\rho(0) = 0$ imply that $t^\star > 0$. 
For $s=\frac{3}2$ and for any $0 < t < t^\star$, we have
\begin{equation*}
	\norm{u_\phi(t)}_{\mathcal{B}^{\frac{3}{2}}} \leq \norm{e^{\Rcal t} \pare{u_{\phi},T_{\phi}}}_{\tilde{L}^\infty_t(\Bcal^{\frac{3}{2}})} < 2C \norm{e^{a\abs{D_x}} (u_0,T_0)}_{\mathcal{B}^{\frac{3}{2}}} \leq 2C \cdot \frac{1}{64C^5} = \frac{1}{32C^4}.
\end{equation*}
For $s=\frac{1}2$ and for any $0 < t < t^\star$, we deduce from \eqref{eq:u+T} that
\begin{align*}
	\rho (t) &= \int_0^t \Vert \pa_z (u_{\phi},T_\phi)(t') \Vert_{\mathcal{B}^\f12} dt' \leq \int_0^t  e^{-\mathcal{R}t'} \Vert e^{\mathcal{R}t'} \pa_z (u_{\phi},T_\phi)(t') \Vert_{\mathcal{B}^\f12} dt' \\
	&\leq \left( \int_0^t  e^{-2\mathcal{R}t'} dt' \right)^{\f12} \left( \int_0^t \Vert e^{\mathcal{R}t'} \pa_z (u_{\phi},T_\phi)(t') \Vert_{\mathcal{B}^\f12}^2  dt' \right)^{\f12} \\
	&\leq \frac{1}{\sqrt{2\Rcal}} \norm{e^{\Rcal t} \pa_z (u_{\phi},T_\phi)}_{\tilde{L}^2_t(\Bcal^{\frac{1}{2})}} \\
	&\leq \frac{4C}{\sqrt{2\Rcal}} \left( \Vert e^{a|D_x|} u_0 \Vert_{\mathcal{B}^\f12} +  \Vert e^{a|D_x|}  T_{0} \Vert_{\mathcal{B}^\f12}\right) < \frac{a}{4\lambda}.
\end{align*}
We deduce that $t^\star = +\infty$ and that \eqref{eq:hydrolimitU} is verified for any $t \in \RR_+$.

\medskip

In order to end the proof of Theorem \ref{th:hydrolimit}, we need to prove Inequality \eqref{eq:hydrolimitdtU}. Applying the operator $\phi$ and then $\Delta_q^h$ to the first equation of \eqref{eq:hydrolimit} and taking the $L^2$ inner product of resulting equation with $\Delta_q^h (\pa_tu)_\phi$ yield 
\begin{align*}
	\Vert \Delta_q^h (\pa_tu)_\phi \Vert_{L^2}^2 &= \psca{\Delta_q^h \pa_z^2 u_\phi,\Delta_q^h (\pa_tu)_\phi }_{L^2}  -  \psca{\Delta_q^h (u\pa_xu)_\phi,\Delta_q^h (\pa_tu)_\phi }_{L^2} \\ & - \psca{ \Delta_q^h (v\pa_zu)_\phi,\Delta_q^h (\pa_tu)_\phi }_{L^2} - \psca{ \Delta_q^h \pa_x p_\phi , \Delta_q^h (\pa_tu)_\phi }_{L^2}.
\end{align*}
An integration by parts gives 
\begin{align*}
	\psca{ \Delta_q^h \pa_z^2 u_\phi,\Delta_q^h (\pa_tu)_\phi }_{L^2} &= \psca{ \Delta_q^h \pa_z^2 u_\phi, \pa_t \Delta_q^hu_\phi + \lambda \dot{\rho}(t) \abs{D_x} \Delta_q^hu_\phi }_{L^2}
	\\
	&= -\pare{\f12 \frac{d}{dt} \norm{\Delta_q^h \pa_z u_\phi}_{L^2}^2 + \lambda \dot{\rho}(t) \norm{\abs{D_x}^{\frac{1}2} \Delta_q^h \pa_z u_\phi}_{L^2}^2},
\end{align*}
and thus, 
\begin{align*}
	\norm{\Delta_q^h (\pa_tu)_\phi}_{L^2}^2 + \f12 \frac{d}{dt} \norm{\Delta_q^h \pa_z u_\phi}_{L^2}^2 + \lambda \dot{\rho}(t) \norm{\abs{D_x}^{\frac{1}2} \Delta_q^h \pa_z u_\phi}_{L^2}^2 \leq  I_1 + I_2 + I_3,
\end{align*}
where
\begin{align*}
	I_1 &= \abs{\psca{ \Delta_q^h (u\pa_xu)_\phi,\Delta_q^h (\pa_tu)_\phi }_{L^2}}\\
	I_2 &=  \abs{\psca{  \Delta_q^h (v\pa_zu)_\phi,\Delta_q^h (\pa_tu)_\phi }_{L^2}}\\
	I_3 &= \abs{\psca{  \Delta_q^h \pa_x p_\phi , \Delta_q^h (\pa_tu)_\phi }_{L^2}}.
\end{align*}

For $I_1$ and $I_2$, we simply have 
\begin{align*}
	I_1 &= \abs{\psca{  \Delta_q^h (u\pa_xu)_\phi,\Delta_q^h (\pa_tu)_\phi }_{L^2}} \leq  C \Vert \Delta_q^h (u\pa_xu)_\phi \Vert_{L^2}^2  + \frac{1}{6} \Vert \Delta_q^h (\pa_tu)_\phi \Vert_{L^2}^2 \\
	I_2 &= \abs{\psca{  \Delta_q^h (v\pa_zu)_\phi,\Delta_q^h (\pa_tu)_\phi }_{L^2}} \leq C\Vert \Delta_q^h (v\pa_zu)_\phi \Vert_{L^2}^2  + \frac{1}{6} \Vert \Delta_q^h (\pa_tu)_\phi \Vert_{L^2}^2.
\end{align*}
Now, using similar calculations as in \eqref{eq:pressure}, we find 
\begin{align*}
	I_3 = \abs{\psca{  \Delta_q^h \pa_x p_\phi , \Delta_q^h (\pa_tu)_\phi }_{L^2}} \leq \norm{\Delta_q^h \pa_xT_\phi}_{L^2} \norm{\Delta_q^h (\pa_tu)_\phi}_{L^2} \leq C \norm{\Delta_q^h \pa_xT_\phi}_{L^2}^2 + \frac{1}6  \norm{\Delta_q^h (\pa_tu)_\phi}_{L^2}^2. 
\end{align*}
Then, we deduce that
\begin{multline*}
	\norm{\Delta_q^h (\pa_tu)_\phi}_{L^2}^2 + \frac{d}{dt} \norm{\Delta_q^h \pa_z u_\phi}_{L^2}^2 + \lambda \dot{\rho}(t) \norm{\abs{D_x}^{\frac{1}2} \Delta_q^h \pa_z u_\phi}_{L^2}^2
	\\ 
	\leq C\left(\Vert \Delta_q^h (u\pa_xu)_\phi \Vert_{L^2}^2  + \Vert \Delta_q^h (v\pa_zu)_\phi \Vert_{L^2}^2 + \Vert \Delta_q^h \pa_xT_\phi \Vert_{L^2}^2 \right).
\end{multline*}
Multiplying the above inequality by $e^{2\mathcal{R}t}$, we get 
\begin{multline*}
	\norm{e^{\mathcal{R}t} \Delta_q^h (\pa_tu)_\phi}_{L^2}^2 + \frac{d}{dt} \norm{e^{\mathcal{R}t} \Delta_q^h \pa_z u_\phi}_{L^2}^2 - 2\Rcal \norm{e^{\mathcal{R}t} \Delta_q^h \pa_z u_\phi}_{L^2}^2 + \lambda \dot{\rho}(t) \norm{e^{\mathcal{R}t} \abs{D_x}^{\frac{1}2} \Delta_q^h \pa_z u_\phi}_{L^2}^2
	\\ 
	\leq C\left(\Vert \Delta_q^h (u\pa_xu)_\phi \Vert_{L^2}^2  + \Vert \Delta_q^h (v\pa_zu)_\phi \Vert_{L^2}^2 + \Vert \Delta_q^h \pa_xT_\phi \Vert_{L^2}^2 \right).
\end{multline*}
Integrating over $[0,t]$, we obtain 
\begin{multline}
	\label{dtDelta}
	\Vert e^{\mathcal{R}t}\Delta_q^h (\pa_tu)_\phi \Vert_{L^2_t(L^2)}^2 + \Vert e^{\mathcal{R}t} \Delta_q^h \pa_z u_\phi \Vert_{L^\infty_t(L^2)}^2 - 2\Rcal \Vert e^{\mathcal{R}t} \Delta_q^h \pa_z u_\phi \Vert_{L^2_t(L^2)}^2 + \lambda \int_0^t \dot{\rho} \norm{e^{\mathcal{R}t'} \abs{D_x}^{\frac{1}2} \Delta_q^h \pa_z u_\phi}_{L^2}^2 dt'
	\\ 
	\leq \Vert \Delta_q^h \pa_ze^{a|D_x|} u_0 \Vert_{L^2}^2 + C \Big(\Vert e^{\mathcal{R}t}\Delta_q^h (u\pa_xu)_\phi \Vert_{L^2_t(L^2)}^2+\Vert e^{\mathcal{R}t}\Delta_q^h (v\pa_zu)_\phi \Vert_{L^2_t(L^2)}^2 + \Vert e^{\mathcal{R}t}\Delta_q^h \pa_xT_\phi \Vert_{L^2_t(L^2)}^2 \Big).
\end{multline}
We recall that Inequality \eqref{squareroot} allows to take the square root of each terms in the above inequality, with the price of a constant multiplier that will be included in the generic constant $C$. So, taking the square root of each terms of \eqref{dtDelta}, multiplying the obtained inequality by $2^{qs}$ and then summing with respect to $q \in \mathbb{Z},$ we obtain
\begin{multline} \label{4.14}
	\Vert e^{\mathcal{R}t} (\pa_tu)_\phi \Vert_{\tilde{L}^2_t(\mathcal{B}^s)} + \Vert e^{\mathcal{R}t} \pa_z u_\phi \Vert_{\tilde{L}^\infty_t(\mathcal{B}^s)}\leq \Vert e^{a|D_x|} \pa_z u_0 \Vert_{\mathcal{B}^s} + 2\Rcal \norm{e^{\Rcal t} \dd_z u_\phi}_{\tilde{L}^2_t(\Bcal^s)}\\ 
	+ C \Big( \Vert e^{\mathcal{R}t} (u\pa_xu)_\phi \Vert_{\tilde{L}^2_t(\mathcal{B}^s)}+\Vert e^{\mathcal{R}t} (v\pa_zu)_\phi \Vert_{\tilde{L}^2_t(\mathcal{B}^s)} + \Vert e^{\mathcal{R}t} \pa_xT_\phi \Vert_{\tilde{L}^2_t(\mathcal{B}^s)} \Big).
\end{multline}
We will accept for now the following estimates. The proof of these estimates will be given later.
\begin{lemma}
	\label{lem:31}
	Under the hypotheses of Theorem \ref{th:hydrolimit}, we have the following inequalities
	\begin{align*}
		\Vert e^{\mathcal{R}t} (u\pa_xu)_\phi \Vert_{\tilde{L}^2_t(\mathcal{B}^s)} &\leq C \Vert  u_\phi \Vert_{\tilde{L}^\infty(\mathcal{B}^\f12)} \Vert e^{\mathcal{R}t} \pa_z u_\phi \Vert_{\tilde{L}^2_t(\mathcal{B}^{s+1})}; \\
		\Vert e^{\mathcal{R}t} (v\pa_zu)_\phi \Vert_{\tilde{L}^2_t(\mathcal{B}^s)} &\leq C \Vert  u_\phi \Vert_{\tilde{L}^\infty(\mathcal{B}^\f12)} \Vert e^{\mathcal{R}t} \pa_z u_\phi \Vert_{\tilde{L}^2_t(\mathcal{B}^{s+1})} + \Vert  u_\phi \Vert_{\tilde{L}^\infty(\mathcal{B}^{s+1})} \Vert e^{\mathcal{R}t} \pa_z u_\phi \Vert_{\tilde{L}^2_t(\mathcal{B}^\f12)}.
	\end{align*}
\end{lemma}

Now, using Estimate \eqref{eq:hydrolimitU}, we have
\begin{align*}
	\norm{u_\phi}_{\tilde{L}^\infty_t(\Bcal^{\frac{1}2})} &\leq 2C \norm{e^{a\abs{D_x}} (u_0,T_0)}_{\Bcal^{\frac{1}2}}
	\\
	\norm{e^{\Rcal t} \dd_z u_\phi}_{\tilde{L}^2_t(\Bcal^{\frac{1}2})} &\leq 2C \norm{e^{a\abs{D_x}} (u_0,T_0)}_{\Bcal^{\frac{1}2}}
	\\
	\norm{e^{\Rcal t} \dd_z u_\phi}_{\tilde{L}^2_t(\Bcal^s)} &\leq 2C \norm{e^{a\abs{D_x}} (u_0,T_0)}_{\Bcal^s}
	\\
	\norm{e^{\Rcal t} \dd_x T_\phi}_{\tilde{L}^2_t(\Bcal^s)} &\leq 2C \norm{e^{a\abs{D_x}} (u_0,T_0)}_{\Bcal^s}
	\\
	\norm{u_\phi}_{\tilde{L}^\infty_t(\Bcal^{s+1})} &\leq 2C \norm{e^{a\abs{D_x}} (u_0,T_0)}_{\Bcal^{s+1}}
	\\
	\norm{e^{\Rcal t} \dd_z u_\phi}_{\tilde{L}^2_t(\Bcal^{s+1})} &\leq 2C \norm{e^{a\abs{D_x}} (u_0,T_0)}_{\Bcal^{s+1}}.
\end{align*}
Inserting all the above estimates into \eqref{4.14}, we finally obtain the existence of a constant $\overline{C} > 0$ such that
\begin{multline*}
	\Vert e^{\mathcal{R}t} (\pa_tu)_\phi \Vert_{\tilde{L}^2_t(\mathcal{B}^s)} + \Vert e^{\mathcal{R}t}  \pa_z u_\phi \Vert_{\tilde{L}^\infty_t(\mathcal{B}^s)}
	\\
	\leq \overline{C} \pare{ \norm{e^{a|D_x|} (u_0,T_0)}_{\mathcal{B}^s} + \norm{e^{a|D_x|} (u_0,T_0)}_{\mathcal{B}^\f12} \norm{e^{a|D_x|} (u_0,T_0)}_{\mathcal{B}^{s+1}} }.
\end{multline*}
Theorem \ref{th:hydrolimit} is proved. \hfill $\square$

\noindent \textit{Proof of Lemma \ref{lem:31}.}

The proof of this lemma is very similar to the proof of Lemmas \ref{lem:uww} and \ref{lem:vww}. We will give the main calculations without going into details. For the first inequality, we decompose
\begin{align*}
	(u\dd_x u)_\phi = (\Tcal^h_u\dd_x u)_\phi + (\Tcal^h_{\dd_x u}u)_\phi + (\Rcal^h(u,\dd_x u))_\phi.
\end{align*}
We have
\begin{align*}
	\norm{e^{\Rcal t} \DD^h_q (\Tcal^h_u\dd_x u)_\phi}_{L^2_t(L^2)} 
	&\leq \sum_{\abs{q-q'} \leq 4} \pare{\int_0^t \norm{S^h_{q'-1} u_\phi}_{L^\infty_xL^2_z}^2 \norm{e^{\Rcal t'} \DD^h_{q'} \dd_x u_\phi}_{L^2_xL^\infty_z}^2 dt'}^{\frac{1}2}
	\\
	&\lesssim \norm{u_\phi}_{L^\infty_t(\Bcal^{\frac{1}2})} \sum_{\abs{q-q'} \leq 4} 2^{q'} \norm{e^{\Rcal t} \DD^h_{q'} \dd_z u_\phi}_{L^2_t(L^2)} \\
	&\lesssim d_q 2^{-qs} \norm{u_\phi}_{L^\infty_t(\Bcal^{\frac{1}2})} \norm{e^{\Rcal t}\dd_z u_\phi}_{\tilde{L}^2_t(\mathcal{B}^{s+1})}
\end{align*}
where $$d_q = \sum_{\abs{q-q'} \leq 4} d_{q'}(\dd_z u_\phi) 2^{q'-q},$$
\begin{align*}
	\norm{e^{\Rcal t} \DD^h_q (\Tcal^h_{\dd_x u}u)_\phi}_{L^2_t(L^2)} 
	&\leq \sum_{\abs{q-q'} \leq 4} \pare{\int_0^t \norm{S^h_{q'-1} \dd_x u_\phi}_{L^\infty_xL^2_z}^2 \norm{e^{\Rcal t'} \DD^h_{q'} u_\phi}_{L^2_xL^\infty_z}^2 dt'}^{\frac{1}2}
	\\
	&\lesssim \sum_{\abs{q-q'} \leq 4} \pare{\int_0^t 2^{2q'} \norm{u_\phi}_{\Bcal^{\frac{1}2}}^2 \norm{e^{\Rcal t'} \DD^h_{q'} \dd_z u_\phi}_{L^2}^2 dt'}^{\frac{1}2}
	\\
	&\lesssim \norm{u_\phi}_{L^\infty_t(\Bcal^{\frac{1}2})} \sum_{\abs{q-q'} \leq 4} 2^{q'} \norm{e^{\Rcal t} \DD^h_{q'} \dd_z u_\phi}_{L^2_t(L^2)} \\
	&\lesssim d_q 2^{-qs} \norm{u_\phi}_{L^\infty_t(\Bcal^{\frac{1}2})} \norm{e^{\Rcal t}\dd_z u_\phi}_{\tilde{L}^2_t(\mathcal{B}^{s+1})}
\end{align*}
where $$d_q = \sum_{\abs{q-q'} \leq 4} d_{q'}(\dd_z u_\phi) 2^{q'-q}.$$
Finally, we have
\begin{align*}
	\norm{e^{\Rcal t} \DD^h_q (\Rcal^h(u,\dd_x u))_\phi}_{L^2_t(L^2)} &\leq \sum_{q' \geq q-3} \pare{\int_0^t \norm{\Delta^h_{q'} \dd_x u_\phi}_{L^\infty_xL^2_z}^2 \norm{e^{\Rcal t'} \DD^h_{q'} u_\phi}_{L^2_xL^\infty_z}^2 dt'}^{\frac{1}2}
	\\
	&\lesssim \sum_{q' \geq q-3} \pare{\int_0^t 2^{2q'} \norm{u_\phi}_{\Bcal^{\frac{1}2}}^2 \norm{e^{\Rcal t'} \DD^h_{q'} \dd_z u_\phi}_{L^2}^2 dt'}^{\frac{1}2}
	\\
	&\lesssim \norm{u_\phi}_{L^\infty_t(\Bcal^{\frac{1}2})} \sum_{q' \geq q-3} 2^{q'} \norm{e^{\Rcal t} \DD^h_{q'} \dd_z u_\phi}_{L^2_t(L^2)} \\
	&\lesssim d_q 2^{-qs} \norm{u_\phi}_{L^\infty_t(\Bcal^{\frac{1}2})} \norm{e^{\Rcal t}\dd_z u_\phi}_{\tilde{L}^2_t(\mathcal{B}^{s+1})},
\end{align*}
where $$d_q = \sum_{q' \geq q-3} d_{q'}(\dd_z u_\phi) 2^{-(q-q')}.$$ Here, for the third estimate, we remark that the sequence $(d_q)_{q\in\ZZ}$ can be considered as the convolution of two summable sequences and is then also summable.

For the second inequality, we can perform the same calculations, while taking into account the ``transfer of regularity'' given in Lemma \ref{lem:tranreg}. \hfill $\square$

}

\section{Global well-posedness of the 2D non-rotating primitive equations in a thin strip} \label{se:PEhydro}

{\color{black}

In this section, we will prove Theorem \ref{th:primitive} and establish the global well-posedness of the system \eqref{eq:hydroPE} for small analytic data. We will use the same technique as in the previous section but will consider a different auxiliary function to control the evolution of the analyticity of the solutions. This auxiliary function is chosen to be adapted to the primitive system \eqref{eq:hydroPE}. Let $\Theta: \RR_+\times \RR \to \RR_+$ such that 
\begin{equation*}
	\Theta(0, \xi) = 0, \quad \Theta(t,\xi) = (a- \lambda \tau(t))|\xi|, \ \forall \, t > 0, \, \forall \, \xi \in \RR,
\end{equation*}
where $\lambda > 0$ and $\tau(t)$ will be determined later. For any function $f \in L^2(\Scal)$, we define
\begin{align*}
	\Theta: f \mapsto f_\Theta; \quad f_{\Theta}(t,x,z) = e^{\Theta(t,D_x)} f(t,x,z) = \mathcal{F}^{-1}_h(e^{\Theta(t,\xi)}\widehat{f}(t,\xi,z)).
\end{align*}

In what follows, for the sake of the simplicity, we will neglect the index $\eps$ and write $(u_{\Theta},v_{\Theta},T_{\Theta})$ instead of $(u_{\Theta}^\eps,v_{\Theta}^\eps,T_{\Theta}^\eps)$. Applying the operator $\Theta$ to the system \eqref{eq:hydroPE}, we obtain
\begin{equation}\label{eq4}
	\quad\left\{\begin{array}{l}
		\displaystyle \partial_t u_{\Theta} + \lambda\dot{\tau}(t)|D_x|u_{\Theta} +(u\partial_x u)_{\Theta} + (v\pa_zu)_{\Theta}-\eps^2\partial_x^2u_{\Theta}-\pa_z^2u_{\Theta}+\partial_x p_{\Theta}=0,\ \ \\
		\displaystyle \eps^2\left(\partial_tv_{\Theta}+(u\partial_x v)_{\Theta}+(v\pa_zv)_{\Theta}-\eps^2\partial_x^2v_{\Theta}-\pa_z^2v_{\Theta} \right)+\pa_zp_{\Theta}=T_{\Theta},\\
		\partial_t T_{\Theta}+(u\partial_x T)_{\Theta}+ (v\pa_z T)_{\Theta}-\Delta T_{\Theta}=0,\\
		\displaystyle \partial_x u_{\Theta}+\pa_zv_{\Theta}=0,\\
		\displaystyle \left(u_{\Theta}, v_{\Theta}, T_{\Theta} \right)|_{z=0}=\left(u_{\Theta}, v_{\Theta},T_{\Theta} \right)|_{z=1} = 0,\\
		\displaystyle \left(u_{\Theta}, v_{\Theta}, T_{\Theta} \right)|_{t=0}=\left(e^{a\abs{D_x}} u_0, e^{a\abs{D_x}} v_0, e^{a\abs{D_x}} T_0 \right).
	\end{array}\right.
\end{equation}
We remark that the pressure term is not really an unknown and can be determined as function of $(u_\Theta,T_\Theta)$ as we did for the hydrostatic limit system (see also \cite{CLT2019} for more details). We recall also that we always use ``$C$'' to denote a generic positive constant which can change from line to line.

Applying the operator $\Delta^h_q$ to the system \eqref{eq4}, then taking the $L^2(\mathcal{S})$ scalar product of the first three equations of the obtained system with $\Delta^h_q u_{\phi}$, $\Delta^h_q v_{\phi} $ and $\Delta^h_q T_{\phi}$ respectively and summing up the first and second equations, we get 
\begin{align} \label{S5eq5}
	\frac{1}{2}\frac{d}{dt} &\Vert \Delta_q^h (u_{\Theta},\eps v_{\Theta})(t) \Vert_{L^2}^2 + \lambda \dot{\tau}(t)\psca{|D_x| \Delta_q^h (u_{\Theta}, \eps v_{\Theta}),\Delta_q^h (u_{\Theta},\eps v_{\Theta})}_{L^2}\\ 
	&+ \Vert \pa_z \Delta_q^h (u_{\Theta},\eps v_{\Theta}) \Vert_{L^2}^2 + \eps^2  \Vert \pa_x \Delta_q^h (u_{\Theta},\eps v_{\Theta}) \Vert_{L^2}^2 \notag\\
	&= -\psca{\Delta_q^h (u\p_x u)_{\Theta}, \Delta_q^h u_{\Theta}}_{L^2} - \psca{\Delta_q^h (v\pa_z u)_{\Theta}, \Delta_q^h u_{\Theta}}_{L^2} - \psca{\na \Delta_q^h p_{\Theta}, \Delta_q^h (u_{\Theta}, v_{\Theta})}_{L^2} \notag\\
	& -\eps^2 \psca{\Delta_q^h (u\p_x v)_{\Theta}, \Delta_q^h v_{\Theta}}_{L^2} - \eps^2\psca{\Delta_q^h (v\pa_z v)_{\Theta}, \Delta_q^h v_{\Theta}}_{L^2} +  \psca{\Delta_q^h T_{\Theta}, \Delta_q^h v_{\Theta}}_{L^2}, \notag
\end{align}
and 
\begin{multline} \label{S5eq6}
	\frac{1}{2}\frac{d}{dt} \Vert \Delta_q^h T_{\Theta}(t) \Vert_{L^2}^2 + \lambda \dot{\tau}(t)\psca{|D_x| \Delta_q^h T_{\Theta},\Delta_q^h T_{\Theta}}_{L^2}+ \Vert \nabla \Delta_q^h T_{\Theta} \Vert_{L^2}^2 \\
	= -\psca{\Delta_q^h (u\p_x T)_{\Theta}, \Delta_q^h T_{\Theta}}_{L^2} - \psca{\Delta_q^h (v\pa_z T)_{\Theta}, \Delta_q^h T_{\Theta}}_{L^2}.
\end{multline} 
As in the previous section, we will multiply \eqref{S5eq5} and \eqref{S5eq6} by $e^{2\mathcal{R}t}$, then we integrate the obtained equations with respect to the time variable and get, for any $0 < \Rcal < \frac{1}{2\Kcal}$,
\begin{multline} \label{S5eq5bis}
	\Vert e^{\mathcal{R}t} \Delta_q^h  (u_{\Theta},\eps v_{\Theta})(t) \Vert_{L^\infty_t(L^2)}^2 + 2\lambda \int_0^t \dot{\tau}(t')\norm{ e^{\mathcal{R}t'} |D_x|^\f12 \Delta_q^h (u_{\Theta}, \eps v_{\Theta})(t')}_{L^2}^2 dt' \\ + \Vert e^{\mathcal{R}t'}\pa_z \Delta_q^h (u_{\Theta},\eps v_{\Theta}) \Vert_{L^2_t(L^2)}^2 + \eps^2  \Vert e^{\mathcal{R}t'}\pa_x \Delta_q^h (u_{\Theta},\eps v_{\Theta}) \Vert_{L^2_t(L^2)}^2\\
	= \norm{\Delta_q^h (u_{\Theta},\eps v_{\Theta})(0)}_{L^2}^2 + F_1 +F_2 +F_3 +F_4 + F_5 + F_6,
\end{multline}
and
\begin{multline} \label{S5eq6bis}
	\norm{e^{\Rcal t} \Delta_q^h T_{\Theta}(t)}_{L^\infty_t (L^2)}^2 + 2\lambda \int_0^t \dot{\tau}(t') \norm{e^{\Rcal t} |D_x|^{\frac12} \Delta_q^h T_{\Theta}}_{L^2}^2 dt' + \norm{e^{\Rcal t} \Delta_q^h \nabla T_{\Theta}(t)}_{L^2_tL^2}^2\\
	= \norm{\Delta_q^h T_{\Theta}(0)}_{L^2}^2 + F_7 + F_8. \hspace{1.5cm}
\end{multline}

From now on, we will fix
\begin{equation}
	\label{eq:deftau}
	\dot{\tau}(t) =  \Vert \pa_z u_{\Theta}(t) \Vert_{\mathcal{B}^{\frac{1}{2}}} + \Vert \pa_z T_{\Theta}(t) \Vert_{\mathcal{B}^{\frac{1}{2}}} + \eps \Vert \pa_z v_{\Theta}(t) \Vert_{\mathcal{B}^{\frac{1}{2}}}.
\end{equation}
Using Lemma \ref{lem:uww}, we get the following controls
\begin{align*}
	\abs{F_1} = 2\abs{\int_0^t \psca{e^{\mathcal{R}t}\Delta_q^h (u\p_x u)_{\Theta}, e^{\mathcal{R}t}\Delta_q^h u_{\Theta}}_{L^2} dt'} 
	\leq C d_q^2 2^{-2qs} \norm{e^{\mathcal{R}t} u_{\Theta}}_{\tilde{L}^2_{t,\dot{\tau}(t)}(\mathcal{B}^{s+\frac{1}{2}})}^2,
\end{align*}
\begin{align*}
	\abs{F_4} = 2\abs{\int_0^t \psca{e^{\mathcal{R}t} \Delta_q^h (u\p_x (\eps v))_{\Theta},e^{\mathcal{R}t} \Delta_q^h (\eps v)_{\Theta}}_{L^2} dt'} \leq Cd_q^2 2^{-2qs} \norm{e^{\mathcal{R}t} (u_{\Theta},\eps v_{\Theta})}_{\tilde{L}^2_{t,\dot{\tau}(t)}(\mathcal{B}^{s+\frac{1}{2}})}^2 
\end{align*}
and
\begin{align*}
	\abs{F_7} = 2\abs{\int_0^t \psca{e^{\Rcal t'}\Delta_q^h (u\p_x T)_{\Theta}, e^{\Rcal t'}\Delta_q^h T_{\Theta}} dt'} \leq C d_q^2 2^{-2qs} \norm{e^{\Rcal t} (u_\Theta,T_{\Theta})}_{\tilde{L}^2_{t,\dot{\tau}(t)}(\mathcal{B}^{s+\frac{1}{2}})}^2.
\end{align*}
Next, Lemma \ref{lem:vww} implies the following controls
\begin{multline*}
	\abs{F_2} = 2\abs{\int_0^t \psca{e^{\mathcal{R}t}\Delta_q^h (v\pa_z u)_{\Theta}, e^{\mathcal{R}t}\Delta_q^h u_{\Theta}}_{L^2}  dt' } 
	\\
	\leq C d_q^2 2^{-2qs} \pare{\norm{e^{\mathcal{R}t} u_{\Theta}}_{\tilde{L}^2_{t,\dot{\tau}(t)}(\mathcal{B}^{s+\frac{1}{2}})}^2 + \norm{u_{\Theta}}_{L^\infty_t(\mathcal{B}^{\frac{3}{2}})} \norm{e^{\mathcal{R}t} \dd_z u_{\Theta}}_{\tilde{L}^2_t(\mathcal{B}^s)}^2},
\end{multline*}
\begin{multline*}
	\abs{F_5} = 2\abs{\int_0^t \psca{e^{\mathcal{R}t}\Delta_q^h (v\pa_z (\eps v))_{\Theta}, e^{\mathcal{R}t} \Delta_q^h (\eps v)_\Theta}_{L^2} dt'} 
	\\ 
	\leq C d_q^2 2^{-2qs} \pare{\norm{e^{\mathcal{R}t} (u_{\Theta},\eps v_{\Theta})}_{\tilde{L}^2_{t,\dot{\tau}(t)}(\mathcal{B}^{s+\frac{1}{2}})}^2 + \norm{u_{\Theta}}_{L^\infty_t(\mathcal{B}^{\frac{3}{2}})} \norm{e^{\mathcal{R}t} \dd_z (\eps v)_{\Theta}}_{\tilde{L}^2_t(\mathcal{B}^s)}^2}
\end{multline*}
and
\begin{multline*}
	\abs{F_8} = 2\abs{\int_0^t \psca{e^{\Rcal t'}\Delta_q^h (v\pa_z T)_{\Theta}, e^{\Rcal t'}\Delta_q^h T_{\Theta}} dt'}
	\\
	\leq C d_q^2 2^{-2qs} \pare{\norm{e^{\Rcal t} (u_{\phi},T_{\phi})}_{\tilde{L}^2_{t,\dot{\tau}(t)}(\mathcal{B}^{s+\frac{1}{2}})}^2 + \norm{u_{\Theta}}_{L^\infty_t(\mathcal{B}^{\frac{3}{2}})} \norm{e^{\mathcal{R}t} \dd_z T_{\Theta}}_{\tilde{L}^2_t(\mathcal{B}^s)}^2}.
\end{multline*}
For the term $F_3$, using the divergence-free property $\partial_x u_{\Theta}+\pa_zv_{\Theta}=0$ and an integration by parts, we deduce that
\begin{align*}
	\abs{F_3} = 2\abs{\int_0^t \psca{e^{\Rcal t'} \na \Delta_q^h p_{\Theta}, e^{\Rcal t'} \Delta_q^h (u_{\Theta}, v_{\Theta})}_{L^2} dt'} = 0.
\end{align*}
For the last term $F_4$, the boundary condition $\left(u_{\Theta}, v_{\Theta} \right)|_{z\in\set{0,1}} = 0$ and the relation $$v_{\Theta} (t,x,y) = -\int_0^z \partial_x u_{\Theta}(t,x,s) ds$$ imply 
\begin{align*}
	 \abs{ \psca{e^{\Rcal t} \Delta_q^h T_{\Theta}, e^{\Rcal t} \Delta_q^h v_{\Theta}}_{L^2}}  &= \abs{\psca{e^{\Rcal t} \Delta_q^h  T_{\Theta}, e^{\Rcal t} \Delta_q^h  \int_0^z -\pa_x u_{\Theta}ds}_{L^2}} \\
	& \leq \Vert \Delta_q^h e^{\Rcal t} \pa_x T_{\Theta} \Vert_{L^2} \Vert \Delta_q^h e^{\Rcal t}  u_{\Theta} \Vert_{L^2} \\ 
	&\leq C \Vert \Delta_q^h e^{\Rcal t} \pa_x T_{\Theta} \Vert_{L^2}^2 + \frac{1}{16}  \Vert \Delta_q^h e^{\Rcal t} u_{\Theta} \Vert_{L^2}^2.
\end{align*}
Then, Poincar\'e and Cauchy-Schwartz inequalities yield 
\begin{equation*}
	\abs{F_6} = 2\abs{ \int_0^t  \psca{e^{\Rcal t'} \Delta_q^h T_{\Theta}, e^{\Rcal t'} \Delta_q^h v_{\Theta}}_{L^2}dt' } \leq  C d_q^2 2^{-2qs} \Vert e^{\mathcal{R}t} \dd_x T_{\Theta} \Vert_{\tilde{L}^2_t(\mathcal{B}^s)}^2 + \frac{1}8 \Vert e^{\mathcal{R}t} \Delta_q^h \pa_z u_{\Theta} \Vert_{L^2_t(L^2)}^2.
\end{equation*}

We now multiply \eqref{S5eq5bis} and \eqref{S5eq6bis} by $2^{2qs}$ and we recall that we can take square root of each terms with the cost of a multiplier which will be included in the generic constant $C$. Summing the resulting inequalities with respect to $q\in \mathbb{Z}$, we have 
\begin{multline*}
	\Vert e^{\mathcal{R}t} (u_{\Theta},\eps v_{\Theta},T_{\Theta}) \Vert_{\tilde{L}^{\infty}_t(\mathcal{B}^s)} + \sqrt{2\lambda} \Vert e^{\mathcal{R}t} (u_{\Theta}, \eps v_{\Theta},T_{\Theta})\Vert_{\tilde{L}^2_{t,\dot{\tau}(t)}(\mathcal{B}^{s+\frac{1}{2}})} + \norm{e^{\Rcal t} \dd_x T_{\Theta}}_{\tilde{L}^2_t(\Bcal^s)}
	\\ 
	+ \Vert e^{\mathcal{R}t} \pa_z (u_{\Theta},\eps v_{\Theta},T_\Theta)\Vert_{\tilde{L}^{2}_t(\mathcal{B}^s)} + \eps \Vert e^{\mathcal{R}t}  \pa_x (u_{\Theta},\eps v_{\Theta}) \Vert_{\tilde{L}^{2}_t(\mathcal{B}^s)}
	\leq \norm{e^{a\abs{D_x}} (u_0,\eps v_0, T_0)}_{\mathcal{B}^s}
	\\ 
	+ C \Vert e^{\mathcal{R}t} (u_{\Theta},\eps v_{\Theta},T_\Theta)\Vert_{\tilde{L}^2_{t,\dot{\tau}(t)}(\mathcal{B}^{s+\frac{1}{2}})}  + \f18 \norm{e^{\Rcal t} \dd_y u_\Theta}_{\tilde{L}^2_t(\Bcal^s)}
	\\
	+ C \norm{u_{\Theta}}_{L^\infty_t(\mathcal{B}^{\frac{3}{2}})}^{\frac{1}2} \norm{e^{\mathcal{R}t} \dd_z (u_{\Theta}, \eps v_\Theta, T_\Theta)}_{\tilde{L}^2_t(\mathcal{B}^s)}.
\end{multline*}
Without loss of generality, we can suppose that $C \geq 2$ and we choose $\sqrt{\lambda} > 3C$. For $\eps \leq \frac{1}{2C}$, we can simply the above inequality as follows
\begin{multline}
	\label{eq:u+v+T02}
	\norm{e^{\mathcal{R}t} (u_{\Theta},\eps v_{\Theta},T_{\Theta})}_{\tilde{L}^{\infty}_t(\mathcal{B}^s)} + \sqrt{\lambda} \norm{e^{\mathcal{R}t} (u_{\Theta}, \eps v_{\Theta},T_{\Theta})}_{\tilde{L}^2_{t,\dot{\tau}(t)}(\mathcal{B}^{s+\frac{1}{2}})} + \norm{e^{\Rcal t} \dd_x T_{\Theta}}_{\tilde{L}^2_t(\Bcal^s)}
	\\ 
	+ \frac{1}{2} \norm{e^{\mathcal{R}t} \pa_z (u_{\Theta},\eps v_{\Theta},T_\Theta)}_{\tilde{L}^{2}_t(\mathcal{B}^s)} \leq \norm{e^{a\abs{D_x}} (u_0,\eps v_0, T_0)}_{\mathcal{B}^s} 
	\\
	+ C \norm{u_{\Theta}}_{L^\infty_t(\mathcal{B}^{\frac{3}{2}})}^{\frac{1}2} \norm{e^{\mathcal{R}t} \dd_z (u_{\Theta}, \eps v_\Theta, T_\Theta)}_{\tilde{L}^2_t(\mathcal{B}^s)}.
\end{multline}

Let 
\begin{equation*} 
	t^\star \stackrel{\tiny def}{=} \sup\set{t>0\ : \ \norm{u_{\Theta}}_{\mathcal{B}^{\frac{3}{2}}} \leq \frac{1}{4C} \ \mbox{ and } \ \tau(t) \leq  \frac{a}{3\lambda}},
\end{equation*}
For initial data such that
\begin{displaymath}
	\left\{
	\begin{aligned}
		\norm{e^{a\abs{D_x}} (u_0,\eps v_0, T_0)}_{\mathcal{B}^{\frac{1}2}} &\leq \frac{a\sqrt{2\Rcal}}{8\lambda}
		\\
		\norm{e^{a\abs{D_x}} (u_0,\eps v_0, T_0)}_{\mathcal{B}^{\frac{3}2}} &\leq \frac{1}{8C}
	\end{aligned}
	\right.
\end{displaymath} 
the fact that $\tau(0) = 0$ and the continuity in time in $\Bcal^s$ imply that $t^\star > 0$. 

\noindent For $s=\frac{3}2$ and for any $0 < t < t^*$, from \eqref{eq:u+v+T02} we have
\begin{displaymath}
	\norm{u_{\Theta}}_{\mathcal{B}^{\frac{3}{2}}} \leq \norm{e^{\mathcal{R}t} (u_{\Theta},\eps v_{\Theta},T_{\Theta})}_{\tilde{L}^{\infty}_t(\mathcal{B}^{\frac{3}{2}})} \leq \norm{e^{a\abs{D_x}} (u_0,\eps v_0, T_0)}_{\mathcal{B}^{\frac{3}{2}}} \leq \frac{1}{8C}.
\end{displaymath}
For $s=\frac{1}2$ and for any $0 < t < t^*$, Inequality \eqref{eq:u+v+T02} also yields
\begin{align*}
	\tau (t) &= \int_0^t \norm{\pa_z (u_{\Theta},\eps v_{\Theta},T_{\Theta})(t')}_{\mathcal{B}^{\frac{1}{2}}} dt' 
	\\
	& \leq \int_0^t  e^{-\mathcal{R}t'} \norm{e^{\mathcal{R}t'} \pa_z (u_{\Theta},\eps v_{\Theta},T_{\Theta})(t')}_{\mathcal{B}^{\frac{1}{2}}} dt' 
	\\
	&\leq \left( \int_0^t  e^{-2\mathcal{R}t'} dt' \right)^{\f12} \left( \int_0^t \norm{e^{\mathcal{R}t'} \pa_z (u_{\Theta},\eps v_{\Theta},T_{\Theta})(t')}_{\mathcal{B}^{\frac{1}{2}}}^2 dt' \right)^{\f12} 
	\\
	&\leq \frac{1}{\sqrt{2\Rcal}} \Vert e^{\mathcal{R}t} \dd_z(u_{\Theta},\eps v_{\Theta},T_{\Theta}) \Vert_{\tilde{L}^2_t(\mathcal{B}^{\f12})} \\
	&\leq \frac{2}{\sqrt{2\Rcal}} \norm{e^{a\abs{D_x}} (u_0,\eps v_0, T_0)}_{\mathcal{B}^{\f12}} < \frac{a}{4\lambda}.
\end{align*}
A continuity argument implies that $t^\star = +\infty$ and allows to close the proof of Theorem \ref{th:primitive}.

}

\section{Convergence to the hydrostatic limit system} \label{se:convergence}

{\color{black}

In this section, we will prove Theorem \ref{th:Convergence} and justify the approximation of the scaled non-rotating primitive equations and the hydrostatic limit system in a two-dimensional thin strip. To this end, we introduce the following difference quantities
\begin{equation*}
	w^{\eps,1} = u^{\eps} - u, \quad w^{\eps,2} = v^{\eps} - v, \quad \theta^{\eps} = T^{\eps}- T, \quad q^{\eps} = p^{\eps}-p
\end{equation*}
where $(u^\eps,v^\eps,T^\eps,p^\eps)$ and $(u,v,T,p)$ are respectively solutions of the systems \eqref{eq:hydroPE} and \eqref{eq:hydrolimit}. We deduce that $(w^{\eps,1},w^{\eps,2},\theta^\eps,q^\eps)$ satisfies the following system
\begin{equation}\label{S1eq10}
	\left\{
	\begin{aligned}
		&\partial_t w^{\eps,1} - \eps^2\partial_x^2 w^{\eps,1}-\pa_z^2 w^{\eps,1} + \partial_x q^\eps= R^{1,\eps}
		\\
		&\eps^2\left(\partial_t w^{\eps,2} - \eps^2\partial_x^2 w^{\eps,2}-\pa_z^2 w^{\eps,2}\right)+\pa_z q^\eps =\theta^\eps + R^{2,\eps}
		\\
		&\partial_t \theta^\eps - \pa_x^2 \theta^\eps - \pa_z^2 \theta^\eps = R^{3,\eps}
		\\
		&\partial_x w^{\eps,1}+\pa_z  w^{\eps,2}=0
		\\
		&\left( w^{\eps,1}, w^{\eps,2}, \theta^\eps \right)|_{t=0}=\left(u_0^\eps-u_0, v_0^{\eps} - v_0,T_0^{\eps}- T_0 \right)
		\\
		&\left( w^{\eps,1}, w^{\eps,2}, \theta^\eps \right)|_{z=0}=\left( w^{\eps,1}, w^{\eps,2}, \theta^\eps \right)|_{z=1} = 0,
	\end{aligned}
	\right.
\end{equation}
where the remaining terms $R^{i,\eps}$, with $ i =1,2,3$, are given by 
\begin{align} \label{eq:Ri}
	\left\{
	\begin{aligned}
		R^{1,\eps} &= \eps^2 \pa_x^2 u -(u^{\eps} \pa_x u^{\eps} - u\pa_xu )- (v^\eps \pa_z u^\eps - v\pa_z u), \\
		R^{2,\eps} &= - \eps^2\left(\pa_t v -\eps^2 \pa_x^2 v -\pa_z^2 v + u^{\eps} \pa_x v^{\eps} + v^\eps \pa_z v^\eps \right), \\
		R^{3,\eps} &= -\left(u^\eps\pa_xT^{\eps} + v^\eps\pa_zT^{\eps} \right) + \left( u\pa_xT + v\pa_zT \right).
	\end{aligned}
	\right.
\end{align}

Let $\varphi: \RR_+\times \RR \to \RR_+$ such that 
\begin{equation*}
	\varphi(0, \xi) = 0, \quad \varphi(t,\xi) = (a- \mu \eta(t))|\xi|, \ \forall \, t > 0, \, \forall \, \xi \in \RR,
\end{equation*}
where $\mu \geq \lambda > 0$ and $\eta(t)$ will be determined later. For any function $f \in L^2(\Scal)$, we define
\begin{align*}
	\varphi: f \mapsto f_\varphi; \quad f_{\varphi}(t,x,z) = e^{\varphi(t,D_x)} f(t,x,z) = \mathcal{F}^{-1}_h(e^{\varphi(t,\xi)}\widehat{f}(t,\xi,z)).
\end{align*}
In what follows, for the sake of the simplicity, we will drop the index $\eps$ and write $(w^{1}_{\varphi},w^2_{\varphi},\theta_\varphi,q_\varphi,R^{i}_\varphi)$ instead of $(w^{\eps,1}_{\varphi},w^{\eps,2}_{\varphi},\theta_{\varphi}^\eps,q^\eps,R^{i,\eps}_\varphi)$. Direct calculations show that $(w^{1}_{\varphi},w^2_{\varphi},\theta_\varphi,q_\varphi)$ satisfies
\begin{equation}
	\label{22}
	\left\{
	\begin{aligned}
		&\partial_t w^{1}_{\varphi} + \mu |D_x| \dot{\eta}(t) w^{1}_{\varphi}  - \eps^2\partial_x^2 w^{1}_{\varphi}-\pa_z^2 w^{1}_{\varphi} + \partial_x q_\varphi = R^{1}_\varphi
		\\
		&\eps^2\left(\partial_t w^2_{\varphi} +\mu |D_x| \dot{\eta}(t) w^{2}_{\varphi} - \eps^2\partial_x^2 w^2_{\varphi}-\pa_z^2 w^2_{\varphi}\right)+\pa_z q_\varphi =\theta_\varphi + R^{2}_\varphi
		\\
		&\partial_t \theta_\varphi + \mu |D_x| \dot{\eta}(t) \theta_{\varphi} - \pa_x^2 \theta_\varphi - \pa_z^2 \theta_\varphi = R^3_\varphi 
		\\
		&\partial_x w^{1}_{\varphi}+\pa_z  w^2_{\varphi}=0
		\\
		&\left( w^{1}_{\varphi},w^2_{\varphi},\theta_\varphi \right)|_{t=0}=\left(u_0^\eps-u_0, v_0^{\eps} - v_0,T_0^{\eps}- T_0 \right)
		\\
		&\left( w^{1}_{\varphi},w^2_{\varphi},\theta_\varphi \right)|_{z=0}=\left( w^{1}_{\varphi},w^2_{\varphi},\theta_\varphi \right)|_{z=1} = 0.
	\end{aligned}
	\right.
\end{equation}
As in the previous sections, we will use ``$C$'' to denote a generic positive constant which can change from line to line.

Applying the dyadic operator $\Delta^h_q$ to the system \eqref{22}, then taking the $L^2(\mathcal{S})$ scalar product of the first, second and the third equations of the obtained system with $\Delta_q^h w^1_{\varphi}$, $\Delta_q^h w^2_{\varphi} $ and $\Delta_q^h \theta_{\varphi}$ respectively, we obtain 
\begin{multline} \label{S6eq5}
	\frac{1}{2}\frac{d}{dt} \Vert \Delta_q^h (w^1_\varphi,\eps w^2_\varphi)(t) \Vert_{L^2}^2 + \mu \dot{\eta}(t)(|D_x| \Delta_q^h (w^1_\varphi,\eps w^2_\varphi),\Delta_q^h (w^1_\varphi,\eps w^2_\varphi))_{L^2}\\ 
	+ \Vert \pa_z \Delta_q^h (w^1_\varphi,\eps w^2_\varphi) \Vert_{L^2}^2 + \eps^2  \Vert \pa_x \Delta_q^h (w^1_\varphi,\eps w^2_\varphi) \Vert_{L^2}^2 = \psca{\Delta_q^h R^1_{\varphi}, \Delta_q^h w^1_{\varphi}}_{L^2}\\ 
	+ \psca{\Delta_q^h R^2_{\varphi}, \Delta_q^h w^2_{\varphi}}_{L^2} - \psca{\Delta_q^h \na q_{\varphi}, \Delta_q^h (w^1_\varphi, w^2_\varphi)}_{L^2} + \psca{\Delta_q^h \theta_{\phi},\Delta_q^h w^2_\phi}_{L^2}
\end{multline}
and 
\begin{equation} \label{S6eq6}
	\frac{1}{2}\frac{d}{dt} \Vert \Delta_q^h \theta_{\varphi}(t) \Vert_{L^2}^2 +  \mu \dot{\eta}(t)\psca{|D_x| \Delta_q^h \theta_{\Theta},\Delta_q^h \theta_{\varphi}}_{L^2}+ \Vert \nabla \Delta_q^h \theta_{\varphi} \Vert_{L^2}^2 = \psca{\Delta_q^h R^3_{\varphi}, \Delta_q^h \theta_{\varphi}}_{L^2}.
\end{equation}
Integrating $\eqref{S6eq5}$ and $\eqref{S6eq6}$ with respect to the time variable, we have 
\begin{multline} \label{S6eq5bis}
	\Vert \Delta_q^h (w^1_\varphi,\eps w^2_\varphi) \Vert_{L^\infty_t(L^2)}^2 + \mu \int_0^t \dot{\eta}(t')\psca{|D_x| \Delta_q^h (w^1_\varphi,\eps w^2_\varphi),\Delta_q^h (w^1_\varphi,\eps w^2_\varphi)}_{L^2} dt' 
	\\ 
	+ \Vert \pa_z \Delta_q^h (w^1_{\varphi},\eps w^2_{\varphi}) \Vert_{L^2_t(L^2)}^2 + \eps^2 2^q \Vert \Delta_q^h (w^1_{\varphi},\eps w^2_{\varphi}) \Vert_{L^2_t(L^2)}^2\\
	\leq \norm{\Delta_q^h (w^1_{\varphi},\eps w^2_{\varphi})(0)}_{L^2}^2 + G_1^q +G_2^q +G_3^q +G_4^q,
\end{multline}
and 
\begin{multline} \label{S6eq6bis}
	\Vert \Delta_q^h \theta_{\varphi} \Vert_{L^\infty_t(L^2)}^2 + \mu \int_0^t  \dot{\eta}(t')\psca{|D_x| \Delta_q^h \theta_{\Theta},\Delta_q^h \theta_{\varphi}}_{L^2} dt' + \Vert \dd_z \Delta_q^h \theta_{\varphi} \Vert_{L^2_t(L^2)}^2 + 2^q \Vert \Delta_q^h \theta_{\varphi} \Vert_{L^2_t(L^2)}^2 
	\\
	\leq \norm{\Delta_q^h \theta_{\varphi}(0)}_{L^2}^2 + G_5^q,
\end{multline}
where the terms $G_i$, $i= 1,\ldots,5,$ will be precised and controlled in what follows.

We will now estimate the linear terms on the right-hand side of the above inequalities. First of all, for pressure term, using the incompressibility property $\partial_x w^1_{\varphi}+\pa_zw^2_{\varphi}=0$, we perform an integration by parts and get
\begin{align}
	\label{eq:G3}
	G_3^q = 2 \abs{\int_0^t \psca{\na \Delta_q^h q_{\varphi}, \Delta_q^h (w^1_{\varphi}, w^2_{\varphi})}_{L^2} dt'} =  2 \abs{\int_0^t \psca{ \Delta_q^h q_{\varphi}, \div (\Delta_q^h (w^1_{\varphi}, w^2_{\varphi}))}_{L^2} dt'} = 0.
\end{align}
For the temperature term, the boundary condition $ \left(u_{\Theta}, v_{\Theta} \right)|_{z=0}=\left(u_{\Theta}, v_{\Theta} \right)|_{z=1} = 0$ and the incompressibility property $\partial_x w^1_{\varphi}+\pa_zw^2_{\varphi}=0$ allow to write $w^2_\varphi$ as function of $w^1_\varphi$
\begin{equation*}
	w^2_\varphi = -\int_0^z \pa_x w^1_{\varphi}(x,z')dz'.
\end{equation*}
Then, using integration by parts, Poincar\'e inequality and Bernstein lemma \ref{le:Bernstein}, we have 
\begin{align}
	\label{eq:G4}
	G_4^q &= 2\abs{\int_0^t \psca{\Delta_q^h \theta_{\varphi},\Delta_q^h w^2_{\varphi}}_{L^2} dt'} = 2\abs{\int_0^t \psca{\Delta_q^h \theta_{\varphi}, -\Delta_q^h\int_0^z \pa_x w^1_{\varphi}(x,s)ds}_{L^2} dt'}
	\\
	& \leq 2\int_0^t \abs{\psca{\Delta_q^h \pa_x \theta_{\varphi},\Delta_q^h \int_0^z w^1_{\varphi}(x,s)ds}_{L^2}} dt' \leq 2\int_0^t \Vert \Delta_q^h \pa_x \theta_{\varphi} \Vert_{L^2} \Vert \Delta_q^h w^1_{\varphi} \Vert_{L^2} dt' \notag
	\\
	&\leq C 2^q \Vert \Delta_q^h \theta_{\varphi} \Vert_{L^2_t(L^2)}^2 + \frac{1}{16} \Vert \dd_z \Delta_q^h w^1_{\varphi} \Vert_{L^2_t(L^2)}^2 \notag
\end{align}
The goal of the main part of this section is to estimate the nonlinear terms on the right-hand side of \eqref{S6eq5bis} and \eqref{S6eq6bis}, which are
\begin{align*}
	G_1^q &= 2\abs{\int_0^t \psca{ \Delta_q^h R^1_{\varphi}, \Delta_q^h w^1_{\varphi}}_{L^2} dt'} \\
	G_2^q &= 2\abs{\int_0^t \psca{\Delta_q^h R^2_{\varphi}, \Delta_q^h w^2_{\varphi}}_{L^2}dt'} \\
	G_5^q &= 2\abs{\int_0^t \psca{\Delta_q^h R^3_{\varphi}, \Delta_q^h \theta_{\varphi}}_{L^2} dt'},
\end{align*}
with $R^i_\varphi$ being defined in \eqref{eq:Ri}. From now on, we will set
\begin{align*}
	\eta(t) = \int_0^t \left( \Vert (\pa_z u^{\eps}_{\Theta}, \eps \pa_x u^{\eps}_{\Theta}, \pa_z T^\eps_\Theta)(t') \Vert_{\mathcal{B}^{\frac{1}{2}}} + \Vert \pa_z (u_{\phi},T_\phi) (t') \Vert_{\mathcal{B}^{\frac{1}{2}}} \right) dt' = \rho(t) + \tau(t),
\end{align*}
where $\rho$ and $\tau$ are defined in the previous sections. We remark that, under the hypotheses of Theorems \ref{th:hydrolimit} and \ref{th:primitive}, we have 
\begin{align*}
	\frac{a}3 \leq \varphi(t,\xi) \leq \min\set{\phi(t,\xi),\Theta(t,\xi)}.
\end{align*}

Before giving the estimates of $G_1^q$, $G_2^q$ and $G_5^q$, we remark that we need a slightly different version of Lemmas \ref{lem:uww} and \ref{lem:vww} where we can relax the condition
$f(t) \geq \norm{\dd_z w_\phi(t)}_{\mathcal{B}^{\frac{1}{2}}}$. 
%Moreover, we also need controls for new terms of the type $\psca{w\dd u, \overline{w}}$ and $\psca{w\dd v, \overline{w}}$. 
More precisely, we will prove the following lemma.
\begin{lemma}
	\label{lem:nonli-conv}
	Let $0 < s < \frac{1}2$ and $\phi: \RR_+\times \RR \to \RR_+$. There exists a constant $C \geq 1$ such that, for any $(u,v,w,\overline{w})$, which are defined on $\RR_+\times \mathcal{S}$, $(u,w,\overline{w})\vert_{\dd\Scal} = 0$ and satisfy, for any $t\geq 0$, 
	\begin{equation*}
		u_\phi(t) \in \mathcal{B}^{\frac{3}{2}}, \quad \dd_z u_\phi(t), \, \dd_z w_\phi(t) \in \mathcal{B}^{\frac{1}{2}}, \quad \int_0^1 \pa_x u \, dz = 0, \quad \mbox{and} \quad u, \, w, \, \overline{w} \in \tilde{L}^2_{t,f(t)} (\mathcal{B}^{s+\frac{1}{2}}),
	\end{equation*}
	with 
	\begin{align*}
		f \in L^1_{\tiny loc}(\RR_+), \quad f(t) \geq \norm{\dd_z u_\phi}_{\mathcal{B}^{\frac{1}{2}}} \quad \mbox{and} \quad v(t,x,z) = -\int_0^z \pa_x u(t,x,z') dz',
	\end{align*} 
	we have, for any $\mathcal{R} \geq 0$ and for any $q\in \ZZ$,
	\begin{align}
		\label{eq:udxww} \int_0^t \abs{\psca{e^{\Rcal t'} \Delta_q^h (u\p_x w)_{\phi}, e^{\Rcal t'} \Delta_q^h \overline{w}_{\phi}}_{L^2}} dt'
		\leq C d_q^2 2^{-2qs} \norm{e^{\mathcal{R}t} w_{\phi}}_{\tilde{L}^2_{t,f(t)}(\mathcal{B}^{s+\frac{1}{2}})} \norm{e^{\mathcal{R}t} \overline{w}_{\phi}}_{\tilde{L}^2_{t,f(t)}(\mathcal{B}^{s+\frac{1}{2}})},
	\end{align}
	\begin{multline}
		\label{eq:vdzww}
		\int_0^t \abs{\psca{e^{\mathcal{R}t'} \Delta_q^h (v\p_z w)_{\phi}, e^{\mathcal{R}t'} \Delta_q^h \overline{w}_{\phi}}_{L^2}} dt'
		\\ 
		\leq C d_q^2 2^{-2qs} \norm{u_\phi}_{L^\infty_t(\mathcal{B}^{\frac{3}{2}})}^{\frac{1}2} \norm{e^{\mathcal{R}t} \dd_z w_{\phi}}_{\tilde{L}^2_t(\mathcal{B}^s)} \norm{e^{\mathcal{R}t} \overline{w}_{\phi}}_{\tilde{L}^2_{t,f(t)}(\mathcal{B}^{s+\frac{1}{2}})},
	\end{multline}
	where $(d_q)_{q\in\ZZ}$ is a positive sequence with $\sum_{q\in\ZZ} d_q = 1$.
\end{lemma}

\begin{proof}[Proof of Lemma \ref{lem:nonli-conv}]
	
	\mbox{}
	
	1. In order to prove Estimate \eqref{eq:udxww}, we can use Lemma \ref{lem:uww} and we only need to modify the calculations that we did for the $A_{2,q}$ on the page \pageref{eq:A2q}. We use Sobolev inclusion $\dot{H}^1_z([0,1]) \hookrightarrow L^\infty_z([0,1])$ and write
	\begin{align*}
		A_{2,q} &\lesssim \sum_{|q-q'| \leq 4} \int_0^t e^{2\Rcal t'} \Vert S^h_{q'-1} \partial_xw_{\phi} \Vert_{L^{\infty}_xL^2_z} \Vert \Delta_{q'}^h \dd_z u_{\phi} \Vert_{L^2} \Vert \Delta_q^h \overline{w}_{\phi} \Vert_{L^2} dt'.
	\end{align*}
	Now, using the definition of $S^h_{q'-1}$ and of the $\Bcal^{\frac{1}2}$-norm and Bernstein lemma \ref{le:Bernstein}, we get
	\begin{align*}
		A_{2,q} &\lesssim \sum_{|q-q'| \leq 4} \int_0^t e^{2\Rcal t'} \pare{\sum_{l \leq q'-2} 2^{\frac{3l}2} \Vert \Delta^h_l w_{\phi} \Vert_{L^2}} d_{q'}(\dd_z u_\phi) 2^{-\frac{q'}2} \Vert \dd_z u_{\phi} \Vert_{\Bcal^{\frac{1}2}} \Vert \Delta_q^h \overline{w}_{\phi} \Vert_{L^2} dt'.
	\end{align*}
	Since $d_{q'}(\dd_z u_\phi) \leq 1$, using Proposition \ref{prop:dqweight}, we can write
	\begin{align*}
		A_{2,q} &\lesssim \sum_{|q-q'| \leq 4} 2^{-\frac{q'}2} \pare{\int_0^t \sum_{l \leq q'-2} 2^{3l} \Vert \dd_z u_{\phi} \Vert_{\Bcal^{\frac{1}2}} \Vert \Delta^h_l e^{\Rcal t'} w_{\phi} \Vert_{L^2}^2 \ dt'}^{\frac{1}2} \pare{\int_0^t \Vert \dd_z u_{\phi}  \Vert_{\Bcal{\frac{1}2}} \Vert \Delta_q^h e^{\Rcal t'} \overline{w}_{\phi} \Vert_{L^2}^2 \ dt'}^{\frac{1}2} \\
		&\lesssim \sum_{|q-q'| \leq 4} 2^{-\frac{q'}2} \pare{\sum_{l \leq q'-2} d_l(w) 2^{l(1-s)}} \norm{e^{\mathcal{R}t} w_{\phi}}_{\tilde{L}^2_{t,f(t)}(\mathcal{B}^{s+\frac{1}{2}})} d_q(\overline{w}) 2^{-q(s+\frac{1}2)} \norm{e^{\mathcal{R}t} \overline{w}_{\phi}}_{\tilde{L}^2_{t,f(t)}(\mathcal{B}^{s+\frac{1}{2}})}.
	\end{align*}
	We remark that the sequence $\set{\overline{d}_{q'}}_{q'\in\ZZ}$ with $$\overline{d}_{q'} = \sum_{l \leq q'-2} d_l(w) 2^{(q'-l)(s-1)}$$
	can be written as a convolution product of two summable sequences if $0 < s < 1$. Thus, we finally obtain
	\begin{align*}
		A_{2,q} &\lesssim d_q^2 2^{-2qs} \norm{e^{\mathcal{R}t} w_{\phi}}_{\tilde{L}^2_{t,f(t)}(\mathcal{B}^{s+\frac{1}{2}})} \norm{e^{\mathcal{R}t} \overline{w}_{\phi}}_{\tilde{L}^2_{t,f(t)}(\mathcal{B}^{s+\frac{1}{2}})},
	\end{align*}
	where $$d_q^2 = d_q(\overline{w}) \sum_{|q-q'| \leq 4} \overline{d}_{q'} 2^{(q-q')(s-\frac{1}2)}.$$
	
	\medskip
	
	2. To prove Estimate \eqref{eq:vdzww}, the same modifications can be done to the terms $B_{2,q}$ and $B_{3,q}$ in the proof of Lemma \ref{lem:vww} on page \pageref{eq:B2q}. We will show these modifications for $B_{2,q}$. We recall that, using Lemma \ref{lem:tranreg}, we have
	\begin{align*}
		B_{2,q} &\lesssim \sum_{|q'-q|\leq 4} \int_0^t e^{2\mathcal{R}t'} \Vert S_{q'-1}^h \pa_z w_{\phi} \Vert_{L^{\infty}_xL^2_z} \Vert \Delta_{q'}^h v_{\phi} \Vert_{L^2_xL^{\infty}_z} \Vert \Delta_q^h \overline{w}_{\phi} \Vert_{L^2} dt' \\
		&\lesssim  \sum_{|q'-q|\leq 4} \int_0^t \sum_{l\leq q'-2} 2^{\frac{l}2} \Vert \Delta^h_l e^{\mathcal{R}t'} \pa_z w_{\phi} \Vert_{L^2} 2^{q'} \Vert \Delta_{q'}^h u_{\phi} \Vert_{L^2} \Vert \Delta_q^h e^{\mathcal{R}t'} \overline{w}_{\phi} \Vert_{L^2} dt'.
	\end{align*}
	Now, using Poincar\'e inequality, we can write
	\begin{align*}
		2^{q'} \Vert \Delta_{q'}^h u_{\phi} \Vert_{L^2} = 2^{\frac{3q'}4} \Vert \Delta_{q'}^h u_{\phi} \Vert_{L^2}^{\frac{1}2} 2^{\frac{q'}4} \Vert \Delta_{q'}^h \dd_z u_{\phi} \Vert_{L^2}^{\frac{1}2} \leq \norm{u_\phi}_{\Bcal^{\frac{3}2}}^{\frac{1}2} \norm{\dd_z u_\phi}_{\Bcal^{\frac{1}2}}^{\frac{1}2}.
	\end{align*}
	Then, Cauchy-Schwartz inequality implies
	\begin{align*}
		B_{2,q} &\lesssim \sum_{|q'-q|\leq 4} \int_0^t \norm{u_\phi}_{\Bcal^{\frac{3}2}}^{\frac{1}2} \sum_{l\leq q'-2} 2^{\frac{l}2} \Vert \Delta^h_l e^{\mathcal{R}t'} \pa_z w_{\phi} \Vert_{L^2} \pare{\norm{\dd_z u_\phi}_{\Bcal^{\frac{1}2}}^{\frac{1}2} \Vert \Delta_q^h e^{\mathcal{R}t'} \overline{w}_{\phi} \Vert_{L^2}} dt' \\
		&\lesssim \norm{u_\phi}_{L^\infty_t(\Bcal^{\frac{3}2})}^{\frac{1}2} \sum_{l\leq q'-2} 2^{\frac{l}2} \pare{\int_0^t \Vert \Delta^h_l e^{\mathcal{R}t'} \pa_z w_{\phi} \Vert_{L^2}^2 \ dt'}^{\frac{1}2} \pare{\int_0^t \norm{\dd_z u_\phi}_{\Bcal^{\frac{1}2}} \Vert \Delta_q^h e^{\mathcal{R}t'} \overline{w}_{\phi} \Vert_{L^2}^2 \ dt'}^{\frac{1}2} \\
		&\lesssim \norm{u_\phi}_{L^\infty_t(\Bcal^{\frac{3}2})}^{\frac{1}2} \pare{\sum_{l \leq q'-2} d_l(w) 2^{l(\frac{1}2-s)}} \norm{e^{\Rcal t} \dd_z w_\phi}_{\tilde{L}^2_t(\Bcal^s)} \pare{d_q(\overline{w}) 2^{-q(s+\frac{1}2)} \norm{e^{\mathcal{R}t} \overline{w}_{\phi}}_{\tilde{L}^2_{t,f(t)}(\mathcal{B}^{s+\frac{1}{2}})}}
	\end{align*}
	We remark that the sequence $\set{\overline{d}_{q'}}_{q'\in\ZZ}$ with $$\underline{d}_{q'} = \sum_{l \leq q'-2} d_l(w) 2^{(q'-l)(s-\frac{1}2)}$$
	can be written as a convolution product of two summable sequences if $0 < s < \frac{1}2$. Thus, we obtain
	\begin{align*}
		B_{2,q} &\lesssim d_q^2 2^{-2qs} \norm{u_\phi}_{L^\infty_t(\Bcal^{\frac{3}2})}^{\frac{1}2} \norm{e^{\Rcal t} \dd_z w_\phi}_{\tilde{L}^2_t(\Bcal^s)} \norm{e^{\mathcal{R}t} \overline{w}_{\phi}}_{\tilde{L}^2_{t,f(t)}(\mathcal{B}^{s+\frac{1}{2}})}.
	\end{align*}
	Finally, we remark that the term $B_{3,q}$ can be controlled in the similar way and we obtain Estimate \eqref{eq:vdzww}.		
\end{proof}

\subsection{Control of $G_1^q$}

We start by observing that we can write
\begin{align*}
	R^1_{\varphi} = (\eps^2 \pa_x^2u)_{\varphi} - (u^\eps \pa_x w^1)_{\varphi} - (w^1 \pa_x u)_{\varphi} - (v^\eps\pa_zw^1)_{\varphi} - (w^2 \pa_z u)_{\varphi},
\end{align*}     
so,   
\begin{align*}
	G_1^q \leq I_1^q + I_2^q + I_3^q + I_4^q + I_5^q,
\end{align*}
where $I_i^q$ will be precised and controlled in what follows.

Using Bernstein lemma \ref{le:Bernstein} and Poincar\'e inequality, we can write
\begin{equation*}
	\norm{\DD^h_q (\dd_xu)_\varphi}_{L^2} \lesssim 2^{-q} \norm{\DD^h_q u_\varphi}_{L^2} \lesssim 2^{-q} \norm{\DD^h_q \dd_z u_\varphi}_{L^2}
\end{equation*}
Then, using Poincar\'e inequality and Bernstein lemma \ref{le:Bernstein}, we have 
\begin{align}
	\label{eq:I1}
	I_1^q &= 2\eps^2 \int_0^t \abs{ \psca{\Delta_q^h (\pa_x^2 u)_{\varphi}, \Delta_q^h w^1_{\varphi}}_{L^2}} dt' = 2\eps^2 \int_0^t  \norm{\Delta_q^h (\pa_x^2 u)_{\varphi}}_{L^2} \norm{\Delta_q^h w^1_{\varphi}}_{L^2} dt' 
	\\
	&\lesssim \eps^2 2^{2q} \norm{\DD^h_q u_\varphi}_{L^2_t(L^2)} \norm{\Delta_q^h w^1_{\varphi}}_{L^2_t(L^2)} \lesssim \eps^2 2^{2q} \norm{\DD^h_q \dd_z u_\varphi}_{L^2_t(L^2)} \norm{\Delta_q^h w^1_{\varphi}}_{L^2_t(L^2)}\notag
	\\
	&\lesssim d_q^2 2^{-2qs} \eps \Vert \pa_z u_{\varphi} \Vert_{\tilde{L}^2_t(\mathcal{B}^{s+\frac{3}2})} \Vert \eps w^1_{\varphi} \Vert_{\tilde{L}^2_t(\mathcal{B}^{s+\frac{1}2})} \lesssim d_q^2 2^{-2qs} \pare{C\eps^2 \Vert \pa_z u_{\varphi} \Vert_{\tilde{L}^2_t(\mathcal{B}^{s+\frac{3}2})}^2 + \frac{1}4 \Vert \eps w^1_{\varphi} \Vert_{\tilde{L}^2_t(\mathcal{B}^{s+\frac{1}2})}^2} \notag
\end{align}
%{\color{red}
	%Summing with respect to $q \in \mathbb{Z}$, we obtain
	%\begin{align*} 
	%	\sum_{q\in\mathbb{Z}} 2^{\frac{q}{2}} \sqrt{I_1^q} &= \eps \sum_{q\in\mathbb{Z}} 2^{\frac{q}{2}}  \sqrt{\int_0^t \abs{ \psca{\Delta_q^h (\pa_x^2 u)_{\varphi}, \Delta_q^h w^1_{\varphi}}_{L^2}} dt'} \lesssim C\eps \Vert \pa_z u_{\varphi} \Vert_{\tilde{L}^2_t(\mathcal{B}^{\frac{3}{2}})} + C \Vert \eps w^1_{\varphi} \Vert_{\tilde{L}^2_t(\mathcal{B}^{\frac{3}{2}})}. 
	%\end{align*}
	%}

For $I_2^q$, Estimate \eqref{eq:udxww} of Lemma \ref{lem:nonli-conv} implies
\begin{align}
	\label{eq:I21}
	I_{2}^q =  2\int_0^t \abs{ \psca{\Delta_q^h(u^\eps \pa_x w^1)_{\varphi},\Delta_q^h w^1_{\varphi}}_{L^2}} dt' \lesssim d_q^2 2^{-2qs}  \Vert w^1_{\varphi} \Vert_{\tilde{L}^2_{t,\dot{\eta}(t)}(\mathcal{B}^{s+\frac{1}2})}^2,
\end{align}
%{\color{red}
	%since $u$ belongs to $\tilde{L}^\infty(B^{1/2})$, so our sequence $d_q$ not depending on time, then
	%\begin{align*}
	%	\sum_{q\in \mathbb{Z}} 2^\frac{q}{2} \sqrt{I_{21}^q} =  \sum_{q\in \mathbb{Z}} 2^\frac{q}{2} \sqrt{\int_0^t \abs{ \psca{\Delta_q^h(u^\eps \pa_x w^1)_{\varphi},\Delta_q^h w^1_{\varphi}}_{L^2}} dt'} \leq C  \Vert w^1_{\varphi} \Vert_{\tilde{L}^2_{t,\dot{\eta}(t)}(\mathcal{B}^{1})},
	%\end{align*}
	%}

The term $I_{3}^q$ can be controlled in the exact same way as we did to prove Estimate \eqref{eq:vdzww}. Indeed, using the Bony decomposition, we can write 
\begin{align*}
	\int_0^t \abs{\psca{\Delta_q^h (w^1\p_x u)_{\varphi}, \Delta_q^h w^1_{\varphi}}_{L^2}} dt'  \leq \overline{B}_{1,q} + \overline{B}_{2,q} + \overline{B}_{3,q},
\end{align*}
with 
\begin{align*}
	\overline{B}_{1,q} &= \int_0^t \abs{\psca{\Delta_q^h ( \Tcal^h_{w^1}\pa_x u)_{\varphi}, \Delta_q^h w^1_{\varphi}}_{L^2}} dt'\\
	\overline{B}_{2,q} &= \int_0^t \abs{\psca{\Delta_q^h ( \Tcal^h_{\pa_xu} w^1)_{\varphi}, \Delta_q^h w^1_{\varphi}}_{L^2}} dt'\\
	\overline{B}_{3,q} &= \int_0^t \abs{\psca{\Delta_q^h (\Rcal^h(w^1,\pa_xu))_{\varphi}, \Delta_q^h w^1_{\varphi}}_{L^2}} dt'.
\end{align*}
Using Poincar\'e inequality and Bernstein lemma \ref{le:Bernstein}, we can write
\begin{align*}
	\overline{B}_{1,q} &\lesssim \sum_{|q'-q|\leq 4} \int_0^t \Vert S_{q'-1}^h w^1_{\varphi} \Vert_{L^\infty} \Vert \Delta_{q'}^h \dd_x u_{\varphi} \Vert_{L^2} \Vert \Delta_q^h w^1_{\varphi} \Vert_{L^2} dt'
	\\
	&\lesssim  \sum_{|q'-q|\leq 4} \int_0^t \sum_{l\leq q'-2} 2^{\frac{l}2} \Vert \Delta^h_l \pa_z w^1_{\varphi} \Vert_{L^2} 2^{q'} \Vert \Delta_{q'}^h u_{\varphi} \Vert_{L^2} \Vert \Delta_q^h w^1_{\varphi} \Vert_{L^2} dt'.
\end{align*}
Thus, we can control $\overline{B}_{1,q}$ in the exact same way as what we did for $B_{2,q}$ above. Similarly, $\overline{B}_{2,q}$ can be controlled in the same way as $B_{1,q}$ and $\overline{B}_{3,q}$ as $B_{3,q}$. We obtain
\begin{align}
	\label{eq:I3} I_{3}^q &= 2\int_0^t \abs{\psca{\Delta_q^h (w^1\p_x u)_{\varphi}, \Delta_q^h w^1_{\varphi}}_{L^2}} dt' \lesssim d_q^2 2^{-2qs} \norm{u_\varphi}_{L^\infty_t(\Bcal^{\frac{3}2})}^{\frac{1}2} \norm{\dd_z w^1_\varphi}_{\tilde{L}^2_t(\Bcal^s)} \norm{w^1_{\varphi}}_{\tilde{L}^2_{t,\dot{\eta}(t)}(\mathcal{B}^{s+\frac{1}{2}})} \\
	&\lesssim d_q^2 2^{-2qs} \pare{\norm{u_\varphi}_{L^\infty_t(\Bcal^{\frac{3}2})} \norm{\dd_z w^1_\varphi}_{\tilde{L}^2_t(\Bcal^s)}^2 + \norm{w^1_{\varphi}}_{\tilde{L}^2_{t,\dot{\eta}(t)}(\mathcal{B}^{s+\frac{1}{2}})}^2}. \notag
\end{align}

For $I_4^q$, using Estimate \eqref{eq:vdzww} of Lemma \ref{lem:nonli-conv}, we also have
\begin{align}
	\label{eq:I42}
	I_4^q &= 2\int_0^t \abs{\psca{ \Delta_q^h(v^\eps\pa_z w^1)_\varphi , \Delta_q^h w^1_\varphi }_{L^2}} dt' 
	\\
	&\lesssim d_q^2 2^{-2qs} \norm{u^\eps_\varphi}_{L^\infty_t(\mathcal{B}^{\frac{3}{2}})}^{\frac{1}2} \norm{\dd_z w^1_{\varphi}}_{\tilde{L}^2_t(\mathcal{B}^s)} \norm{w^1_{\varphi}}_{\tilde{L}^2_{t,\dot{\eta}(t)}(\mathcal{B}^{s+\frac{1}{2}})} \notag
	\\
	&\lesssim d_q^2 2^{-2qs} \pare{\norm{u^\eps_\varphi}_{L^\infty_t(\mathcal{B}^{\frac{3}{2}})}  \norm{\dd_z w^1_\varphi}_{\tilde{L}^2_t(\Bcal^s)}^2 + \norm{w^1_{\varphi}}_{\tilde{L}^2_{t,\dot{\eta}(t)}(\mathcal{B}^{s+\frac{1}{2}})}^2}. \notag
\end{align}

For the term $I_5^q$, we will not use Lemma \ref{lem:nonli-conv} and we will use a slightly different control. We first recall the following horizontal decomposition into paraproducts and remainders
\begin{align*}
	I_5^q &= 2\int_0^t \abs{\psca{ \Delta_q^h(w^2 \pa_z u)_\varphi , \Delta_q^h w^1_\varphi }_{L^2}} dt'
	\\
	&= 2\int_0^t \abs{\psca{\Delta_q^h( T_{w^2}^h\pa_z u + T^h_{\pa_zu} w^2 + R^h(w^2,\pa_z u))_\varphi, \Delta_q^h w^1_\varphi}_{L^2}} dt'
	\\ 
	&\leq  2 \pare{I_{5,1}^q + I_{5,2}^q+I_{5,3}^q},
\end{align*}
where
\begin{align*}
	I_{5,1}^q &= \int_0^t \abs{\psca{\Delta_q^h( T_{w^2}^h \pa_zu)_\varphi, \Delta_q^h w^1_\varphi}_{L^2}} dt' \\
	I_{5,2}^q &= \int_0^t \abs{\psca{\Delta_q^h(  T^h_{\pa_zu} w^2)_\varphi, \Delta_q^h w^1_\varphi}_{L^2}} dt' \\
	I_{5,3}^q &= \int_0^t \abs{\psca{\Delta_q^h( R^h(w^2,\pa_z u))_\varphi, \Delta_q^h w^1_\varphi}_{L^2}} dt'.
\end{align*}
For $I_{5,1}^q$, Lemma \ref{lem:tranreg} yields
\begin{align*}
	\Vert S_{q'-1}^h w^2_{\varphi} \Vert_{L^\infty} \lesssim \sum_{l\leq q'-2} \Vert \Delta^h_l w^1_\varphi \Vert_{L^\infty} \lesssim \sum_{l\leq q'-2} 2^{\frac{3l}{2}} \Vert \Delta^h_l w^1_\varphi \Vert_{L^2}.
\end{align*}
Then, we can write
\begin{align*}
	I_{5,1}^q & \lesssim \sum_{|q'-q|\leq 4} \int_0^t \Vert S_{q'-1}^h w^2_\varphi \Vert_{L^2} \Vert \Delta_{q'}^h \pa_z u_\varphi \Vert_{L^2} \Vert \Delta_q^h w^1_\varphi \Vert_{L^2} dt' \\
	&\lesssim \sum_{|q'-q|\leq 4} \int_0^t \pare{\sum_{l\leq q'-2} 2^{\frac{3l}{2}} \Vert \Delta^h_l w^1_\varphi \Vert_{L^2}} \Vert \Delta_{q'}^h \pa_z u_\varphi \Vert_{L^2} \Vert \Delta_q^h w^1_\varphi \Vert_{L^2} dt' .
\end{align*}
Thus, $I_{5,1}^q$ can be controlled in the same way as $A_{2,q}$ in the proof of Estimate \eqref{eq:udxww} of Lemma \ref{lem:nonli-conv} and if $0 < s < 1$, we have
\begin{align}
	\label{eq:I51}
	I_{5,1}^q \lesssim d_q^2 2^{-2qs} \Vert w^1_{\varphi} \Vert_{\tilde{L}^2_{t,\dot{\eta}(t)}(\mathcal{B}^{s+\frac{1}2})}^2.
\end{align}
For $I_{5,2}^q$, using Lemma \ref{lem:tranreg}, we write
\begin{align*}
	\norm{\Delta_{q'}^h w^2_\varphi}_{L^2_xL^\infty_z} \lesssim 2^{q'} \norm{\Delta_{q'}^h w^1_\varphi}_{L^2}.
\end{align*}
Since
\begin{align*}
	\Vert S_{q'-1}^h \pa_z u_\varphi \Vert_{L^\infty_xL^2_z} \lesssim \norm{\pa_z u_\varphi}_{\Bcal^{\frac{1}2}},
\end{align*}
as in the estimate of the term $A_{1,q}$ in the proof of Lemma \ref{lem:uww}, we have,
\begin{align}
	\label{eq:I52}
	I_{5,2}^q &\lesssim \sum_{|q'-q|\leq 4} \int_0^t \Vert S_{q'-1}^h \pa_z u_\varphi \Vert_{L^\infty_xL^2_z} \Vert \Delta_{q'}^h w^2_\varphi \Vert_{L^2_xL^\infty_z} \Vert \Delta_q^h w^1_\varphi \Vert_{L^2} dt' 
	\\
	&\lesssim \sum_{|q'-q|\leq 4} \int_0^t \Vert \pa_z u_\varphi \Vert_{\mathcal{B}^{\f12}} 2^{q'} \Vert \Delta_{q'}^h w^1_\varphi \Vert_{L^2} \Vert \Delta_q^h w^1_\varphi \Vert_{L^2}dt' \notag
	\\
	&\lesssim  \sum_{|q'-q|\leq 4} 2^{q'} \left( \int_0^t \Vert \pa_z u_\varphi \Vert_{\mathcal{B}^{\f12}}  \Vert \Delta_{q'}^h w^1_\varphi \Vert_{L^2}^2 dt' \right)^\f12  \left( \int_0^t \Vert \pa_z u_\varphi \Vert_{\mathcal{B}^{\f12}}  \Vert \Delta_{q}^h w^1_\varphi \Vert_{L^2}^2 dt' \right)^\f12 \notag
	\\
	&\lesssim d_q^2 2^{-2qs} \Vert w^1_\varphi \Vert_{\tilde{L}^2_{t,\dot{\eta}(t)}(\mathcal{B}^{s+\frac{1}2})}^2. \notag
\end{align}
The last term $I_{5,3}^q$ can also be treated in the similar way as the term $A_{3,q}$  in the proof of Lemma \ref{lem:uww} and we will get
\begin{align} \label{eq:I5,3}
	I_{5,3}^q \lesssim d_q^2 2^{-2qs} \Vert w^1_\varphi \Vert_{\tilde{L}^2_{t,\dot{\eta}(t)}(\mathcal{B}^{s+\frac{1}2})}^2.
\end{align}
Summing Estimate \eqref{eq:I51}, \eqref{eq:I52} and \eqref{eq:I5,3} will imply that 
\begin{align}
	\label{eq:I5}
	I_5^q \lesssim d_q^2 2^{-2qs} \Vert w^1_\varphi \Vert_{\tilde{L}^2_{t,\dot{\eta}(t)}(\mathcal{B}^{s+\frac{1}2})}^2.
\end{align}

\medskip

Summing Estimates \eqref{eq:I1}, \eqref{eq:I21}, \eqref{eq:I3}, \eqref{eq:I42}, and \eqref{eq:I5}, we obtain 
\begin{align} \label{eq:G1}
	G_1^q \lesssim d_q^2 2^{-q} \Big(C \eps^2 \Vert \pa_z u_{\varphi} \Vert_{\tilde{L}^2_t(\mathcal{B}^{s+\frac{3}2})}^2 + \frac{1}4 \Vert \eps w^1_{\varphi} \Vert_{\tilde{L}^2_t(\mathcal{B}^{s+\frac{1}2})}^2
	+ C\norm{u^\eps_\varphi}_{L^\infty_t(\mathcal{B}^{\frac{3}{2}})}\norm{\dd_z w^1_\varphi}_{\tilde{L}^2_t(\Bcal^s)}^2 + \norm{w^1_{\varphi}}_{\tilde{L}^2_{t,\dot{\eta}(t)}(\mathcal{B}^{s+\frac{1}{2}})}^2 \Big).
\end{align}

\medskip

\subsection{Control of $G_2^q$}

We recall that
\begin{align*}
	R^2_\varphi &= - \eps^2\left(\pa_t v  -\pa_z^2 v -\eps^2 \pa_x^2 v + u^{\eps} \pa_x v^{\eps} + v^\eps \pa_z v^\eps \right)_\varphi.
\end{align*}
Then, we can bound $G_2^q$ as follows
\begin{align*}
	G_2^q \leq J_1^q + J_2^q + J_3^q + J_4^q + J_5^q,
\end{align*}
where $J_i^q$ will be precised and controlled in what follows.

Using Young inequality, Sobolev inclusion $\dot{H}^1_z([0,1]) \hookrightarrow L^\infty_z([0,1])$ and Bernstein lemma \ref{le:Bernstein}, we have
\begin{align*}
	J_1^q &=2\eps^2 \int_0^t \abs{\psca{\Delta_q^h(\pa_t v)_\varphi ,\Delta_q^h w^2_\varphi}_{L^2}} dt' \lesssim \eps^2 \Vert \Delta^h_q(\pa_t u)_\varphi \Vert_{L^2_t(L^\infty_xL^2_z)} \Vert w^2_\varphi \Vert_{L^2_t(L^2_xL^\infty_z)}
	\\
	& \lesssim \eps^2 2^{\frac{q}2} \Vert \Delta^h_q(\pa_t v)_\varphi \Vert_{L^2_t(L^2)} \Vert \Delta^h_q \dd_z w^2_\varphi \Vert_{L^2_t(L^2)} \lesssim \eps^2 2^{\frac{q}2} \Vert \Delta^h_q(\pa_t \dd_z v)_\varphi \Vert_{L^2_t(L^2)} \Vert \Delta^h_q \dd_z w^2_\varphi \Vert_{L^2_t(L^2)}
	\\ 
	&\lesssim \eps^2 2^{\frac{q}2} \Vert \Delta^h_q(\pa_t \dd_x u)_\varphi \Vert_{L^2_t(L^2)} \Vert \Delta^h_q \dd_z w^2_\varphi \Vert_{L^2_t(L^2)} \lesssim  \eps^2 2^{\frac{3q}2} \Vert \Delta^h_q(\pa_t u)_\varphi \Vert_{L^2_t(L^2)} \Vert \Delta^h_q \dd_z w^2_\varphi \Vert_{L^2_t(L^2)}	\\
	&\lesssim d_q^2 2^{-2qs} \eps^2 \pare{C\Vert (\pa_t u)_\varphi \Vert_{\tilde{L}^2_t(\Bcal^{s + \frac{3}2})}^2 + \frac{1}{100} \Vert \dd_z w^2_\varphi \Vert_{\tilde{L}^2_t(\Bcal^{s})}^2}.
\end{align*}

Similarly, we have
\begin{align*}
	J_2^q &= 2\eps^2 \int_0^t \abs{\psca{\Delta_q^h(\pa_z^2 v)_\varphi ,\Delta_q^h w^2_\varphi}_{L^2}} dt' \lesssim d_q^2 2^{-2qs} \eps^2 \pare{C\Vert \pa_z u_\varphi \Vert_{\tilde{L}^2_t(\Bcal^{s + \frac{3}2})}^2 + \frac{1}{100} \Vert \dd_z w^2_\varphi \Vert_{\tilde{L}^2_t(\Bcal^{s})}^2},
\end{align*}
and
\begin{align*}
	J_3^q &= 2\eps^4 \int_0^t \abs{\psca{ \Delta_q^h(\pa_x^2 v)_\varphi ,\Delta_q^h w^2_\varphi}_{L^2}}dt' \lesssim d_q^2 2^{-2qs} \eps^4 \pare{C\Vert \pa_z u_\varphi \Vert_{\tilde{L}^2_t(\Bcal^{s + \frac{5}2})}^2 + \frac{1}{100} \Vert \dd_z w^2_\varphi \Vert_{\tilde{L}^2_t(\Bcal^{s})}^2}.
\end{align*}

For $J_4^q$, we can use similar estimates as in the proof of Estimate \eqref{eq:vdzww} of Lemma \ref{lem:nonli-conv}. We first apply Bony's decomposition for the horizontal variable and write
\begin{align*}
	J_{4}^q = 2 \eps^2\int_0^t \abs{\psca{\Delta_q^h(u^\eps \pa_x v^\eps)_{\varphi},\Delta_q^h w^2_\varphi}_{L^2}} dt' \leq 2\pare{J_{41}^q + J_{42}^q + J_{43}^q},
\end{align*}
with
\begin{align*}
	J_{41}^q &= \eps^2 \int_0^t \abs{\psca{\Delta_q^h(T^h_{u^\eps} \pa_x v^\eps)_{\varphi},\Delta_q^h w^2_\varphi}_{L^2}} dt' \\
	J_{42}^q &= \eps^2 \int_0^t \abs{\psca{\Delta_q^h(T^h_{\pa_x v^\eps} u^\eps)_{\varphi},\Delta_q^h w^2_\varphi}_{L^2}} dt'\\
	J_{43}^q &= \eps^2 \int_0^t \abs{\psca{\Delta_q^h(R^h(u^\eps,\pa_x v^\eps))_{\varphi},\Delta_q^h w^2_\varphi}_{L^2}} dt'.
\end{align*}
We remark that Lemma \ref{lem:tranreg} and Bernstein lemma \ref{le:Bernstein} imply
\begin{align*}
	\Vert \Delta_{q'}^h \pa_x v^\eps_\varphi \Vert_{L^2_xL^\infty_z} \leq 2^{2q'} \norm{\Delta^h_q u^\eps_\varphi}_{L^2}.
\end{align*}
We also remark that, using Poincar\'e inequality, Bernstein lemma \ref{le:Bernstein} and then Cauchy-Schwartz inequality, we have
\begin{align*}
	\Vert S_{q'-1}^h u^\eps_\varphi \Vert_{L^\infty_xL^2_z} \lesssim \sum_{l \leq q'-2} 2^{\frac{l}2} \norm{\Delta^h_l u^\eps_\varphi}_{L^2} \lesssim \sum_{l \leq q'-2} 2^{\frac{l}2} \norm{\Delta^h_l u^\eps_\varphi}_{L^2}^{\frac{1}2} \norm{\Delta^h_l \dd_z u^\eps_\varphi}_{L^2}^{\frac{1}2} \lesssim \Vert u_\varphi^\eps \Vert_{\mathcal{B}^{\f12}}^{\frac{1}2} \Vert \dd_z u_\varphi^\eps \Vert_{\mathcal{B}^{\f12}}^{\frac{1}2}.
\end{align*}
Then, using Young inequality, we get
\begin{align*}
	J_{41}^q &\lesssim \eps^2 \sum_{|q'-q| \leq 4} \int_0^t  \Vert S_{q'-1}^h u^\eps_\varphi \Vert_{L^\infty_xL^2_z} \Vert \Delta_{q'}^h \pa_x v^\eps_\varphi \Vert_{L^2_xL^\infty_z} \Vert \Delta_q^h w^2_\varphi \Vert_{L^2} dt' \\
	&\lesssim \eps^2 \sum_{|q'-q| \leq 4} \int_0^t \Vert u_\varphi^\eps \Vert_{\mathcal{B}^{\f12}}^{\frac{1}2} \Vert \dd_z u_\varphi^\eps \Vert_{\mathcal{B}^{\f12}}^{\frac{1}2} 2^{2q'} \Vert \Delta_{q'}^h u^\eps_\varphi \Vert_{L^2} \Vert \Delta_q^h w^2_\varphi \Vert_{L^2} dt'.
\end{align*}
Cauchy-Schwartz and Poincar\'e inequalities finally imply
\begin{align*}
	J_{41}^q &\lesssim d_q^2 2^{-2qs} \eps^2 \pare{\Vert u_\varphi^\eps \Vert_{L^\infty(\mathcal{B}^{\f12})} \norm{\dd_z u_\varphi}_{\tilde{L}^2_t(\Bcal^{s+\frac{3}2})}^2 + \Vert w^2_\varphi \Vert_{\tilde{L}^2_{t,\dot{\eta}(t)}(\mathcal{B}^{s+\frac{1}2})}^2}
\end{align*}
Similar estimates also lead to
\begin{align*}
	J_{42}^q &\lesssim d_q^2 2^{-2qs} \eps^2 \pare{\Vert u_\varphi^\eps \Vert_{L^\infty(\mathcal{B}^{\f12})} \norm{\dd_z u_\varphi}_{\tilde{L}^2_t(\Bcal^{s+\frac{3}2})}^2 + \Vert w^2_\varphi \Vert_{\tilde{L}^2_{t,\dot{\eta}(t)}(\mathcal{B}^{s+\frac{1}2})}^2}
	\\
	J_{43}^q &\lesssim d_q^2 2^{-2qs} \eps^2 \pare{\Vert u_\varphi^\eps \Vert_{L^\infty(\mathcal{B}^{\f12})} \norm{\dd_z u_\varphi}_{\tilde{L}^2_t(\Bcal^{s+\frac{3}2})}^2 + \Vert w^2_\varphi \Vert_{\tilde{L}^2_{t,\dot{\eta}(t)}(\mathcal{B}^{s+\frac{1}2})}^2}.
\end{align*}
Summing the above inequalities, we obtain
\begin{align*}
	J_{4}^q &\lesssim d_q^2 2^{-2qs} \eps^2 \pare{\Vert u_\varphi^\eps \Vert_{L^\infty(\mathcal{B}^{\f12})} \norm{\dd_z u_\varphi}_{\tilde{L}^2_t(\Bcal^{s+\frac{3}2})}^2 + \Vert w^2_\varphi \Vert_{\tilde{L}^2_{t,\dot{\eta}(t)}(\mathcal{B}^{s+\frac{1}2})}^2}.
\end{align*}

Now, for $J_5^q$, using the divergence-free property, we can write $v^\eps \pa_z v^\eps = v^\eps \dd_x u^\eps$. Thus, $J_5^q$ can be controlled in the exact same way as $J_4^q$ and we get
\begin{align*}
	J_{5}^q &= 2 \eps^2\int_0^t \abs{\psca{\Delta_q^h(v^\eps \pa_z v^\eps)_{\varphi},\Delta_q^h w^2_\varphi}_{L^2}} dt' 
	\\
	&\lesssim d_q^2 2^{-2qs} \eps^2 \pare{\Vert u_\varphi^\eps \Vert_{L^\infty(\mathcal{B}^{\f12})} \norm{\dd_z u_\varphi}_{\tilde{L}^2_t(\Bcal^{s+\frac{3}2})}^2 + \Vert w^2_\varphi \Vert_{\tilde{L}^2_{t,\dot{\eta}(t)}(\mathcal{B}^{s+\frac{1}2})}^2}.
\end{align*}

We deduce that
\begin{multline}
	\label{eq:G2}
	G_2^q \lesssim d_q^2 2^{-2qs} \eps^2 \Big( C\Vert (\pa_t u)_\varphi \Vert_{\tilde{L}^2_t(\Bcal^{s + \frac{3}2})}^2 + C\Vert \pa_z u_\varphi \Vert_{\tilde{L}^2_t(\Bcal^{s + \frac{3}2})}^2 + C\eps^2\Vert \pa_z u_\varphi \Vert_{\tilde{L}^2_t(\Bcal^{s + \frac{5}2})}^2
	\\ 
	+ \frac{1}{20} \Vert \dd_z w^2_\varphi \Vert_{\tilde{L}^2_t(\Bcal^{s})}^2
	+ \Vert u_\varphi^\eps \Vert_{L^\infty(\mathcal{B}^{\f12})} \norm{\dd_z u_\varphi}_{\tilde{L}^2_t(\Bcal^{s+\frac{3}2})}^2 + \Vert w^2_\varphi \Vert_{\tilde{L}^2_{t,\dot{\eta}(t)}(\mathcal{B}^{s+\frac{1}2})}^2 \Big).
\end{multline}

\medskip

\subsection{Control of $G_5^q$}

To estimate the last term $G_5^q$, we will write 
\begin{align*}
	R^3_{\varphi} = -(u^\eps \pa_x \theta + w^1 \pa_x T)_{\varphi} - (v^\eps \pa_z \theta + w^2 \pa_z T)_{\varphi},
\end{align*}   
and so,   
\begin{align*}
	G_5^q = 2\int_0^t \abs{\psca{\Delta_q^h R^3_{\varphi}, \Delta_q^h \theta_{\varphi}}_{L^2}} dt'  \leq 2\pare{L_1^q + L_2^q + L_3^q + L_4^q},
\end{align*}
where
\begin{align*}
	L_1^q &= \int_0^t \abs{\psca{\Delta_q^h (u^\eps \pa_x \theta )_{\varphi}, \Delta_q^h \theta_{\varphi}}_{L^2}} dt'\\
	L_2^q &= \int_0^t \abs{\psca{ \Delta_q^h(w^1 \pa_x T )_{\varphi},\Delta_q^h \theta_{\varphi}}_{L^2}} dt'\\
	L_3^q &= \int_0^t \abs{\psca{ \Delta_q^h(v^\eps \pa_z\theta)_{\varphi}, \Delta_q^h \theta_{\varphi}}_{L^2}} dt' \\
	L_4^q &= \int_0^t \abs{\psca{ \Delta_q^h(w^2 \pa_z T)_{\varphi}, \Delta_q^h \theta_{\varphi}}_{L^2}} dt'.
\end{align*}

Using Lemma \ref{lem:nonli-conv}, we immediately obtain
\begin{align*}
	L_1^q \lesssim d_q^2 2^{-2qs} \Vert \theta_{\varphi} \Vert_{\tilde{L}^2_{t,\dot{\eta}(t)}(\mathcal{B}^{s+\frac{1}2})}^2
\end{align*}
and
\begin{align*}
	L_3^q \lesssim d_q^2 2^{-2qs} \pare{ \norm{u^\eps_\varphi}_{L^\infty_t(\mathcal{B}^{\frac{3}{2}})} \norm{\dd_z \theta_{\varphi}}_{\tilde{L}^2_t(\mathcal{B}^s)}^2 +  \norm{\theta_{\varphi}}_{\tilde{L}^2_{t,\dot{\eta}(t)}(\mathcal{B}^{s+\frac{1}{2}})}^2}.
\end{align*}
The term $L_2^q$ can be controlled in the same way as we did for the term $I_{3}^q$ (page \pageref{eq:I3}). We have
\begin{align*}
	L_2^q \lesssim d_q^2 2^{-2qs} \pare{\norm{T^\eps_\varphi}_{L^\infty_t(\mathcal{B}^{\frac{3}{2}})}  \norm{\dd_z w^1_\varphi}_{\tilde{L}^2_t(\Bcal^s)}^2 + \norm{\theta_{\varphi}}_{\tilde{L}^2_{t,\dot{\eta}(t)}(\mathcal{B}^{s+\frac{1}{2}})}^2}.
\end{align*}
Finally, the term $L_4^q$ can be controlled in the same way as we did for the term $I_{5,1}^q$ (page \pageref{eq:I51}) and we obtain
\begin{align*}
	L_4^q \lesssim d_q^2 2^{-2qs} \pare{\Vert w^1_{\varphi} \Vert_{\tilde{L}^2_{t,\dot{\eta}(t)}(\mathcal{B}^{s+\frac{1}2})}^2 + \Vert \theta_{\varphi} \Vert_{\tilde{L}^2_{t,\dot{\eta}(t)}(\mathcal{B}^{s+\frac{1}2})}^2}.
\end{align*}

Thus, we deduce that 
\begin{multline} 
	\label{eq:G5}
	G_5^q \lesssim d_q^2 2^{-2qs} \Big( C\norm{w^1_{\varphi}}_{\tilde{L}^2_{t,\dot{\eta}(t)}(\mathcal{B}^{s+ \frac{1}2})}^2
	+ C\norm{\theta_{\varphi}}_{\tilde{L}^2_{t,\dot{\eta}(t)}(\mathcal{B}^{s+\frac{1}2})}^2
	\\ 
	+ C\norm{T_{\varphi}}_{L^\infty_{t}(\mathcal{B}^{\frac{3}{2}})} \norm{\pa_z w^1_{\varphi}}_{\tilde{L}^2_{t}(\mathcal{B}^s)}^2 
	+ C\norm{u_{\varphi}}_{L^\infty_{t}(\mathcal{B}^{\frac{3}{2}})} \norm{\pa_z \theta_{\varphi}}_{\tilde{L}^2_{t}(\mathcal{B}^s)}^2 \Big).
\end{multline}

\bigskip

\noindent \emph{End of the proof of Theorem \ref{th:Convergence}.} We put Estimates \eqref{eq:G1}, \eqref{eq:G2}, \eqref{eq:G3}, \eqref{eq:G4} and \eqref{eq:G5} into Inequalities \eqref{S6eq5bis} and \eqref{S6eq6bis}, we multiply \eqref{S6eq6bis} by $2C$ and then we sum the obtained inequalities. Now, we recall that we can take the square-root of each term of the resulting inequality, with the cost of a constant multiplier which will be included in the generic constant $C$. Without loss of generality, we can always consider that $C \geq 2$. We obtain
\begin{align}
	\label{eq:conv01}
	&\Vert \Delta_q^h (w^1_\varphi,\eps w^2_\varphi,\theta_\varphi) \Vert_{L^\infty_t(L^2)} + \sqrt{\mu} \pare{\int_0^t \dot{\eta}(t')\norm{|D_x|^{\frac{1}2} \Delta_q^h (w^1_\varphi,\eps w^2_\varphi,\theta_\varphi)}_{L^2} dt'}^\f12
	\\ 
	&\qquad + \frac{3}4 \Vert \pa_z \Delta_q^h (w^1_{\varphi},\eps w^2_{\varphi},\theta_\varphi) \Vert_{L^2_t(L^2)}
	+ \eps 2^{\frac{q}2} \Vert \Delta_q^h (w^1_{\varphi},\eps w^2_{\varphi}) \Vert_{L^2_t(L^2)} + 2^{\frac{q}2} \Vert \Delta_q^h \theta_{\varphi} \Vert_{L^2_t(L^2)} \notag
	\\
	&\qquad \qquad \leq 2C \norm{\Delta_q^h (w^1_{\varphi},\eps w^2_{\varphi},\theta_\varphi)(0)}_{L^2} + d_q 2^{-qs} \Big(C \eps \Vert \pa_z u_{\varphi} \Vert_{\tilde{L}^2_t(\mathcal{B}^{s+\frac{3}2})} + \frac{1}2 \Vert \eps w^1_{\varphi} \Vert_{\tilde{L}^2_t(\mathcal{B}^{s+\frac{1}2})} \notag
	\\ 
	&\qquad \qquad + C\norm{u^\eps_\varphi}_{L^\infty_t(\mathcal{B}^{\frac{3}{2}})}^{\f12} \norm{\dd_z w^1_\varphi}_{\tilde{L}^2_t(\Bcal^s)} + C\eps \Vert (\pa_t u)_\varphi \Vert_{\tilde{L}^2_t(\Bcal^{s + \frac{3}2})} + C\eps \Vert \pa_z u_\varphi \Vert_{\tilde{L}^2_t(\Bcal^{s + \frac{3}2})} \notag
	\\
	&\qquad \qquad + C\eps^2\Vert \pa_z u_\varphi \Vert_{\tilde{L}^2_t(\Bcal^{s + \frac{5}2})} + \frac{1}{4} \eps \Vert \dd_z w^2_\varphi \Vert_{\tilde{L}^2_t(\Bcal^{s})} + \eps \Vert u_\varphi^\eps \Vert_{L^\infty(\mathcal{B}^{\f12})}^{\f12} \norm{\dd_z u_\varphi}_{\tilde{L}^2_t(\Bcal^{s+\frac{3}2})} \notag
	\\
	&\qquad \qquad + C\norm{(w^1_{\varphi},\eps w^2_\varphi,\theta_{\varphi})}_{\tilde{L}^2_{t,\dot{\eta}(t)}(\mathcal{B}^{s+ \frac{1}2})} + C\norm{T_{\varphi}}_{L^\infty_{t}(\mathcal{B}^{\frac{3}{2}})}^{\f12} \norm{\pa_z w^1_{\varphi}}_{\tilde{L}^2_{t}(\mathcal{B}^s)} + C\norm{u_{\varphi}}_{L^\infty_{t}(\mathcal{B}^{\frac{3}{2}})}^{\f32} \norm{\pa_z \theta_{\varphi}}_{\tilde{L}^2_{t}(\mathcal{B}^s)} \Big).  \notag
\end{align}
Now, from the results of the previous sections, we deduce the existence of a constant $M \geq 1$ such that
\begin{align*}
	\norm{u^\eps_\varphi}_{L^\infty_t(\mathcal{B}^{\frac{1}{2}})} + \norm{u^\eps_\varphi}_{L^\infty_t(\mathcal{B}^{\frac{3}{2}})}
	+ \Vert (\pa_t u)_\varphi \Vert_{\tilde{L}^2_t(\Bcal^{s + \frac{3}2})} + \Vert \pa_z u_\varphi \Vert_{\tilde{L}^2_t(\Bcal^{s + \frac{3}2})} + \Vert \pa_z u_\varphi \Vert_{\tilde{L}^2_t(\Bcal^{s + \frac{5}2})} \leq M
\end{align*}
and we deduce also that, for $\norm{(u_0,T_0)}_{\Bcal^{\frac{3}2}}$ small enough,
\begin{equation*}
	\norm{T_\varphi}_{L^\infty_t(\mathcal{B}^{\frac{3}{2}})} + \norm{u_\varphi}_{L^\infty_t(\mathcal{B}^{\frac{3}{2}})} \leq \frac{1}{4C}.
\end{equation*}
Multiplying \eqref{eq:conv01} by $2^{qs}$ and then summing with respect to $q$, we obtain
\begin{align*}
	\label{eq:conv02}
	&\Vert (w^1_\varphi,\eps w^2_\varphi,\theta_\varphi) \Vert_{L^\infty_t(\Bcal^s)} + \sqrt{\mu} \norm{(w^1_\varphi,\eps w^2_\varphi,\theta_\varphi)}_{\tilde{L}^2_{t,\dot{\eta}(t)}(\Bcal^{s + \frac{1}2})}
	\\ 
	&\qquad + \frac{1}4 \Vert \pa_z (w^1_{\varphi},\eps w^2_{\varphi},\theta_\varphi) \Vert_{\tilde{L}^2_t(\Bcal^s)}
	+ \frac{\eps}2 \Vert (w^1_{\varphi},\eps w^2_{\varphi}) \Vert_{\tilde{L}^2_t(\Bcal^{s + \frac{1}2})} + \Vert \theta_{\varphi} \Vert_{\tilde{L}^2_t(\Bcal^{s + \frac{1}2})} \notag
	\\
	&\qquad \qquad \leq 2C \norm{(w^1_{\varphi},\eps w^2_{\varphi},\theta_\varphi)(0)}_{\Bcal^s} + C M \eps  + C\norm{(w^1_{\varphi},\eps w^2_\varphi,\theta_{\varphi})}_{\tilde{L}^2_{t,\dot{\eta}(t)}(\mathcal{B}^{s+ \frac{1}2})}.  \notag
\end{align*}

Then by taking $\mu \geq C^2$, we can complete the proof Theorem \ref{th:Convergence}. \hfill $\square$

}

\end{document}